\documentclass[preprint,12pt,numbers,sort&compress]{elsarticle}
\usepackage{soul}
\usepackage{psfrag}
\usepackage{a4wide}
\usepackage{rotating}
\usepackage{color}
\usepackage{amssymb}
\usepackage{amsmath}
\usepackage{amsthm}
\usepackage{verbatim}

\newcommand{\f}[1]{\mathbf{#1}}
\newcommand{\ab}[1]{\boldsymbol{#1}}
\def\bfm#1{\boldsymbol{#1}}
\newcommand{\bb}[1]{\bfm{#1}}

\newcommand{\R}{\mathbb R}

\newcommand{\V}{\mathcal{V}}

\newcommand{\W}{\mathcal{W}}
\newcommand{\glob}{y}
\newcommand{\loc}{\zeta}

\DeclareMathOperator{\Span}{span}

\theoremstyle{definition}
\newtheorem{ex}{Example}

\newproof{pf}{proof}
\advance\textheight by 0.4cm
\advance\topmargin by -0.2cm

\definecolor{gold}{rgb}{1,0.7,0}
\definecolor{dred}{rgb}{0.92,0,0}
\definecolor{dgreen}{rgb}{0,0.6,0}

\begin{document}

\begin{frontmatter}

\title{Isogeometric collocation on planar multi-patch domains}

\author[lnz]{Mario Kapl\corref{cor}}
\ead{mario.kapl@ricam.oeaw.ac.at}

\author[slo1]{Vito Vitrih}
\ead{vito.vitrih@upr.si}
 
\address[lnz]{Johann Radon Institute for Computational and Applied Mathematics, \\Austrian Academy of Sciences, Linz, Austria}

\address[slo1]{IAM and FAMNIT, University of Primorska, Koper, Slovenia}

\cortext[cor]{Corresponding author}

\begin{abstract}
We present an isogeometric framework based on collocation to construct a $C^2$-smooth approximation of the solution of the Poisson's equation over planar bilinearly parameterized 
multi-patch domains. The construction of the used globally $C^2$-smooth discretization space for the partial differential equation is simple and works uniformly for all possible 
multi-patch configurations.
The basis of the $C^2$-smooth space can be described as the span of three different types of 
locally supported functions corresponding to the single patches, edges and vertices of the multi-patch domain. For the selection of the collocation points, which is important for the 
stability and convergence of the collocation problem, two different choices are numerically investigated.
The first approach employs the tensor-product Greville abscissae as collocation points, and shows for the multi-patch case the same convergence behavior as for the one-patch 
case~\cite{IsoCollocMethods2010}, which is suboptimal in particular for odd spline degree. The second approach generalizes the concept of superconvergent points 
from the one-patch case (cf.~\cite{SuperConvergent2015,GomezLorenzisVariationalCollocation,MonSanTam2017})
to the multi-patch case. Again, these points possess better convergence properties than Greville 
abscissae in case of odd spline degree. 
\end{abstract}


\end{frontmatter}

\section{Introduction}

Isogeometric Analysis (IgA) is an approach for numerically solving a partial differential equation (PDE) by using the same spline or NURBS space for describing 
the geometry and for representing the solution of the considered PDE, cf.~\cite{ANU:9260759,CottrellBook,HuCoBa04}. Usually, the developed PDE solvers are based on finding 
a solution of the variational (weak) form of the PDE, since using the weak form requires less smooth functions compared to directly employing the 
strong form. However, this technique needs the evaluation of integrals by using some particular quadrature 
rules, and the accuracy of the solution depends on the quality of the numerical integration. In contrast, solving the strong form of the PDE via collocation eliminates 
integration but requires spaces of higher regularity.

While in case of a one-patch domain, the higher continuity of the functions can be
easily guaranteed by using spline functions with the desired smoothness within the patch, the construction of $C^s$-smooth ($s \geq 1$) isogeometric spline spaces over 
multi-patch domains is challenging, and is the task of current research, see e.g. \cite{BeMa14,BlMoVi17,CoSaTa16,KaBuBeJu16,KaSaTa17a,KaSaTa17c,KaSaTa19b,KaViJu15,Pe15-2,
mourrain2015geometrically,NgPe16,ToSpHu17b,ToSpHu17,ChAnRa18} for $s=1$ and \cite{KaVi17a,KaVi17b,KaVi17c,KaVi19a,ToSpHiHu16} for $s=2$.
The design of $C^s$-smooth multi-patch spline spaces is linked to the concept of geometric continuity of multi-patch surfaces, 
cf. \cite{HoLa93,Pe02}, due to the fact that \emph{an isogeometric function is $C^{s}$-smooth on a multi-patch domain if and only if its associated multi-patch graph 
surface is $G^{s}$-smooth}, cf. \cite{Pe15,KaViJu15}.

So far, the problem of isogeometric collocation has been mostly explored on one-patch domains, see e.g. \cite{SuperConvergent2015,IsoCollocMethods2010,
ShearTimoshenko2012,GomezLorenzisVariationalCollocation,MonSanTam2017,CostComparison2013,Reali2015,EnzoKiendlLorenzis2019,EnzoKiendlLorenzis2017,EvansHughesReali2018,
FahrenderLorenzisGomez2018}.  
For analyzing the convergence behavior of a particular approach, the errors are generally computed with respect to $L^{\infty}$, $W^{1,\infty}$ and $W^{2,\infty}$~norm, 
or equivalently with respect to $L^2$, $H^1$ and $H^2$~norm, respectively. The study of isogeometric collocation methods has started in~\cite{IsoCollocMethods2010} by 
using the Greville and Demko abscissae as collocation points. In comparison to the Galerkin approach, both choices show a suboptimal convergence behavior with respect 
to $L^2 (L^{\infty})$ norm, namely of orders~$\mathcal{O}(h^{p})$ and $\mathcal{O}(h^{p-1})$ under $h$-refinement for even and odd spline degree~$p$, respectively, and 
additionally a suboptimal convergence order~$\mathcal{O}(h^{p-1})$ with respect to $H^1$ ($W^{1,\infty}$) norm but just in case of odd spline degree~$p$. 

In~\cite{SuperConvergent2015}, an isogeometric collocation method has been presented, which is based on the computation and on the use of specific collocation points 
called  
\emph{superconvergent points}. The proposed technique possesses only for even spline degree~$p$ in case of the $L^2$ ($L^{\infty}$) 
norm a suboptimal convergence behavior of order~$\mathcal{O}(h^{p})$, and is optimal for all other cases. This is achieved at the expense of solving an overdetermined linear 
system due to the fact that the number of superconvergent points is larger than the number of degrees of freedom. The two methods~\cite{GomezLorenzisVariationalCollocation} and 
\cite{MonSanTam2017} select different specific subsets of the superconvergent points to achieve that the number of collocation points coincides with the 
number of degrees of freedom. While the points\footnote{In~\cite{GomezLorenzisVariationalCollocation}, the superconvergent points are also denoted as 
\emph{Cauchy Galerkin points}, since collocation at these points produces the Galerkin solution exactly.} in~\cite{GomezLorenzisVariationalCollocation} are chosen in an alternating 
way, which reduces for the case of odd spline degree~$p$ the convergence order also to $\mathcal{O}(h^{p})$ in the $L^2$ ($L^{\infty}$) norm, the points in \cite{MonSanTam2017} are 
selected in a clustered way to maintain the convergence behavior of~\cite{SuperConvergent2015}. 

In~\cite{SuperConvergent2015,IsoCollocMethods2010,FahrenderLorenzisGomez2018,MonSanTam2017,CostComparison2013}, the concept of isogeometric collocation has been mostly 
used for solving the Poisson's equation. Further problems and PDEs of interest in the framework of isogeometric collocation are amongst others linear and nonlinear 
elasticity, elastostatistics and elastodynamics~
\cite{SuperConvergent2015,Elatostatics2012,EvansHughesReali2018,FahrenderLorenzisGomez2018,GomezLorenzisVariationalCollocation,JiaAnitescuZhangRabczuk2019,Reali2015,CostComparison2013}, 
plate/shell problems~\cite{GomezLorenzisVariationalCollocation,EnzoKiendlLorenzis2017} and beam problems~\cite{ShearTimoshenko2012,EnzoKiendlLorenzis2019}. 

As already mentioned before, most of the existing isogeometric collocation methods have in common that they are mainly restricted to the case of one-patch domains. 
Two techniques which also deal with the case of multi-patch domains are~\cite{Elatostatics2012,JiaAnitescuZhangRabczuk2019}. While in~\cite{Elatostatics2012} standard 
NURBS functions are employed, the work in~\cite{JiaAnitescuZhangRabczuk2019} is based on the use of PHT-splines~\cite{DeChLi08}. Since for both methods the used spline 
spaces are in general just $C^0$-smooth across the patch interfaces, special techniques for the collocation on the interfaces are required. 

The goal of this paper is to compute a globally $C^2$-smooth approximation of the solution of the Poisson's equation on planar bilinearly parameterized multi-patch 
domains by means of isogeometric collocation. For this purpose, a globally $C^2$-smooth discretization space is generated and used. 
The construction of this $C^2$-smooth multi-patch space is based on and extends the work~\cite{KaVi19a}.
There, a 
subspace of the entire~$C^2$-smooth spline space~\cite{KaVi17b} has been generated, which possesses a simpler structure as the entire $C^2$-smooth space but still maintains the full 
approximation properties. The construction of the locally supported basis functions is simple via explicit formulae or by computing the null space of a small system of 
linear equations. 
Whereas the basis construction in~\cite{KaVi19a} 
works uniformly for all possible multi-patch configurations, the number of basis functions depends on the initial geometry. This work overcomes the 
latter limitation by constructing a $C^2$-smooth spline space whose dimension is independent of the initial geometry by e.g. additionally enforcing that the functions 
have to be $C^4$-smooth at the inner vertices and at the boundary vertices of patch valencies greater or equal to three of the multi-patch domain.

The generated $C^2$-smooth spline space can be directly employed to perform isogeometric collocation on the multi-patch domain, and no special treatment for the patch 
interfaces as in \cite{Elatostatics2012,JiaAnitescuZhangRabczuk2019} is needed. For assembling the system matrices of the collocation problem two different choices of 
collocation points are studied. On the one hand, we employ the tensor-product Greville abscissae for the single patches, and on the other hand we generalize the so-called 
superconvergent points for the one-patch case~(cf. \cite{SuperConvergent2015,GomezLorenzisVariationalCollocation,MonSanTam2017}) to the multi-patch case. 
Like in the one-patch case, the 
superconvergent points possess a better convergence behavior as the Greville points for odd spline degree, namely of orders $\mathcal{O}(h^p)$, $\mathcal{O}(h^p)$ and 
$\mathcal{O}(h^{p-1})$ in comparison to $\mathcal{O}(h^{p-1})$, $\mathcal{O}(h^{p-1})$ and $\mathcal{O}(h^{p-1})$ with respect to the $L^2$ 
$(L^{\infty})$, $H^1$ ($W^{1,\infty}$) and $H^2$ ($W^{2,\infty}$) norm.

The remainder of the paper is organized as follows. Section~\ref{C2_isogeometric_splinespace} presents the construction of a particular $C^2$-smooth isogeometric spline space over 
bilinear multi-patch domains, which will be used as a discretization space for solving the Poisson's equation by means of isogeometric collocation. The multi-patch isogeometric 
collocation method as well as two different strategies for choosing the collocation points are described in Section~\ref{sec:collocation}. On the basis of several numerical examples, 
the convergence behavior under $h$-refinement of the collocation method is studied in Section~\ref{section_Numerical_examples}, and demonstrates the potential of the method for the use 
in IgA. Finally, we conclude the paper.

\section{A $C^2$-smooth isogeometric spline space} \label{C2_isogeometric_splinespace}

The construction of a particular $C^2$-smooth multi-patch isogeometric spline space will be presented, which will be based on and will extend the work~\cite{KaVi19a}. 
The $C^2$-smooth space will be used in Section~\ref{sec:collocation} to build an isogeometric collocation method for solving the Poisson's equation over bilinear multi-patch domains. 
We will start with the introduction of the multi-patch setting, which will be used throughout the paper.

\subsection{The multi-patch setting} \label{subsec:multipatch}
Let $\mathcal{I}_{\Omega}$ be an index set, and let $\Omega$ and $\Omega^{(i)}$, $i \in \mathcal{I}_{\Omega}$, be open domains in $\R^2$, such that 
$\overline{\Omega} = \cup_{i \in \mathcal{I}_{\Omega}} \overline{\Omega^{(i)}}$. Furthermore, let $\Omega^{(i)}$,
$i \in \mathcal{I}_{\Omega}$, be quadrangular patches, which are mutually disjoint, and the closures of any two of them have either an empty intersection, possess 
exactly one common vertex or share the whole common edge. Additionally, the deletion of any vertex of the multi-patch domain~$\overline{\Omega}$ does not split 
$\overline{\Omega}$ into subdomains, whose union would be unconnected. We will further assume that each patch $\overline{\Omega^{(i)}}$ is parameterized by a bilinear, 
bijective and regular geometry mapping~$\ab{F}^{(i)}$,
\begin{align*}
 \ab{F}^{(i)}: [0,1]^{2}  \rightarrow \R^{2}, \quad 
 \bb{\xi} =(\xi_1,\xi_2) \mapsto
 \ab{F}^{(i)}(\bb{\xi}) = \ab{F}^{(i)}(\xi_1,\xi_2), \quad i \in \mathcal{I}_{\Omega},
\end{align*}
such that $\overline{\Omega^{(i)}} = \ab{F}^{(i)}([0,1]^{2})$.
We will also use the splitting of the multi-patch domain $\overline{\Omega}$ into the single patches~$\Omega^{(i)}$, $i \in \mathcal{I}_{\Omega}$, edges~$\Gamma^{(i)}$, 
$i \in \mathcal{I}_{\Gamma}$ and vertices~$\ab{v}^{(i)}$, $i \in \mathcal{I}_{\Xi}$, i.e.
\begin{equation*} 
\displaystyle
\overline{\Omega} = \bigcup_{i \in \mathcal{I}_{\Omega}} \Omega^{(i)}  \; \dot{\cup}  \bigcup_{i \in \mathcal{I}_{\Gamma}} \Gamma^{(i)} \; \dot{\cup} \bigcup_{i \in \mathcal{I}_{\Xi}} \bfm{v}^{(i)},
\end{equation*}
where $\mathcal{I}_{\Gamma}$ and $\mathcal{I}_{\Xi}$ are the index sets of the indices of the edges~$\Gamma^{(i)}$ and vertices~$\ab{v}^{(i)}$, respectively. We further 
have $\mathcal{I}_{\Gamma} = \mathcal{I}_{\Gamma_I} \dot{\cup} \mathcal{I}_{\Gamma_B}$ and $\mathcal{I}_{\Xi} =  \mathcal{I}_{\Xi_I} \dot{\cup} \mathcal{I}_{\Xi_B}$ 
with $\mathcal{I}_{\Xi_B} =   \mathcal{I}_{\Xi_{1}} \dot{\cup} \mathcal{I}_{\Xi_{2}} \dot{\cup} \mathcal{I}_{\Xi_{3}}$, where $\mathcal{I}_{\Gamma_I}$ and $\mathcal{I}_{\Gamma_B}$ 
contain all indices of interfaces and boundary edges, respectively, $\mathcal{I}_{\Xi_I}$ and $\mathcal{I}_{\Xi_B}$ collect all indices of inner and boundary vertices, 
respectively, and $\mathcal{I}_{\Xi_{1}}$, $\mathcal{I}_{\Xi_{2}}$ and $\mathcal{I}_{\Xi_{3}}$ are further the index sets of the indices of boundary vertices of patch 
valencies~one, two and greater or equal to three, respectively. In addition, the patch valency of a vertex~$\ab{v}^{(i)}$ will be denoted by $\nu_i$.

Let $\mathcal{S}_h^{p,r}([0,1])$ be the univariate spline space of degree~$p$, regularity~$r$ and mesh size~$h=\frac{1}{k+1}$ on the unit interval~$[0,1]$, constructed from 
the uniform open knot vector
\begin{equation*}  
(t_0^{p,r},\ldots,t_{2p+k(p-r)+1}^{p,r})=(\underbrace{0,\ldots,0}_{(p+1)-\mbox{\scriptsize times}},
\underbrace{\textstyle \frac{1}{k+1},\ldots,\frac{1}{k+1}}_{(p-r) - \mbox{\scriptsize times}},\ldots, 
\underbrace{\textstyle \frac{k}{k+1},\ldots ,\frac{k}{k+1}}_{(p-r) - \mbox{\scriptsize times}},
\underbrace{1,\ldots,1}_{(p+1)-\mbox{\scriptsize times}}),
\end{equation*}
where $k$ is the 
number of different inner knots. Furthermore, let $\mathcal{S}_h^{\ab{p},\ab{r}}([0,1]^2)$ be the tensor-product spline space $\mathcal{S}_h^{p,r}([0,1]) \otimes 
\mathcal{S}_h^{p,r}([0,1])$ on the unit-square~$[0,1]^2$. We denote the B-splines of the spaces~$\mathcal{S}_h^{p,r}([0,1])$ and $\mathcal{S}_h^{\ab{p},\ab{r}}([0,1]^2)$ by
$N_{j}^{p,r}$ and $N_{j_1,j_2}^{\ab{p},\ab{r}}=N_{j_1}^{p,r}N_{j_2}^{p,r}$, respectively, with $j,j_1,j_2=0,1,\ldots,n-1$, and $n=p+1+k(p-r)$.
We assume that $p \geq 5$, $2 \leq r \leq p-3$, and that the number of inner knots satisfies $k \geq \frac{9-p}{p-r-2}$. 
These assumptions are necessary to ensure that 
the constructed $C^2$-smooth spline spaces in Section~\ref{subsec:C2space} will be $h$-refineable and well-defined, see also~\cite{KaVi17c,KaVi19a}.
Since the geometry mappings~$\ab{F}^{(i)}$, $i \in \mathcal{I}_{\Omega}$, are bilinearly parameterized, we trivially have that
\[
\ab{F}^{(i)} \in \mathcal{S}_{h}^{\ab{p},\ab{r}}([0,1]^2) \times \mathcal{S}^{\ab{p},\ab{r}}_{h}([0,1]^2).
\]

\begin{figure}[htp] 
\centering
\includegraphics[width=4.9cm]{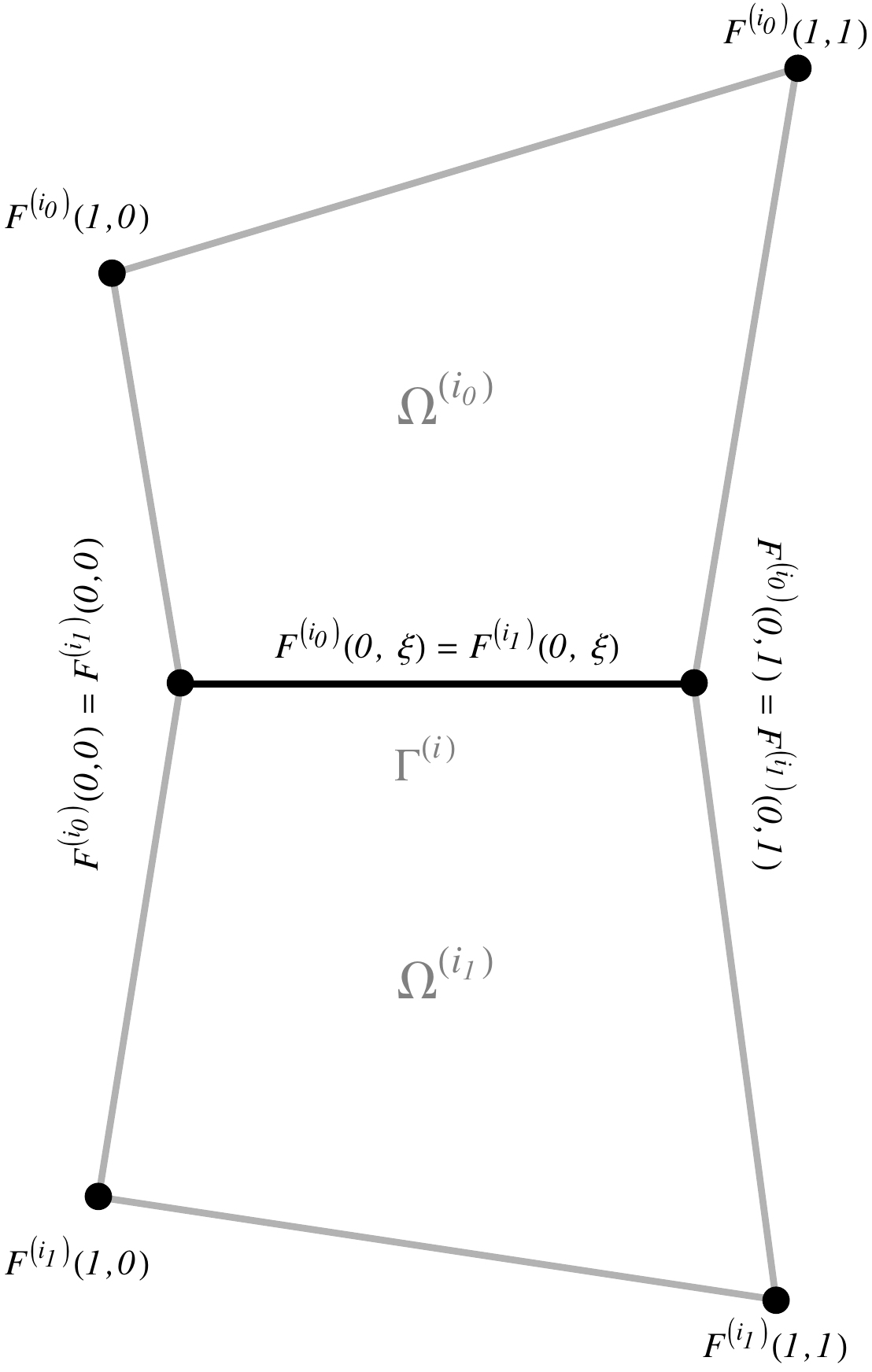}
\hskip1em
\includegraphics[width=9.5cm]{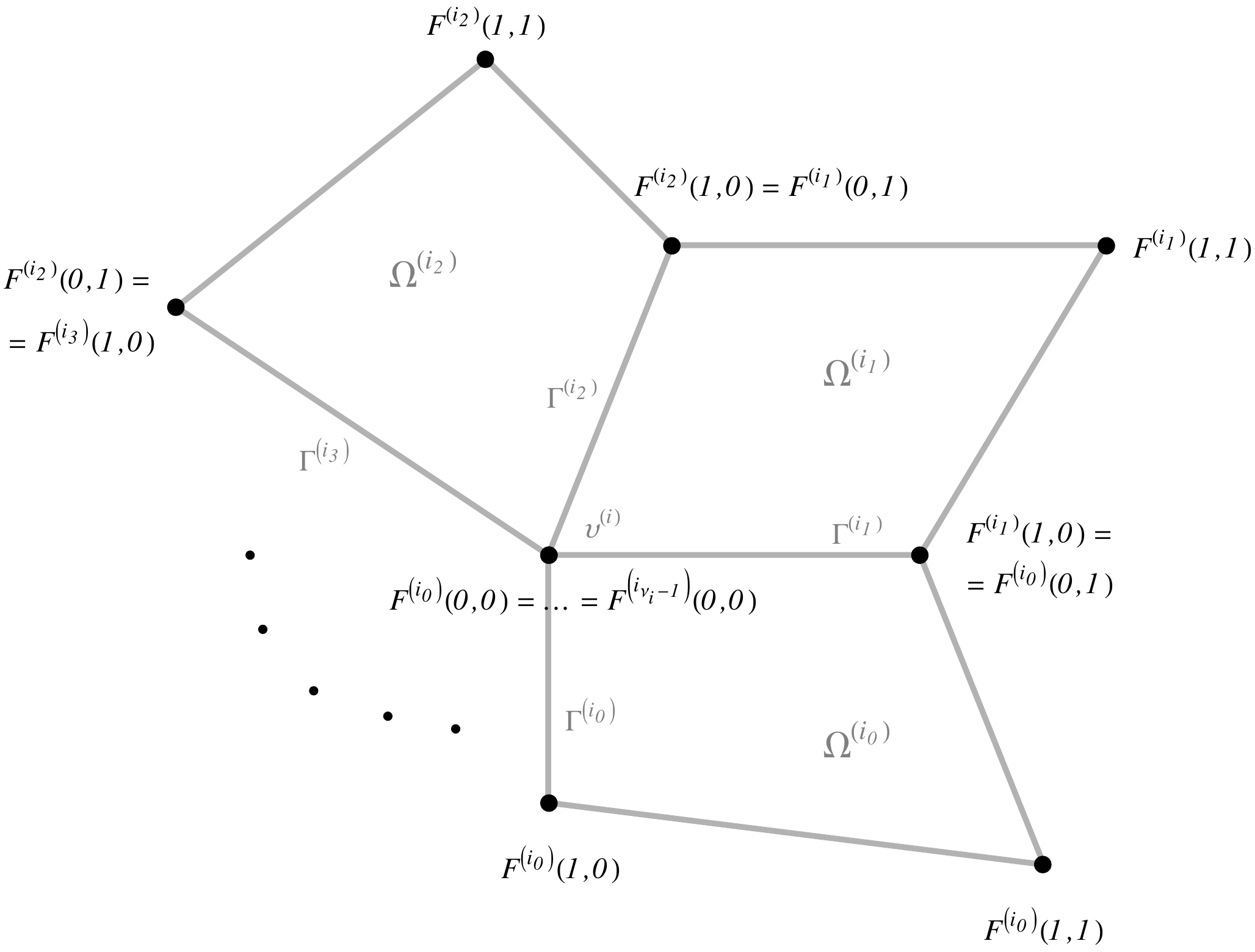}
\caption{Left: Considering two neighboring patches~$\Omega^{(i_0)}$ and $\Omega^{(i_1)}$, we can always assume that the two corresponding geometry mappings 
$\ab{F}^{(i_0)}$ and $\ab{F}^{(i_1)}$ are parameterized as shown. Right: For the patches $\Omega^{(i_\ell)}$, $\ell =0, 1,\ldots, \nu_{i} -1$, which enclose the vertex~$\bfm{v}^{(i)}$, 
the corresponding geometry mappings $\ab{F}^{(i_\ell)}$ can be always reparameterized to be given as visualized.} 
\label{fig:situation_two_patches}
\end{figure}

Moreover, we assume for two neighboring patches~$\Omega^{(i_0)}$ and $\Omega^{(i_1)}$ with the common interface~$\Gamma^{(i)} \subset \overline{\Omega^{(i_0)}} \cap 
\overline{\Omega^{(i_1)}}$ that the associated geometry mappings~$\ab{F}^{(i_0)}$ and $\ab{F}^{(i_1)}$ are parameterized as given in Fig.~\ref{fig:situation_two_patches}~(left), 
and assume for the neighboring patches enclosing an inner or boundary vertex~$\ab{v}^{(i)}$ of patch valency~$\nu_i \geq 3$ that the associated geometry mappings are parameterized 
as shown in Fig.~\ref{fig:situation_two_patches}~(right). For the latter case, we additionally assume that the patches enclosing the vertex~$\ab{v}^{(i)}$ are 
relabeled by~$\Omega^{(i_0)}, \ldots,\Omega^{(i_{\nu_i -1})}$, and that the common interfaces\footnote{In case of a 
boundary vertex~$\ab{v}^{(i)}$, we denote the boundary edges by $\Gamma^{(i_0)}$ and $\Gamma^{(i_{\nu_i})}$.} of the two-patch subdomains $\overline{\Omega^{(i_\ell)}} 
\cup \overline{\Omega^{(i_{\ell+1})}}$, $\ell=0,1,\ldots,\nu_i -1$, are relabeled by $\Gamma^{(i_{\ell+1})}$. In case of an inner vertex~$\ab{v}^{(i)}$ we 
further consider the lower index~$\ell$ of the indices $i_{\ell}$ modulo~$\nu_i$,  which just means that the lower index~$\ell$ is replaced by the remainder of the division of 
$\ell$ by $\nu_i$. Note that the local reparameterizations as shown in Fig.~\ref{fig:situation_two_patches} are always possible (if necessary), and hence can be assumed to 
be satisfied without loss of generality. 

\subsection{$C^2$-smoothness across patch interfaces}

The space of $C^2$-smooth isogeometric functions on~$\Omega$ is given as
\begin{equation*}
\V = \left\{ \phi \in C^{2}(\overline{\Omega}) : \; \phi |_{\overline{\Omega}^{(i)}} \in {\mathcal{S}_{h}^{\ab{p},\ab{r}}([0,1]^{2})} \circ (\ab{F}^{(i)})^{-1}, \; 
i \in \mathcal{I}_{\Omega}   \right\}.
\end{equation*}
This space can be characterized by means of the concept of geometric continuity of multi-patch surfaces, cf.~\cite{HoLa93,Pe02}. \emph{An isogeometric 
function~$\phi$ belongs to the space~$\mathcal{V}$ if and only if for any two neighboring patches~$\Omega^{(i_0)}$ and $\Omega^{(i_1)}$ with the common 
interface~$\overline{\Gamma^{(i)}} = \overline{\Omega^{(i_0)}} \cap 
\overline{\Omega^{(i_1)}}$, $i \in \mathcal{I}_{\Gamma_{I}}$, the associated graph surface patches
$\left[\begin{array}{c}
  \ab{F}^{(i_0)} \\
  \phi \circ \ab{F}^{(i_0)}
 \end{array} \right]$
and 
$\left[\begin{array}{c}
  \ab{F}^{(i_1)} \\
  \phi \circ \ab{F}^{(i_1)}
 \end{array} \right]$
are $G^2$-smooth across their common interface}, see e.g.~\cite{Pe15,KaViJu15}. An equivalent condition to the $G^2$-smoothness of the two neighboring graph surface patches is that 
the two associated spline functions~$g^{(i_0)}=\phi \circ \ab{F}^{(i_0)}$ and $g^{(i_1)}=\phi \circ \ab{F}^{(i_1)}$ satisfy
\begin{equation} \label{eq:g0}
 \underbrace{g^{(i_0)}(0,\xi) = g^{(i_1)}(0,\xi)}_{:=g_0(\xi)},
\end{equation}
\begin{equation} \label{eq:g1}
\underbrace{\frac{D_{\xi_2}g^{(i_0)}(0,\xi)- \beta^{(i_0)}_{\Gamma^{(i)}}(\xi) g'_0(\xi)}{\alpha^{(i_0)}_{\Gamma^{(i)}}(\xi)} = 
\frac{D_{\xi_2}g^{(i_1)}(0,\xi)- \beta^{(i_1)}_{\Gamma^{(i)}}(\xi) g'_0(\xi)}{\alpha^{(i_1)}_{\Gamma^{(i)}}(\xi)}}_{:=g_1(\xi)},
\end{equation}
and
\begin{eqnarray}
 \frac{D_{\xi_1 \xi_1}g^{(i_0)}(0,\xi) - (\beta^{(i_0)}_{\Gamma^{(i)}}(\xi))^2 g''_0(\xi) - 2 \alpha^{(i_0)}_{\Gamma^{(i)}} \beta^{(i_0)}_{\Gamma^{(i)}}(\xi) g'_1(\xi) 
 }{(\alpha^{(i_0)}(\xi))^{2}} = \nonumber \\[-0.3cm] 
 \label{eq:g2} \\[-0.3cm] 
 \frac{D_{\xi_1 \xi_1}g^{(i_1)}(0,\xi) - (\beta^{(i_1)}_{\Gamma^{(i)}}(\xi))^2 g''_0(\xi) - 2 \alpha^{(i_1)}_{\Gamma^{(i)}} \beta^{(i_1)}_{\Gamma^{(i)}}(\xi) g'_1(\xi) 
 }{(\alpha^{(i_1)}(\xi))^{2}}, \nonumber
\end{eqnarray}
for all~$\xi \in [0,1]$, cf.~\cite[Lemma~1]{KaVi17c}, where $\alpha^{(\tau)}_{\Gamma^{(i)}}$ and $\beta^{(\tau)}_{\Gamma^{(i)}}$, $\tau \in \{i_0,i_1 \}$, are linear polynomials 
given by
\[
 \alpha_{\Gamma^{(i)}}^{(\tau)}(\xi) = c_1 \det [D_{\xi_1}\ab{F}^{(\tau)}(0,\xi), D_{\xi_2}\ab{F}^{(\tau)}(0,\xi)], 
\]
and
\[
 \beta_{\Gamma^{(i)}}^{(\tau)}(\xi) = \frac{D_{\xi_1} \ab{F}^{(\tau)}(0,\xi) \cdot D_{\xi_2}\ab{F}^{(\tau)}(0,\xi)}{||D_{\xi_2}\ab{F}^{(\tau)}(0,\xi)||^{2}},
\]
respectively, with $c_1 \in \R$ such that
\[
 || \alpha_{\Gamma^{(i)}}^{(i_0)}+1 ||^2_{L^2} + || \alpha_{\Gamma^{(i)}}^{(i_1)}-1 ||^2_{L^2}
\]
is minimized, cf. \cite{KaVi19a}. Note that $\alpha^{(i_0)}_{\Gamma^{(i)}} < 0$ and $\alpha^{(i_1)}_{\Gamma^{(i)}} > 0$, since the geometry mappings~$\ab{F}^{(i_0)}$ and 
$\ab{F}^{(i_1)}$ are regular.

The space~$\mathcal{V}$ has been studied in different configurations in~\cite{KaVi17a, KaVi17b,KaVi17c}. There, it has been observed that the space~$\mathcal{V}$ possesses 
a complex structure and its dimension depends on the initial geometry. In the next subsection, we will consider instead a subspace of the space~$\mathcal{V}$, 
whose structure is simpler and whose dimension is independent of the initial geometry. 

\subsection{A $C^2$-smooth isogeometric space} \label{subsec:C2space}

We will describe the construction of a particular $C^2$-smooth space~$\W \subseteq \mathcal{V}$, which will be similar to the one in~\cite{KaVi19a}. There will be 
two main differences between the newly constructed space~$\mathcal{W}$ and the space from~\cite{KaVi19a}. Firstly, the space~$\W$ will not be anymore restricted to 
containing just $C^2$-smooth functions fulfilling the homogeneous boundary conditions of order~$2$ given by
\begin{equation*}
 u(\ab{x}) = \frac{\partial u}{\partial \ab{n}}(\ab{x}) = \triangle u(\ab{x})  = 0 , \quad \ab{x} \in \partial \Omega,
 \end{equation*}
which have been used in~\cite{KaVi19a} as boundary conditions for solving the triharmonic equation via its weak form and standard Galerkin discretization. 
Secondly, the dimension of the space~$\W$ will now be independent of the initial geometry. 

Analogous to~\cite{KaVi19a}, the space~$\mathcal{W}$ is defined as the direct sum
\begin{equation} \label{eq:Wh}
    \W =  \left(\bigoplus_{i \in \mathcal{I}_{\Omega} } \W_{\Omega^{(i)}}\right) \oplus 
    \left(\bigoplus_{i \in \mathcal{I}_{\Gamma}} \W_{\Gamma^{(i)}}\right) 
    \oplus \left(\bigoplus_{i \in \mathcal{I}_{\Xi}} \W_{\bfm{v}^{(i)}}\right),
\end{equation}
where $\W_{\Omega^{(i)}}$, $\W_{\Gamma^{(i)}}$ and $\W_{\bfm{v}^{(i)}}$ are subspaces corresponding to the single patches~$\Omega^{(i)}$, edges~$\Gamma^{(i)}$ and 
vertices~$\ab{v}^{(i)}$, respectively. While the construction of the subspaces~$\W_{\Omega^{(i)}}$ and $\W_{\Gamma^{(i)}}$ will work analogously to or similar as 
in~\cite{KaVi19a}, and hence will be kept short below, the construction of the subspaces~$ \W_{\bfm{v}^{(i)}}$ will differ. 

\subsubsection{The patch subspace~$\W_{\Omega^{(i)}}$ and the edge subspace~$\W_{\Gamma^{(i)}}$}

As in~\cite{KaVi19a}, the patch subspaces~$\W_{\Omega^{(i)}}$, $i \in \mathcal{I}_{\Omega}$, are given as
\begin{equation*}
\mathcal{W}_{\Omega^{(i)}} =  \Span \{ \phi_{\Omega^{(i)}; j_1,j_2} |\;  j_1,j_2=3,4,\ldots,n-4  \},
\end{equation*}
with the functions
\begin{equation} \label{eq:defphiOmega}
 \phi_{\Omega^{(i)};j_1,j_2} (\bfm{x}) = 
    \begin{cases}
   (N_{j_1,j_2}^{\ab{p},\ab{r}}\circ (\ab{F}^{(i)})^{-1})(\bfm{x}) \;
   \mbox{ if }\f \, \bfm{x} \in\Omega^{(i)},
    \\ 0 \quad \mbox{ otherwise.}
    \end{cases}
\end{equation}

For the edge subspaces~$\W_{\Gamma^{(i)}}$, $i \in \mathcal{I}_{\Gamma}$, we have to distinguish between the case of a boundary edge~$\Gamma^{(i)}$, i.e. 
$i \in \mathcal{I}_{\Gamma_{B}}$, and the case of an interface~$\Gamma^{(i)}$, $i \in \mathcal{I}_{\Gamma_{I}}$. Let us start with the case of a boundary edge~$\Gamma^{(i)}$, 
$i \in \mathcal{I}_{\Gamma_{B}}$. Assume without loss of generality, that the boundary edge~$\Gamma^{(i)}$ is contained in the geometry mapping $\ab{F}^{(i_0)}$ and is 
given by $\ab{F}^{(i_0)}(\{0 \} \times (0,1))$. Then the subspace~$\W_{\Gamma^{(i)}}$ possesses the form
\begin{align*} 
\mathcal{W}_{\Gamma^{(i)}} =&  \Span \{ \phi_{\Gamma^{(i)}_B; j_1,j_2} |\; j_2=5-j_1,\ldots,n+j_1-6, \; j_1=0,1,2 \} ,
\end{align*}
where $\phi_{\Gamma^{(i)}_B;j_1,j_2}$ is defined similar to~\eqref{eq:defphiOmega} as
\begin{equation} \label{eq:defphiGammaBoundary}
 \phi_{\Gamma^{(i)}_B;j_1,j_2} (\bfm{x}) = 
    \begin{cases}
   (N_{j_1,j_2}^{\ab{p},\ab{r}}\circ (\ab{F}^{(i_0)})^{-1})(\bfm{x}) \;
   \mbox{ if }\f \, \bfm{x} \in \overline{\Omega^{(i_0)}},
    \\ 0 \quad \mbox{ otherwise. }
    \end{cases}
\end{equation}
In case of an interface~$\Gamma^{(i)} \subset \overline{\Omega^{(i_0)}} \cap \overline{\Omega^{(i_1)}}$, $i \in \mathcal{I}_{\Gamma_{I}}$, and assuming without loss of generality 
that the associated geometry mappings~$\ab{F}^{(i_0)}$ and $\ab{F}^{(i_1)}$ are parameterized as in~Fig.~\ref{fig:situation_two_patches}~(left), the edge 
subspace~$\mathcal{W}_{\Gamma^{(i)}}$ is defined as
\begin{equation*} 
{\mathcal{W}}_{\Gamma^{(i)}} = \Span \{ {\phi}_{\Gamma^{(i)}_I;j_1,j_2} | \;  j_2=5-j_1,\ldots,n_{j_1}+j_1-6, \; j_1=0,1,2\}, 
\end{equation*}
with the functions 
\begin{equation}  \label{eq:basisFunctionsGenericEdge}
{\phi}_{\Gamma^{(i)}_I;j_1,j_2}(\bfm{x})  = 
\begin{cases}
  (g_{\Gamma^{(i)}; j_1,j_2}^{(i_0)} \circ (\ab{F}^{(i_0)})^{-1})(\bfm{x}) \;
\mbox{ if }\f \, \bfm{x} \in \overline{\Omega^{(i_0)}},
\\[0.15cm] 
  (g_{\Gamma^{(i)}; j_1,j_2}^{(i_1)} \circ (\ab{F}^{(i_1)})^{-1} )(\bfm{x}) \;
\mbox{ if }\f \, \bfm{x} \in \overline{\Omega^{(i_1)}},
\end{cases}
\end{equation}
where
\begin{align}  \label{eq:basisFunctionsGenericG} 
g_{\Gamma^{(i)}; 0,j_2}^{(\tau)} (\xi_1,\xi_2) & = 
N_{j_2}^{p,r+2}(\xi_2) M_0(\xi_1) + \beta_{\Gamma^{(i)}}^{(\tau)}(\xi_2) (N_{j_2}^{p,r+2})'(\xi_2)  M_1(\xi_1)   \nonumber \\
& \; + \left( \beta^{(\tau)}_{\Gamma^{(i)}}(\xi_2)\right)^2 (N_{j_2}^{p,r+2})''(\xi_2)  M_2(\xi_1),\nonumber \\  
g_{\Gamma^{(i)}; 1,j_2}^{(\tau)} (\xi_1,\xi_2) & = \displaystyle \frac{p}{h}  \left( {\alpha}^{(\tau)}_{\Gamma^{(i)}}(\xi_2) {N}_{j_2}^{p-1,r+1}(\xi_2)  
M_1(\xi_1)  \right.  \\
& \left. \; + 2\,  {\alpha}^{(\tau)}_{\Gamma^{(i)}} (\xi_2)\beta^{(\tau)}_{\Gamma^{(i)}}(\xi_2) ({N}_{j_2}^{p-1,r+1})'(\xi_2)  M_2(\xi_1)\right), \nonumber \\
 g_{\Gamma^{(i)}; 2,j_2}^{(\tau)}  (\xi_1,\xi_2) & =  \displaystyle  \frac{p(p-1)}{h^2} \left( {\alpha}^{(\tau)}_{\Gamma^{(i)}}(\xi_2)\right)^2 
{N}_{j_2}^{p-2,r}(\xi_2) M_2(\xi_1)   \nonumber,
\end{align}
for $\tau \in \{i_0,i_1\} $,
with
\begin{equation*} 
  M_0(\xi) = \sum_{j=0}^2 N_j^{p,r}(\xi) , \;
  M_1(\xi) = \frac{p}{h} \left( N_1^{p,r}(\xi) + 2 N_2^{p,r}(\xi)  \right), \;
  M_2(\xi) = \frac{p(p-1)}{h^2} N_2^{p,r}(\xi), 
\end{equation*}  
and 
\begin{equation*} 
 n_0 = \dim\left( \mathcal{S}_{h}^{p,r+2}([0,1])\right),\;
 n_1 = \dim \left( \mathcal{S}_{h}^{p-1,r+1}([0,1])\right),\;
 n_2 = \dim \left( \mathcal{S}_{h}^{p-2,r}([0,1])\right).
\end{equation*}

All functions of the patch subspaces~$\W_{\Omega^{(i)}}$, $i \in \mathcal{I}_{\Omega}$, and of the edge subspaces~$\W_{\Gamma^{(i)}}$, $i \in \mathcal{I}_{\Gamma}$, 
are $C^2$-smooth on~$\overline{\Omega}$, since the Eqs.~\eqref{eq:g0}--\eqref{eq:g2} are satisfied for all 
interfaces~$\overline{\Gamma^{(j)}}$, $j \in \mathcal{I}_{\Gamma_{I}}$. In case of the patch subspaces~$\W_{\Omega^{(i)}}$ and of the edge subspaces~$\W_{\Gamma^{(i)}}$ 
for $i \in \mathcal{I}_{\Gamma_{B}}$, the equations are even trivially satisfied, since the functions have vanishing values, gradients and Hessians on all 
interfaces~$\overline{\Gamma^{(j)}}$, $j \in \mathcal{I}_{\Gamma_{I}}$. In case of the edge subspaces~$\W_{\Gamma^{(i)}}$ for $i \in \mathcal{I}_{\Gamma_{I}}$, whereas 
the equations are also trivially fulfilled for the interfaces~$\overline{\Gamma^{(j)}}$, $j \in \mathcal{I}_{\Gamma_{I}} \setminus \{ i\}$, the particular construction of 
the functions~$\phi_{\Gamma^{(i)}_I;j_1,j_2}$ ensures that the Eqs.~\eqref{eq:g0}--\eqref{eq:g2} are satisfied for the non-trivial case of the 
corresponding interface~$\overline{\Gamma^{(i)}}$, too, cf.~\cite{KaVi17c,KaVi19a}.

\subsubsection{The vertex subspace~$\W_{\ab{v}^{(i)}}$}

The vertex subspaces~$\W_{\ab{v}^{(i)}}$, $i \in \mathcal{I}_{\Xi}$, constructed in this work, will differ from the ones in~\cite{KaVi19a}. The new construction of the 
subspaces will ensure that the dimension of the vertex subspaces will now be independent of the initial geometry. Furthermore, the design of the vertex 
subspaces~$\W_{\ab{v}^{(i)}}$, $i \in \mathcal{I}_{\Xi_B}$, will not be anymore restricted to any boundary conditions. We will distinguish between four cases, namely 
if the vertex~$\ab{v}^{(i)}$ is an inner vertex, i.e. $i \in \mathcal{I}_{\Xi_I}$, a boundary vertex of patch valency~one, i.e. $i \in \mathcal{I}_{\Xi_1}$, 
a boundary vertex of patch valency~two, i.e. $i \in \mathcal{I}_{\Xi_2}$, or a boundary vertex of patch valency greater or equal to three, i.e. $i \in \mathcal{I}_{\Xi_3}$. 

Let us start with the case of an inner vertex~$\ab{v}^{(i)}$, that is $i \in \mathcal{I}_{\Xi_I}$. Recall that we assume without loss of generality that the 
geometry mappings~$\ab{F}^{(i_{\ell})}$ containing the vertex~$\ab{v}^{(i)}$ are relabeled and parameterized 
as described and shown in Section~\ref{subsec:multipatch} and Fig.~\ref{fig:situation_two_patches}~(right), respectively. Again, the lower index~$\ell$ of the 
indices $i_{\ell}$ will be considered modulo~$\nu_i$. First of all, we consider analogous to~\cite{KaVi19a}, the isogeometric function 
$\phi_{\ab{v}^{(i)}}: \Omega \to \R $, which is defined as
\begin{equation}  \label{eq:defphiXi}
  \phi_{\ab{v}^{(i)}} (\ab{x}) = 
  \begin{cases}
   (g_{i_\ell} \circ (\ab{F}^{(i_\ell)})^{-1})(\bfm{x}) \;
\mbox{ if }\f \, \bfm{x} \in \overline{\Omega^{(i_\ell)}},\; \ell=0,1,\ldots,{\nu}_i -1,
\\
0 \quad \mbox{ otherwise,}
\end{cases}
\end{equation}
where the functions~$g_{i_{\ell}}$ are given as
\begin{align}
 g_{i_{\ell}}(\xi_1,\xi_2)   = & \underbrace{\sum_{j_1=0}^2 \sum_{j_2=0}^{4-j_1} a^{\Gamma^{(i_\ell)}}_{j_1,j_2} \, 
  g_{\Gamma^{(i_\ell); j_1,j_2}}^{(i_\ell)} (\xi_2,\xi_1)}_{:=g_{i_\ell}^{\Gamma^{(i_\ell)}}(\xi_1,\xi_2)}  
  +  \underbrace{\sum_{j_1=0}^2 \sum_{j_2=0}^{4-j_1} a^{\Gamma^{(i_{\ell+1})}}_{j_1,j_2} \, 
  g_{\Gamma^{(i_{\ell+1}); j_1,j_2}}^{(i_\ell)} (\xi_1,\xi_2) }_{:=g_{i_\ell}^{\Gamma^{(i_{\ell+1})}}(\xi_1,\xi_2)} \, - \nonumber \\[-0.3cm] 
 \label{eq:g_vertex} \\[-0.3cm]   
 & \underbrace{\sum_{j_1=0}^2 \sum_{j_2=0}^{2} a^{(i_\ell)}_{j_1,j_2} N_{j_1,j_2}^{\ab{p},\ab{r}} (\xi_1,\xi_2)}_{:=g_{i_\ell}^{\Omega^{(i_\ell)}}(\xi_1,\xi_2)}, 
 \nonumber
  \end{align}
with coefficients $a^{\Gamma^{(i_\ell)}}_{j_1,j_2}, a_{j_1,j_2}^{(i_\ell)} \in \R$, and with the 
functions~$g_{\Gamma^{(i_\ell); j_1,j_2}}^{(i_\ell)}$ and $g_{\Gamma^{(i_\ell+1); j_1,j_2}}^{(i_\ell)}$ given by~\eqref{eq:basisFunctionsGenericG}.

In \cite[Lemma~3]{KaVi19a}, it was shown, that the function~$\phi_{\ab{v}^{(i)}}$ is $C^2$-smooth on~$\overline{\Omega}$ if the 
coefficients~$a^{\Gamma^{(i_\ell)}}_{j_1,j_2}$ and $a_{j_1,j_2}^{(i_\ell)}$ are selected in such a way that the following 
homogeneous linear system
\begin{equation} \label{eq:vertex_homogeneous_system}
  \partial_{\xi_1}^i  
  \partial_{\xi_2}^j 
   \left( g_{i_\ell}^{\Gamma^{(i_{\ell+1})}} - g_{i_\ell}^{\Gamma^{(i_\ell)}} \right)  (\bfm{0}) = 0 \quad \mbox{ and } \quad
   \partial_{\xi_1^{}}^i  
   \partial_{\xi_2^{}}^j 
   \left( g_{i_\ell}^{\Gamma^{(i_{\ell+1})}} - g_{i_\ell}^{\Omega^{(i_\ell)}} \right)  (\bfm{0}) = 0, 
\end{equation}
with $ 0 \leq i,j \leq 2$ and $\ell =0,1, \ldots, \nu_{i}-1$, is satisfied. The construction of a $C^2$-smooth isogeometric function~$\phi_{\ab{v}^{(i)}}$ is based 
on the idea of generating the spline function~$g_{i_\ell}=\phi_{\ab{v}^{(i)}} \circ \ab{F}^{(i_{\ell})}$ for each patch~$\overline{\Omega^{(i_{\ell})}}$, $\ell=0,\ldots, 
\nu_{i}-1$, by using appropriate linear combinations of functions~$g_{\Gamma^{(i_{\ell})};j_1,j_2}^{(i_{\ell})}$ and 
$g_{\Gamma^{(i_{\ell+1})};j_1,j_2}^{(i_{\ell})}$, $j_1=0,1,2$, $j_2=0,1,\ldots,4-j_1$, to ensure $C^2$-smoothness across the interfaces~$\Gamma^{(i_\ell)}$ 
and $\Gamma^{(i_{\ell+1})}$, respectively, and by subtracting a linear combination of those standard B-splines~$N^{\ab{p},\ab{r}}_{j_1,j_2}$, $j_1,j_2=0,1,2$, which 
have been added twice, to get a well defined $C^2$-smooth function~$\phi_{\ab{v}^{(i)}}$, cf. \cite[Section~4.4]{KaVi19a}. 

In~\cite{KaVi19a}, the vertex subspace~$\W_{\ab{v}^{(i)}}$ has been defined as the span of all functions~\eqref{eq:defphiXi} whose spline 
functions~$g_{i_{\ell}}$ fulfill the homogeneous linear system~\eqref{eq:vertex_homogeneous_system}. A drawback of the resulting space is its geometry-dependent dimension. 
Below, we will construct a vertex subspace~$\W_{\ab{v}^{(i)}}$, whose dimension will be independent of the 
initial geometry. This will be achieved by additionally enforcing that the functions of the generated space~$\W_{\ab{v}^{(i)}}$ have to be $C^4$-smooth at the vertex~$\ab{v}^{(i)}$. 
This strategy can be seen as an extension of approach~\cite{KaSaTa17c}, where $C^1$-smooth functions have been generated in the vicinity of a 
vertex~$\ab{v}^{(i)}$ by additionally enforcing $C^2$-smoothness of the functions at the vertex.

Let $\psi_{j_1,j_2}: \overline{\Omega} \to \R$, $j_1,j_2=0,\ldots,4$ with $j_1+j_2 \leq 4$ be functions which are $C^2$-smooth on $\overline{\Omega}$ and even $C^4$-smooth at the 
vertex~$\ab{v}^{(i)}$ such that
\[
 \partial_{x_1}^{m_1} \partial_{x_2}^{m_2} \psi_{j_1,j_2} (\ab{v}^{(i)}) = \sigma^{j_1+j_2} \delta^{m_1}_{j_1} \delta^{m_2}_{j_2},
\]
where $\delta^{m}_{j}$ is the Kronecker delta, and $\sigma$ is a scaling factor given by
\[
 \sigma = \left(\frac{h}{p \, \nu_i} \sum_{\ell=0}^{\nu_i -1} ||\nabla \ab{F}^{(i_{\ell})}(\ab{0}) || \right)^{-1},
\]
cf.~\cite[Definition~14]{KaSaTa17c}. We define the isogeometric function~$\phi_{\ab{v}^{(i)}_{I};j_1,j_2}$, $j_1,j_2=0,1,\ldots,4$, $j_1+j_2 \leq 4$, as
\begin{equation} \label{eq:phi_vertex_I}
 \phi_{\ab{v}^{(i)}_I;j_1,j_2}(\ab{x}) = \phi_{\ab{v}^{(i)}}(\ab{x}),
\end{equation}
where the function~$\phi_{\ab{v}^{(i)}}$ is specified by fulfilling the interpolation conditions
\begin{equation} \label{eq:interpolation_conditions}
 \partial^{m_1}_{x_1} \partial^{m_2}_{x_2} \phi_{\ab{v}^{(i)}} (\ab{v}^{(i)}) = \partial_{x_1}^{m_1} \partial_{x_2}^{m_2} \psi_{j_1,j_2} (\ab{v}^{(i)}), \mbox{ }
 m_1,m_2=0,1,\ldots,4, \mbox{ }m_1+m_2 \leq 4.
\end{equation}
The coefficients~$a^{\Gamma^{(i_\ell)}}_{j_1,j_2}$ and $a_{j_1,j_2}^{(i_\ell)}$ of the spline functions~$g_{i_{\ell}}$ of the isogeometric function~$\phi_{\ab{v}^{(i)}_{I};j_1,j_2}$
are uniquely determined by the conditions~\eqref{eq:interpolation_conditions}, and can be computed via the following equivalent interpolation problem
\small
\begin{align} 
   \partial_{\xi_1}^{m_1}  \partial_{\xi_2}^{m_2} \, g_{i_\ell}^{\Gamma^{(i_\ell)}}  (\bfm{0}) &=
     \partial_{\xi_1}^{m_1}   \partial_{\xi_2}^{m_2} \left( \psi_{j_1,j_2} \circ \bfm{F}^{(i_\ell)} \right)  (\bfm{0}),\ \;
   0 \leq m_1 \leq 4, \;  0 \leq m_2 \leq 2, \;  m_1+m_2 \leq 4, \nonumber \\[-0.3cm] 
 & \label{eq:systemEq} \\[-0.3cm]   
        \partial_{\xi_1}^{m_1}  \partial_{\xi_2}^{m_2} \, g_{i_\ell}^{\Omega^{(i_{\ell})}} (\bfm{0}) &=
     \partial_{\xi_1}^{m_1}   \partial_{\xi_2}^{m_2} \left( \psi_{j_1,j_2} \circ \bfm{F}^{(i_\ell)} \right) (\bfm{0}), \ \;
   0 \leq m_1,m_2 \leq 2, \nonumber
\end{align}
\normalsize
for $\ell = 0,1,\ldots, \nu_{i}-1$. 
The presented construction of the isogeometric function~$\phi_{\ab{v}^{(i)}_{I};j_1,j_2}$ works and leads to a well-defined $C^2$-smooth function on~$\overline{\Omega}$, 
which is additionally $C^4$-smooth at the vertex~$\ab{v}^{(i)}$. This is true because of the following four reasons: First, since the functions 
$g_{\Gamma^{(i_{\ell})};0,j_2}^{(i_{\ell})}$, $j_2=0,1,\ldots, 4$, $g_{\Gamma^{(i_{\ell})};1,j_2}^{(i_{\ell})}$, $j_2=0,1,2, 3$, and $g_{\Gamma^{(i_{\ell})};2,j_2}^{(i_{\ell})}$, 
$j_2=0,1, 2$, of $g_{i_{\ell}}^{\Gamma^{(i_{\ell})}}$ interpolate the traces of the form $N_{j_2}^{p,r+2}$, the first derivatives of the form~$N_{j_2}^{p-1,r+1}$ and the second 
derivatives of the form~$N_{j_2}^{p-2,}$, respectively, across the interface~$\Gamma^{(i_{\ell})}$, cf.~\cite[Section~7]{KaVi17c}, the first set of interpolation conditions 
in~\eqref{eq:systemEq} uniquely determine all coefficients $a^{\Gamma^{(i_\ell)}}_{j_1,j_2}$. Moreover, as the functions~$g_{i_{\ell}}^{\Gamma^{(i_\ell)}}$ ensure $C^2$-smoothness 
across the interfaces~$\Gamma^{(i_{\ell})}$, we also obtain 
\small
\begin{equation} \label{eq:systemEq2}
  \partial_{\xi_1}^{m_1}  \partial_{\xi_2}^{m_2} \, g_{i_\ell}^{\Gamma^{(i_{\ell+1})}}  (\bfm{0}) =
     \partial_{\xi_1}^{m_1}   \partial_{\xi_2}^{m_2} \left( \psi_{j_1,j_2} \circ \bfm{F}^{(i_\ell)} \right)  (\bfm{0}),
\end{equation}
\normalsize
for $0 \leq m_1 \leq 2$, $0 \leq m_2 \leq 4$, $ m_1+m_2 \leq 4$ and $\ell=0,1,\ldots, \nu_{i}-1$. Second, all coefficients~$a_{j_1,j_2}^{(i_\ell)}$ are specified by the 
second set of interpolation conditions in~\eqref{eq:systemEq}. Third, thanks to the fulfillment of conditions~\eqref{eq:systemEq} and \eqref{eq:systemEq2} and the fact 
that the function~$g_{i_{\ell}}$ is given by~\eqref{eq:g_vertex}, we get 
\begin{equation*} 
 \partial_{\xi_1}^{m_1}  \partial_{\xi_2}^{m_2} \, g_{i_\ell}   (\bfm{0}) =
    \partial_{\xi_1}^{m_1}   \partial_{\xi_2}^{m_2} \left( \psi_{j_1,j_2} \circ \bfm{F}^{(i_\ell)} \right) (\bfm{0}),
\end{equation*}
for $0 \leq m_1 \leq 4$, $0 \leq m_2 \leq 4$, $m_1+m_2 \leq 4$ and $\ell=0,1,\ldots, \nu_{i}-1$,
which is further equivalent to satisfying conditions~\eqref{eq:interpolation_conditions}. Moreover, we obtain that the function~$\phi_{\ab{v}^{(i)}_{I};j_1,j_2}$ is well defined 
and $C^2$-smooth on $\overline{\Omega}$, where the latter property is the direct consequence of satisfying the homogeneous linear system~\eqref{eq:vertex_homogeneous_system}.
Fourth, the $C^4$-smoothness of the function~$\phi_{\ab{v}^{(i)}_{I};j_1,j_2}$ at vertex~$\ab{v}^{(i)}$ follows directly from the fulfillment of the interpolation 
conditions~\eqref{eq:interpolation_conditions}.

Finally, we define for each inner vertex~$\ab{v}^{(i)}$, i.e. $i \in \mathcal{I}_{\Xi_{I}}$, the vertex subspace~$\W_{\ab{v}^{(i)}}$ as
\[
 \W_{\ab{v}^{(i)}} =  \Span \{ \phi_{\ab{v}^{(i)}_{I};j_1,j_2} |\; \; j_1,j_2=0,1,\ldots,4, \; j_1+j_2 \leq 4 \}.
\]

Let us continue with the case of a boundary vertex~$\ab{v}^{(i)}$ of patch valency~$\nu_{i} \geq 3$, that is $i \in \mathcal{I}_{\Xi_3}$. This can be treated similarly to 
the case of an inner vertex as above. The vertex subspaces~$\W_{\ab{v}^{(i)}}$, $i \in \mathcal{I}_{\Xi_3}$, are defined as
 \[
 \W_{\ab{v}^{(i)}} =  \Span \{ \phi_{\ab{v}^{(i)}_{3};j_1,j_2} |\; \; j_1,j_2=0,1,\ldots,4, \; j_1+j_2 \leq 4 \},
\]
where the functions~$\phi_{\ab{v}^{(i)}_{3};j_1,j_2}$ are again determined just via
\begin{equation} \label{eq:phi_vertex_3}
 \phi_{\ab{v}^{(i)}_3;j_1,j_2}(\ab{x}) = \phi_{\ab{v}^{(i)}}(\ab{x}),
\end{equation}
by interpolating the conditions~\eqref{eq:interpolation_conditions}. The only two differences in the construction of the functions are that the lower index~$\ell$ of the indices 
$i_{\ell}$ is not taken modulo~$\nu_i$, and that for the patches~$\Omega^{(i_0)}$ and $\Omega^{(i_{\ell-1})}$ the functions~$g^{(i_0)}_{\Gamma^{(i_0)};j_1,j_2}$ and 
$g^{(i_{\nu_{i}-1})}_{\Gamma^{(i_{\nu_i})};j_1,j_2}$ in~\eqref{eq:g_vertex}, which are given on the boundary edges~$\Gamma^{(i_0)}$ or $\Gamma^{(i_{\nu_{i}})}$, are just defined 
as standard B-splines, i.e.
\[
 g^{(i_0)}_{\Gamma^{(i_0)};j_1,j_2}(\xi_2,\xi_1) = N_{j_1,j_2}^{\ab{p},\ab{r}}(\xi_2,\xi_1) \quad \mbox{ and } \quad
 g^{(i_{\nu_{i}-1})}_{\Gamma^{(i_{\nu_i})};j_1,j_2}(\xi_1,\xi_2) = N_{j_1,j_2}^{\ab{p},\ab{r}}(\xi_1,\xi_2).
\]
Analogously to the case of an inner vertex, one can argue that the functions~$\phi_{\ab{v}^{(i)}_3;j_1,j_2}$ are $C^2$-smooth on~$\overline{\Omega}$ and even $C^4$-smooth at the 
vertex~$\ab{v}^{(i)}$.

In case of a boundary vertex~$\ab{v}^{(i)}$ of patch valency~$\nu_i=2$, i.e. $i \in \mathcal{I}_{\Xi_2}$, we assume for the two neighboring 
patches~$\Omega^{(i_0)}$ and $\Omega^{(i_1)}$ enclosing the vertex~$\ab{v}^{(i)}$, that the corresponding geometry mappings~$\ab{F}^{(i_0)}$ and $\ab{F}^{(i_1)}$ are parameterized as 
shown in Fig.~\ref{fig:situation_two_patches}~(left), that we have $\ab{v}^{(i)} = \ab{F}^{(i_0)}(\ab{0}) = \ab{F}^{(i_1)}(\ab{0})$, and that the common interface is denoted by 
$\Gamma^{(j_0)}$, i.e. $\overline{\Gamma^{(j_0)}}=\overline{\Omega^{(i_0)}} \cap \overline{\Omega^{(i_1)}}$. Then, we construct the vertex~subspace~$\W_{\ab{v}^{(i)}}$ as
\begin{align*} \label{eq:spaceWXB}
\mathcal{W}_{\bfm{v}^{(i)}} =&  \Span \{ \phi_{\bfm{v}^{(i)}_2; j_1,j_2} |\;   
j_1=0,1,\ldots,4, \; j_2 =  0,1,\ldots, 4-\min(j_1,2) \} ,
\end{align*}
where $\phi_{\bfm{v}^{(i)}_2; j_1,j_2}$ is given as
\begin{equation} \label{eq:phi_vertex_2}
 \phi_{\bfm{v}^{(i)}_2;j_1,j_2} (\bfm{x}) = 
    \begin{cases} 
     \phi_{\Gamma^{(j_0)}_{I};j_1,j_2}(\ab{x}) \quad \mbox{ if }j_1=0,1,2,\\
     \phi_{\Omega^{(i_{0})};j_1+\lfloor\frac{j_2}{2} \rfloor,j_2 \bmod 2}(\ab{x}) \quad \mbox{ if }j_1=3,\\
     \phi_{\Omega^{(i_{1})};j_1-1+\lfloor\frac{j_2}{2}\rfloor,j_2 \bmod 2}(\ab{x}) \quad \mbox{ if }j_1=4,
    \end{cases}
\end{equation}
with the functions~$\phi_{\Omega^{(i_{\ell})};j_1,j_2}$ and $\phi_{\Gamma^{(j_0)}_{I};j_1,j_2}$ as defined in~\eqref{eq:defphiOmega} and \eqref{eq:basisFunctionsGenericEdge}, 
respectively. The functions~$\phi_{\ab{v}^{(i)}_2;j_1,j_2}$ are $C^2$-smooth on $\overline{\Omega}$, since the functions~$\phi_{\Gamma^{(j_0)}_{I};j_1,j_2}$, $j_1=0,1,2$, 
$j_2=0,\ldots,4-j_1$, are $C^2$-smooth on $\overline{\Omega}$ by construction, and the functions $\phi_{\Omega^{(i_{\ell})};j_1,j_2}$, $j_1=3,4$, $j_2=0,1, \ldots 1-j_1+3$, $\ell=0,1$,
possess vanishing values, gradients and Hessians along all interfaces~$\overline{\Gamma^{(j)}}$, $j \in \mathcal{I}_{\Gamma_I}$

Finally, for the case of a boundary vertex $\bfm{v}^{(i)}$ of patch valency $\nu_i=1$, i.e. $i \in \mathcal{I}_{\Xi_1}$, and assuming without loss of generality that 
the vertex~$\ab{v}^{(i)}$ is given by $\ab{v}^{(i)} = \ab{F}^{(i_0)}(\ab{0})$, the vertex subspace~$\mathcal{W}_{\bfm{v}^{(i)}}$ 
is of the form
\begin{align*} 
\mathcal{W}_{\bfm{v}^{(i)}} =&  \Span \{ \phi_{\bfm{v}^{(i)}_1; j_1,j_2} |\;   
j_1,j_2 =  0,1,\ldots, 4, \; j_1+j_2 \leq 4 \} ,
\end{align*}
where $\phi_{\bfm{v}^{(i)}_1;j_1,j_2}$ is defined as 
\begin{equation} \label{eq:phi_vertex_1}
 \phi_{\bfm{v}^{(i)}_1;j_1,j_2} (\bfm{x}) = 
    \begin{cases}
   (N_{j_1,j_2}^{p,r}\circ (\ab{F}^{(i_0)})^{-1})(\bfm{x}) \;
   \mbox{ if }\f \, \bfm{x} \in \overline{\Omega^{(i_0)}},
    \\ 0 \quad \mbox{ otherwise. }
    \end{cases}
\end{equation}
The functions~$\phi_{\bfm{v}^{(i)}_1;j_1,j_2}$ are trivially $C^2$-smooth on~$\overline{\Omega}$, since the functions have vanishing values, gradients and Hessians along all 
interfaces~$\overline{\Gamma^{(j)}}$, $j \in \mathcal{I}_{\Gamma_I}$.

One could generate for the case of a boundary vertex~$\ab{v}^{(i)}$ of patch valency~$\nu_i=2$ or $\nu_{i}=1$ a different vertex subspace~$\W_{\ab{v}^{(i)}}$ as described above 
by just following also for these cases the strategy explained for a boundary vertex of patch valency~$\nu_i \geq 3$. However, we will use in our numerical examples, 
see Section~\ref{section_Numerical_examples}, the construction presented above, since it is simple and does not require the solving of a system of linear equations.

\subsection{Basis and dimension of the space~$\W$} \label{subsec:basis_dimension}

Since the space~$\W$ is the direct sum~\eqref{eq:Wh}, a basis of the space~$\W$ can be obtained by just constructing bases of the individual 
subspaces~$\W_{\Omega^{(i)}}$, $i \in \mathcal{I}_{\Omega}$, $\W_{\Gamma^{i}}$, $i \in \mathcal{I}_{\Gamma}$, and $\W_{\ab{v}^{(i)}}$, $i \in \mathcal{I}_{\Xi}$. 
As for each collection of functions~\eqref{eq:defphiOmega}, \eqref{eq:defphiGammaBoundary}, \eqref{eq:basisFunctionsGenericEdge}, \eqref{eq:phi_vertex_I}, 
\eqref{eq:phi_vertex_3}, \eqref{eq:phi_vertex_2} or \eqref{eq:phi_vertex_1}, spanning one of the subspaces~$\W_{\Omega^{(i)}}$, $\W_{\Gamma^{(i)}}$ or $\W_{\ab{v}^{(i)}}$, the single 
functions are linearly independent by construction, the functions form a basis of the corresponding subspace. Moreover, the intersection of any two of the individual 
subspaces does not have any common function except the zero function. Therefore, the entire collection of all functions, i.e. 
\[
  \phi_{\Omega^{(i)};j_1,j_2}, \; j_1,j_2=3,4,\ldots,n-4, \; i \in \mathcal{I}_{\Omega};
\]
\[
 \phi_{\Gamma^{(i)}_B;j_1,j_2},\;j_1=0,1,2, \;j_2=5-j_1,\ldots, n+j_1-6, \;i \in \mathcal{I}_{\Gamma_B};
\]
\[
 \phi_{\Gamma^{(i)}_I;j_1,j_2},\;j_1=0,1,2, \;j_2=5-j_1,\ldots, n_{j_{1}}+j_1-6, \;i \in \mathcal{I}_{\Gamma_I};
\]
\[
 \phi_{\ab{v}^{(i)}_{I};j_1,j_2},\;j_1,j_2=0,1,\ldots,4, \;j_1+j_2 \leq 4, \; i \in \mathcal{I}_{\Xi_I};
 \]
 \[
 \phi_{\ab{v}^{(i)}_{3};j_1,j_2},\;j_1,j_2=0,1,\ldots,4, \;j_1+j_2 \leq 4, \; i \in \mathcal{I}_{\Xi_3};
\]
\[
\phi_{\ab{v}^{(i)}_{2};j_1,j_2},\;j_1=0,1,\ldots,4, \;j_2=0,1,\ldots,4-\min(j_1,2),\; i \in \mathcal{I}_{\Xi_2};\]
and
\[
 \phi_{\ab{v}^{(i)}_{1};j_1,j_2},\;j_1,j_2=0,1,\ldots,4, \;j_1+j_2 \leq 4, \; i \in \mathcal{I}_{\Xi_1}; 
\]
builds a basis of the space~$\W$, cf.~\cite[Theorem~1]{KaVi19a}. 

In addition, the decomposition~\eqref{eq:Wh} of the space~$\W$ implies that 
\[
 \dim \W = \sum_{i \in \mathcal{I}_{\Omega}} \dim \W_{\Omega^{(i)}} + \sum_{i \in \mathcal{I}_{\Gamma}} \dim \W_{\Gamma^{(i)}} + 
 \sum_{i \in \mathcal{I}_{\Xi}} \dim \W_{\ab{v}^{(i)}}.
\]
Since the dimensions of the single subspaces~$\W_{\Omega^{(i)}}$, $i \in \mathcal{I}_{\Omega}$, $\W_{\Gamma^{i}}$, $i \in \mathcal{I}_{\Gamma}$, and 
$\W_{\ab{v}^{(i)}}$, $i \in \mathcal{I}_{\Xi}$, are given by
\[
 \dim \W_{\Omega^{(i)}} = (n-6)^2,
\]
\[
 \dim \W_{\Gamma^{(i)}} = \begin{cases}
                           3(n-8) \; \mbox{ if } i \in \mathcal{I}_{\Gamma_{B}}, \\
                           3(n-2k-9) \;  \mbox{ if } i \in \mathcal{I}_{\Gamma_{I}},
                          \end{cases}
\]
and
\[
 \dim \W_{\ab{v}^{(i)}} = \begin{cases}
                           15 \;  \mbox{ if }i \in \mathcal{I}_{\Xi_I} \cup \mathcal{I}_{\Xi_1} \cup \mathcal{I}_{\Xi_3}, \\
                           18 \;  \mbox{ if }i \in \mathcal{I}_{\Xi_2},
                          \end{cases}
\]
respectively, the dimension of~$\W$ is equal to
\begin{equation}  \label{eq:dimension_formula} 
\dim \W = | \mathcal{I}_{\Omega}| (n-6)^2 + 3 | \mathcal{I}_{\Gamma^I}| (n-2k-9) +3  | \mathcal{I}_{\Gamma^B}| (n-8) + 
15 \left( | \mathcal{I}_{\Xi_I}|  + | \mathcal{I}_{\Xi_1}|  + | \mathcal{I}_{\Xi_3}|  \right)  + 18 | \mathcal{I}_{\Xi_2}|.
\end{equation} 
Therefore, in contrast to~\cite{KaVi19a}, the dimension of the $C^2$-smooth space~$\W$ is independent of the initial geometry.

\section{Isogeometric collocation} \label{sec:collocation}

We will present the idea of using the concept of isogeometric collocation for solving the Poisson's equation on planar multi-patch domains.

\subsection{Problem statement and isogeometric formulation}  \label{subsec:problem_statement}

Let $f:\Omega \to \R$, $f_1: \partial \Omega \to \R$. We are interested in finding  $u:\overline{\Omega} \to \R$, $u \in C^2(\overline{\Omega})$, which solves the 
Poisson's equation
\begin{align} \label{eq:Poisson}
\triangle u (\bfm{x}) & =   f(\bfm{x}), \quad \bfm{x} \in \Omega,  \nonumber \\
& \\[-0.4cm]
u (\bfm{x}) & = f_1(\bfm{x}) ,   \quad \bfm{x} \in \partial \Omega. \nonumber
\end{align}
The idea is to employ isogeometric collocation to compute a $C^2$-smooth approximation $u_h \in \mathcal{W}$ of the solution~$u$. This requires the use of global collocation 
points~$\bfm{\glob}_j$, $j \in \mathcal{J}$, which are separated into two distinct sets of collocation points, namely of inner collocation points $\bfm{\glob}_j^I$,  
$j \in \mathcal{J}_I$, and of boundary collocation points $\bfm{\glob}_j^B$,  $j \in \mathcal{J}_B$. Inserting these points into \eqref{eq:Poisson}, we obtain 
\begin{align} \label{eq:PoissonPoints}
\triangle u (\bfm{\glob}_j^I) & =   f(\bfm{\glob}_j^I), \quad j \in \mathcal{J}_I,\nonumber \\
& \\[-0.4cm]
u (\bfm{\glob}_j^B) & = f_1(\bfm{\glob}_j^B) ,   \quad j \in \mathcal{J}_B.  \nonumber 
\end{align}
To use the isogeometric approach for solving problem~$\eqref{eq:PoissonPoints}$, we have to express each global collocation point~$\bfm{\glob}_{j}$ with respect to local coordinates via
$$
\bfm{\loc}^{(\delta(j));I}_{j} = \left( \bfm{F}^{(\delta(j))}\right)^{-1} \left( \bfm{\glob}_j^I \right) , \;  j \in \mathcal{J}_I, \quad \mbox{ and } \quad
\bfm{\loc}^{(\delta(j));B}_{j} = \left( \bfm{F}^{(\delta(j))}\right)^{-1} \left( \bfm{\glob}_j^B \right) , \;  j \in \mathcal{J}_B,
$$
where 
\begin{equation*}  
\delta: \mathcal{J} \to \mathcal{I}_\Omega, \quad
\delta(j) = \min \{ i \in \mathcal{I}_\Omega, \; \bfm{\glob}_j \in \overline{\Omega^{(i)}} \}.
\end{equation*}
Having the local collocation points $\bfm{\loc}^{(\delta(j));I}_{j}$ and $\bfm{\loc}^{(\delta(j));B}_{j}$, then equations \eqref{eq:PoissonPoints} can be transformed into 
\small
\begin{equation*}
\frac{1}{\left| \det J^{(\delta(j))} \left(\bfm{\loc}^{(\delta(j));I}_{j}\right) \right|}  \left( \nabla \circ \left( N^{(\delta(j))}(\bb{\xi})  
\nabla (u \circ \bfm{F}^{(\delta(j))})(\bb{\xi}) \right)\right){\bigg|}_{\bfm{\xi} =\bfm{\loc}^{(\delta(j));I}_{j} }= 
f\left(\ab{F}^{(\delta(j))}\left(\bfm{\loc}^{(\delta(j));I}_{j}\right)\right), \mbox{ }j \in \mathcal{J}_{I},
\end{equation*}
\normalsize
\begin{equation*}
u \left( \bfm{F}^{(\delta(j))}  \left(\bfm{\loc}^{(\delta(j));B}_{j} \right)\right)  = f_1\left(\ab{F}^{(\delta(j))}\left(\bfm{\loc}^{(\delta(j));B}_{j}\right)\right),
 \mbox{ }j \in \mathcal{J}_{B},
\end{equation*}
where $J^{(\ell)}$ is the Jacobian of $\ab{F}^{(\ell)}$ and 
\begin{equation*}
N^{(\ell)}(\bb{\xi}) = \left(J^{(\ell)}(\bb{\xi})\right)^{-T}  \left(J^{(\ell)}(\bb{\xi})\right)^{-1} \left|\det J^{(\ell)}(\bb{\xi})\right|.
\end{equation*}
This gives rise to a linear system for the unknown coefficients $c_i$ of the approximation $u_h = \sum_{i \in \mathcal{I}} c_i \phi_i$, where 
$\mathcal{I} = \{0,1, \ldots, \dim \mathcal{W}-1 \}$, and $\{\phi_i\}_{i\in \mathcal{I}}$ is the basis of $\mathcal{W}$ presented in Section~\ref{subsec:basis_dimension}. 

\subsection{Selection of collocation points} \label{sec:SelectionCollocPoints}

The choice of the collocation points plays an important role in the stability and convergence behavior of the numerical solution. Below, we will first present the general concept 
of choosing these points, and will then describe in more detail two particular selections of the collocation points. 

\subsubsection{The general concept}

The idea is to choose for each patch~$\overline{\Omega^{(i)}}$, $i \in \mathcal{I}_{\Omega}$, global collocation points $\ab{F}^{(i)}(\bfm{\loc}_{\bfm{j}})$, 
which are determined by the local collocation points
\begin{equation} \label{eq:local_points}
\bfm{\loc}_{\bfm{j}} = (\loc_{j_1},\loc_{j_2}) \in [0,1]^2,
\end{equation}
which are again composed of univariate collocation points $\loc_j \in [0,1]$. In case that we obtain the same global collocation point for more than one patch, we keep 
the corresponding point just for the patch~$\overline{\Omega^{(i)}}$ with the smallest index~$i$, to avoid repetitions of points. Moreover, we split the resulting points into inner 
and boundary collocation points, cf.~Section~\ref{subsec:problem_statement}.

In the one-patch case, mainly two different types of (univariate) collocation points have been considered. First, standard collocation points such as Greville abscissae or 
Demko abscissae are used, see e.g.~\cite{IsoCollocMethods2010,EnzoKiendlLorenzis2017,Reali2015}, and possess a convergence behavior with 
respect\footnote{Instead of studying the convergence behavior 
with respect to the $L^2$, $H^1$ and $H^2$ norm, often equivalently, the $L^{\infty}$, $W^{1,\infty}$ and $W^{2,\infty}$-norm is considered in the literature.} to 
the $L^2$, $H^1$ and $H^2$ norm of orders~$\mathcal{O}(h^{p-1})$, $\mathcal{O}(h^{p-1})$ and $\mathcal{O}(h^{p-1})$ and of orders $\mathcal{O}(h^{p})$, $\mathcal{O}(h^{p})$ 
and $\mathcal{O}(h^{p-1})$ for odd and even spline degree~$p$, respectively. These convergence rates are suboptimal for the $L^2$ norm and in the case of odd spline 
degree~$p$ also for the $H^1$ norm in comparison to the Galerkin approach, where the rates are of orders~$\mathcal{O}(h^{p+1})$, $\mathcal{O}(h^{p})$ and 
$\mathcal{O}(h^{p-1})$. 

The second type of collocation points are the so-called 
superconvergent points \cite{SuperConvergent2015,GomezLorenzisVariationalCollocation,MonSanTam2017}, 
which provide in case of odd spline degree~$p$ better convergence rates than the Greville abscissae. In more detail, in case of odd spline degree, the convergence 
behavior for the approach~\cite{GomezLorenzisVariationalCollocation} is of orders $\mathcal{O}(h^{p})$, $\mathcal{O}(h^{p})$ and $\mathcal{O}(h^{p-1})$ with respect to 
$L^2$, $H^1$ and $H^2$ norm, and the convergence rates for the methods~\cite{SuperConvergent2015,MonSanTam2017} are even of optimal orders, 
i.e. $\mathcal{O}(h^{p+1})$, $\mathcal{O}(h^{p})$ and $\mathcal{O}(h^{p-1})$ with respect to $L^2$, $H^1$ and $H^2$ norm.
The basic concept of the superconvergent points is to choose as collocation points 
those points (or more precisely approximations of those points), which are the roots of the Galerkin residual~$D^2(u-u_h)$ of the considered problem. 
Since the resulting set of superconvergent points is approximately twice as large as the number of degrees of freedom, different strategies have 
been proposed in \cite{SuperConvergent2015,GomezLorenzisVariationalCollocation,MonSanTam2017} to use the superconvergent points as collocation points for solving the 
collocation problem~\eqref{eq:PoissonPoints}. In~\cite{SuperConvergent2015}, the entire set of points is used, which leads to an overdetermined linear system that is 
solved by means of the least-squares method. In contrast, in \cite{GomezLorenzisVariationalCollocation,MonSanTam2017}, particular subsets of superconvergent points are 
selected to obtain the same number of collocation points as the number of degrees of freedom. While the superconvergent points are selected 
in~\cite{GomezLorenzisVariationalCollocation} in an alternating way, the points are chosen in \cite{MonSanTam2017} in a clustered way. 

Below, we will extend the use of Greville abscissae and of superconvergent points as collocation points to the case of planar multi-patch domains. One main difference 
will be that splines of maximal smoothness, i.e. $r=p-1$, as usually employed in the one-patch case, cannot be used for the multi-patch case. In the latter case, the 
underlying spline spaces~$\mathcal{S}^{p,r}_h$ have to fulfill $2 \leq r \leq p-3$ with $p \geq 5$ to allow the construction of $h$-refineable $C^2$-smooth spline 
spaces~$\W$, cf. Section~\ref{subsec:multipatch} and Section~\ref{subsec:C2space}. For the sake of simplicity, we will restrict ourselves to the cases 
$(p,r) \in \{ (5,2), (6,2) , (6,3)\}$. In addition, the number of collocation points will be slightly larger than the number of degrees of freedom, cf. 
Section~\ref{subsub:Greville} and Section~\ref{subsub:superconvergent}. Therefore, the resulting linear system will be overdetermined and will be solved as 
in~\cite{SuperConvergent2015} by means of the least-squares method. As already observed in~\cite{SuperConvergent2015}, we first have to solve those 
equations in~\eqref{eq:PoissonPoints} which have been determined by the boundary collocation points, i.e. the second set of 
equations, and then use the result to solve the remaining system of linear equations. 

\subsubsection{Greville abscissae} \label{subsub:Greville}

\begin{figure}[htb]
\centering
\includegraphics[width=10cm]{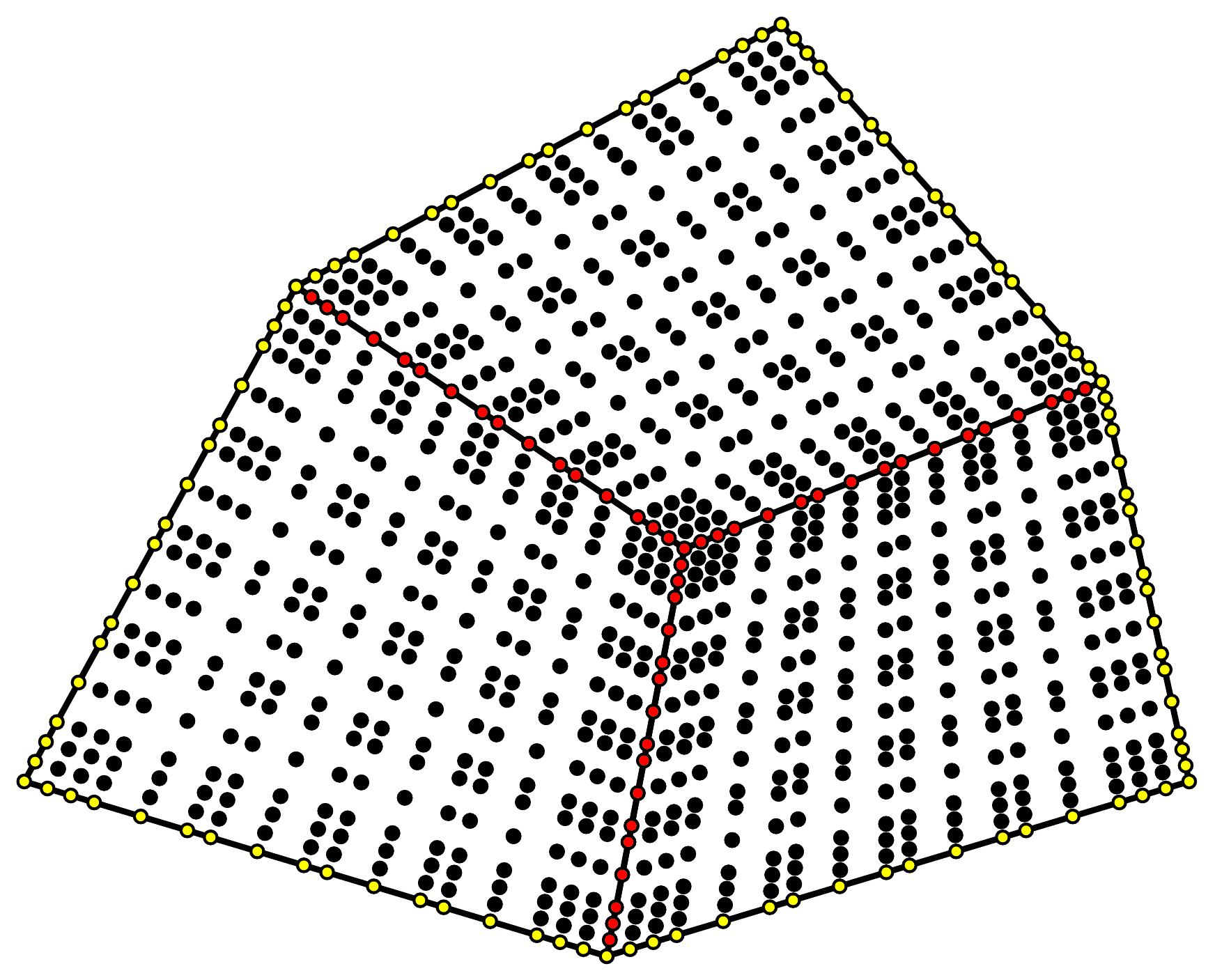} 
\caption{Set of all mapped Greville points $\ab{F}^{(i)}(\bfm{\loc}_{\bfm{j}}^{(5,2)})$ (without repetitions) for a bilinearly parameterized three-patch domain and $k=4$, 
where the boundary collocation points and the collocation points on the interfaces are specified in yellow and red, respectively.}
\label{fig:three-patch_Greville}
\end{figure}

For each patch~$\overline{\Omega^{(i)}}$, $i \in \mathcal{I}_{\Omega}$, the local collocation points~$\bfm{\loc}_{\ab{j}}$, $\ab{j} \in \{0,1,\ldots,n-1\}^2$, 
from~\eqref{eq:local_points} are chosen as the tensor-product Greville points
\[
 \bfm{\loc}_{\ab{j}}^{\ab{p},\ab{r}} = (\loc_{j_1}^{p,r},\loc_{j_2}^{p,r})
\]
with $\loc_{j}^{p,r}=\frac{t_{j+1}^{p,r} + \ldots t_{j+p}^{p,r}}{p}$, $j \in \{0,1,\ldots, n-1 \}$. The resulting global collocation 
points~$\ab{F}^{(i)}(\bfm{\loc}^{\ab{p},\ab{r}}_{\ab{j}})$ for a particular example of a bilinearly parameterized three-patch domain are shown in 
Fig.~\ref{fig:three-patch_Greville}. 
As already mentioned in the previous section, the number of (global) collocation points (without repetitions of same points) is slightly larger than the number of basis 
functions of the $C^2$-smooth space~$\W$. But the quotient of the number of collocation points and of the number of basis functions converges to $1$ when $k$ grows as 
the following example will demonstrate.

\begin{ex}
Let $\overline{\Omega}$ be a multi-patch domain with exactly one inner vertex of edge and patch valency $\nu$, see e.g. the three-, five- and six-patch domain in 
Fig~\ref{fig:bilinear_domains}~(a)--(c). In addition, let $p=5$ and $r=2$, which implies $n=6+3k$. Since $\dim \mathcal{W} = \nu (9 k^2+21 k+12)+15$ 
(compare formula~\eqref{eq:dimension_formula}) and the number of collocation points (without repetitions), i.e. $|\mathcal{J}|$, is equal to 
$\nu n^2 - \nu n +1 = \nu (9 k^2+33k +30) +1$,
we obtain
\[
 \lim_{k \to \infty} \frac{|\mathcal{J}|}{\dim \W} = \lim_{k \to \infty} \, \frac{\nu (9 k^2+33k +30) +1}{\nu (9 k^2+21 k+12)+15}  = 1.
\]
Note that the same conclusion is not only true for $(p,r)=(6,2)$ and $(p,r)=(6,3)$ but also for any multi-patch domain, since the number of collocation points is just 
slightly larger than the number of basis functions along the edges and in the vicinity of inner vertices and boundary vertices of patch valencies greater or equal to two.
\end{ex}

The numerical results in Section~\ref{section_Numerical_examples} will show that we will obtain for the multi-patch case the same convergence behavior as for the 
one-patch case, namely with respect to the $L^2$, $H^1$ and $H^2$ norm, convergence rates of orders~$\mathcal{O}(h^{p-1})$, $\mathcal{O}(h^{p-1})$ and 
$\mathcal{O}(h^{p-1})$ and of orders $\mathcal{O}(h^{p})$, $\mathcal{O}(h^{p})$ and $\mathcal{O}(h^{p-1})$ for odd and even spline degree~$p$, respectively.

\subsubsection{Superconvergent points} \label{subsub:superconvergent}

\begin{figure}
\centering
\begin{tabular}{c}
\includegraphics[width=16.2cm,clip]{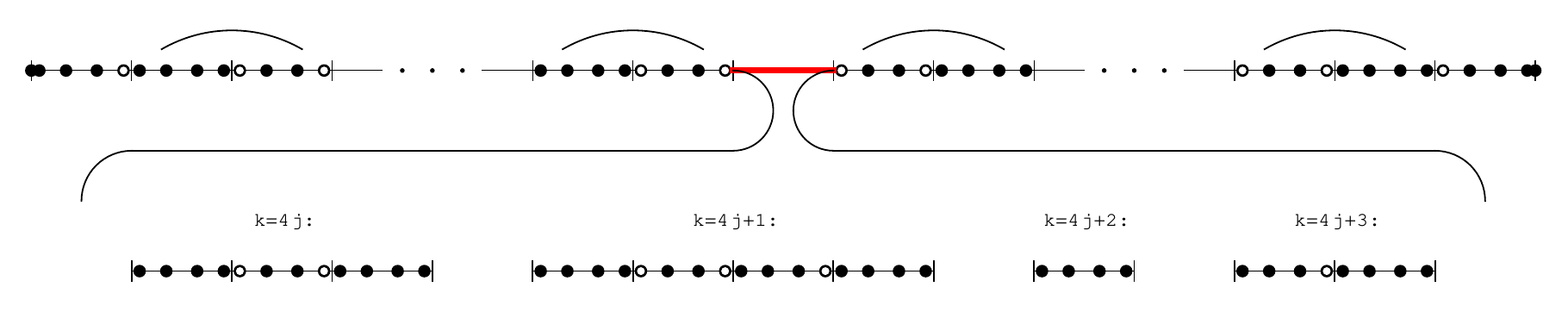} 
\end{tabular}
\caption{All superconvergent points (black and white points) and clustered superconvergent points (only black points) for the case $(p,r)=(5,2)$. The center part (in red) 
is given by the second row depending on the number~$k$ of different inner knots.}
\label{fig:clusteredSuperconv}
\end{figure}

\begin{figure}
\centering
\begin{tabular}{c}
\includegraphics[width=16.2cm,clip]{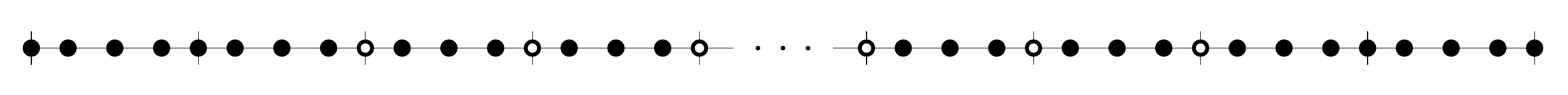} \\
\includegraphics[width=16.2cm,clip]{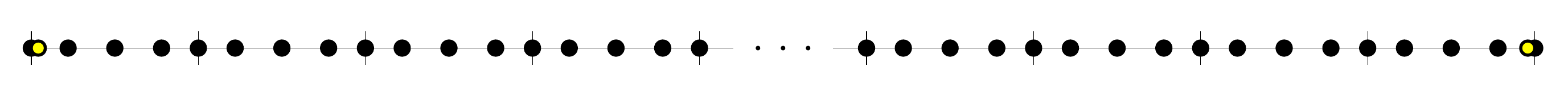} 
\end{tabular}
\caption{First row: All superconvergent points (black and white points) and clustered superconvergent points (only black points) for $(p,r)=(6,3)$. Second row: All 
superconvergent points (black and yellow points), which coincide with the clustered superconvergent points, for the case $(p,r)=(6,2)$, where the yellow points specify 
the two added Greville abscissae.}
\label{fig:clusteredSuperconv6}
\end{figure}

The superconvergent points (cf. \cite{SuperConvergent2015,GomezLorenzisVariationalCollocation,MonSanTam2017}) are defined as the roots of the Galerkin 
residual~$D^2(u-u_h)$ of the considered problem. 
Estimates of the univariate superconvergent points 
can be computed by solving a simple 1D Poisson's equation such as 
\begin{equation} \label{eq:1DPoisson}
u''(x) =f(x), \quad x \in (0,1), \quad {\rm and} \quad u(0)=u(1)=0, 
\end{equation}
with some particular function~$f$. The superconvergent points have been mainly considered so far just for splines of maximal smoothness, i.e., $r=p-1$. Since the 
$C^2$-smooth spline spaces~$\W$ require underlying spline spaces~$\mathcal{S}^{p,r}_{h}$ with $2 \leq r \leq p-3$ and $p \geq 5$, cf. Section~\ref{subsec:multipatch} and 
Section~\ref{subsec:C2space}, we will generate these points for some of these spaces, namely for the cases $(p,r) \in \{ (5,2), (6,2) , (6,3)\}$. To first  
estimate the superconvergent points, we just follow the approach presented in \cite{GomezLorenzisVariationalCollocation} based on a particular 1D Poisson's 
problem~\eqref{eq:1DPoisson}. Further studying the asymptotic behavior of the resulting estimates 
of the superconvergent points for the different levels of refinement, we obtain that the superconvergent points on each knot span with respect to the reference 
interval~$[-1,1]$ are given as the roots of the polynomial $15 x^4-12 x^2 + 1$ for $p=5$ and $r=2$ and as the roots of the polynomial $99x^5-130x^3+31x$ for $p=6$ and $r=2,3$.
The roots of these polynomials, and hence the locations of the superconvergent points on each knot span with respect to the reference interval~$[-1,1]$ are presented in 
Table~\ref{tab:superconvergent_points}.
All superconvergent points for 
$(p,r)=(5,2)$ and $(p,r) \in \{(6,2),(6,3) \}$ are shown in Fig.~\ref{fig:clusteredSuperconv} and \ref{fig:clusteredSuperconv6}, respectively.
\begin{table}[htb]
\begin{center}
\begin{tabular}{c|c}  
 ($p$,$r$) &  Superconvergent points \\[0.1cm]
 \hline\\[-0.3cm]
(5,2) &  $\pm\sqrt{\frac{1}{15} (6 + \sqrt{21})}$, \quad $\pm\sqrt{\frac{1}{15} (6 - \sqrt{21})}$ \\[0.2cm]
(6,2) &  $0$, \quad $\pm 1$, \quad $ \pm \sqrt{\frac{31}{99}}$ \\[0.2cm]
(6,3) &  $0$, \quad $\pm 1$, \quad $ \pm \sqrt{\frac{31}{99}}$ \\
\end{tabular} 
\end{center}
\caption{Locations of the superconvergent points on each knot span w.r.t. the reference interval $[-1,1]$ for the cases $(p,r) \in \{ (5,2),(6,2),(6,3)\}$.}
\label{tab:superconvergent_points}
\end{table}
In case of $(p,r) \in \{(6,2),(6,3) \}$, the knots of the underlying spline spaces~$\mathcal{S}_{h}^{p,r}$ are superconvergent points, and therefore, we count these 
points in case of inner knots just once for two neighboring segments. In case of $(p,r)=(5,2)$, since the boundary points of the whole domain are not contained 
in the set of all superconvergent points, they have to be added to the set, see Fig.~\ref{fig:clusteredSuperconv}, to be able to impose Dirichlet boundary conditions.

The differences $\Delta_{(p,r)}$ between the number of all superconvergent points and the dimension of the underlying spline space~$\mathcal{S}_{h}^{p,r}$ (i.e. 
$\dim \mathcal{S}_{h}^{p,r} = p+1+k(p-r)$) are given by
\begin{eqnarray*}
& \Delta_{(5,2)} = 4(k+1)+2 - (6+3k) = k,\\
& \Delta_{(6,2)} = 4(k+1)+1 - (7+4k) = -2,\\
& \Delta_{(6,3)} = 4(k+1) +1 - (7+3k) = k -2.
\end{eqnarray*}
The idea is to select from the superconvergent points in a clustered way a subset of collocation points such that the number of collocation points is equal to the number of basis 
functions of $\mathcal{S}_{h}^{p,r}$, and such that the selection of points does not ruin the convergence properties as in \cite{MonSanTam2017} for the case of splines of 
maximal smoothness. While in the case of $(p,r)=(6,2)$, we have even to add two points to the set of superconvergent points, in the cases $(p,r)=(5,2)$ and $(p,r)=(6,3)$ we have to 
omit $k$ and $k-2$ points, respectively. We will call the resulting sets of collocation points clustered superconvergent points similarly to~\cite{MonSanTam2017}. 

Let us first consider the case~$(p,r)=(6,2)$. We choose as additional points for the set of superconvergent points just the second and the second last Greville abscissa. 
This strategy has been already used for selecting additional points close to the boundary for Dirichlet boundary problems in case of splines of maximal smoothness, 
cf. \cite{GomezLorenzisVariationalCollocation,MonSanTam2017}. Note that for $(p,r)=(6,2)$ the set of all superconvergent points coincides with the set of clustered 
superconvergent points.

Recall that for the cases $(p,r)=(5,2)$ and $(p,r)=(6,3)$ we have to omit $k$ and $k-2$ points, respectively. We denote by $S_k$ the indices of the points which will be 
omitted. Furthermore, we can assume $k \geq 2$. For $(p,r)=(5,2)$, we choose depending on the number~$k$ of different inner knots the set $S_k$ as 
$$
S_k = \widetilde{S}_k \cup \left\{  \begin{array}{ll}   
\{2k+1,2k+4\}, & k=4j, \\  
\{2k-1,2k+2,2k+6\}, & k=4j+1, \\  
\{\}, & k=4j+2, \\  
\{2k+2\}, & k=4j+3,
\end{array}
\right.
$$
with 
$$
  \widetilde{S}_k = \{4,4k+1\} \cup \bigcup_{i=0}^{\lfloor\frac{k-6}{4} \rfloor} \left( \{ 9+8i, 12+8i \} \cup  \{ 4k+5-(9+8i), 4k+5-(12+8i) \} \right).
$$
In case of $(p,r)=(6,3)$, we omit all inner knots except the first and the last one, which leads to
$$
 S_k =  \{  9+4i, \; i=0,1,\ldots,k-3\}.
$$

All superconvergent points as well as the clustered superconvergent points are shown in Fig.~\ref{fig:clusteredSuperconv} and \ref{fig:clusteredSuperconv6} 
for the cases $(p,r)=(5,2)$ and $(p,r) \in \{ (6,2),(6,3) \}$, respectively. The following example will demonstrate that these points possess for the one-patch case 
the same convergence behavior with respect to the $L^2$, $H^1$ and $H^2$ norm as the sets of all superconvergent points and of the clustered superconvergent points in  
the case of splines of maximal smoothness. Moreover, it will be shown that also for the case of Greville abscissae we get the same convergence rates for the splines of 
maximal smoothness and for the splines with a reduced regularity.
\begin{ex} \label{ex:onepatch}
We perform isogeometric collocation on the one-patch domain in Fig.~\ref{fig:onepatch_domain}~(first row), which is just the bilinearly parameterized unit 
square~$[0,1]^2$, and compare the convergence behavior under mesh-refinement for different spline spaces and different sets of collocation points. 
We use underlying spline spaces $\mathcal{S}^{p,r}_h([0,1]^2)$ with $p=5,6$ and $2\leq r \leq p-3$ with mesh sizes $h=\frac{1}{5}, \frac{1}{10}, 
\frac{1}{20},\frac{1}{40}$ for $(p,r)=(5,2)$ and  $h=\frac{1}{4}, \frac{1}{8}, \frac{1}{16},\frac{1}{32}$ for $p=6$ and $r=2,3$. Fig.~\ref{fig:onepatch_domain} shows 
the resulting relative $L^2$, $H^1$ and $H^2$ errors for the exact solution   
\[
 u(x_1,x_2)=\cos (4x_1-2) \sin \left(4x_2-\frac{2}{3}\right),
\]
(see first row), by using the Greville points (second row), all superconvergent points (third row) and the clustered superconvergent points (fourth row) as collocation 
points. We get as in the case of splines of maximal regularity, see~\cite{SuperConvergent2015,MonSanTam2017}, the same rates of convergence depending on the used 
spline degree~$p$. That is for $p=5$ in case of Greville points the same rates of order $\mathcal{O}(h^{4})$ in the $L^2$, $H^1$ and $H^2$ norm, and in case of all 
or clustered superconvergent points the optimal rates of orders $\mathcal{O}(h^{6})$, $\mathcal{O}(h^{5})$ and $\mathcal{O}(h^{4})$ in the $L^2$, $H^1$ and $H^2$-norm, 
respectively. For $p=6$ we have the same rates independent of one of the three sets of collocation points, that is $\mathcal{O}(h^{6})$, $\mathcal{O}(h^{6})$ and 
$\mathcal{O}(h^{5})$ in the $L^2$, $H^1$ and $H^2$ norm, respectively.
 \begin{figure}
\centering\footnotesize
\begin{tabular}{cc}
\includegraphics[width=4.0cm,clip]{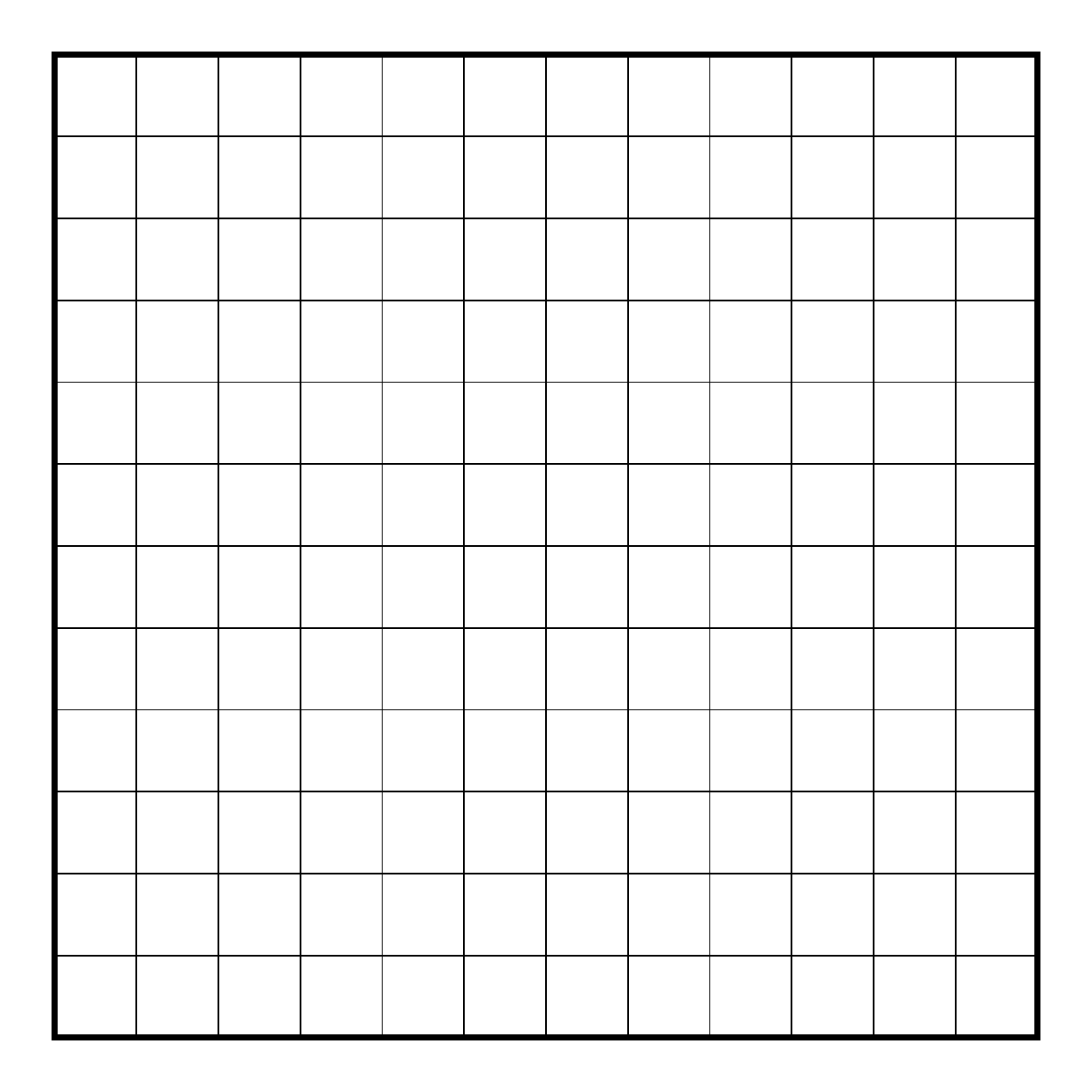} &
\includegraphics[width=6.0cm,clip]{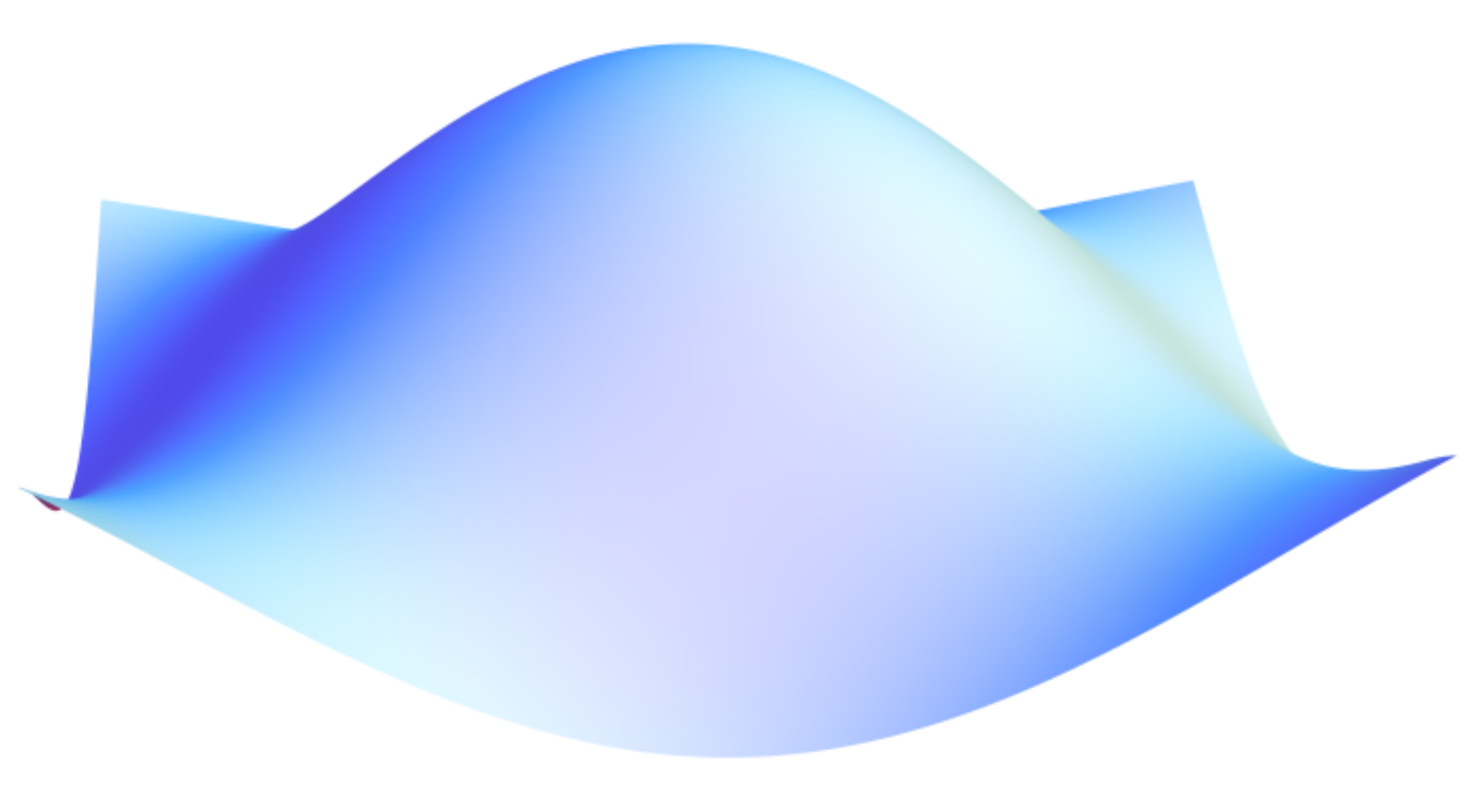} \\
Computational domain & Exact solution 
\end{tabular}
\begin{tabular}{ccc}
\includegraphics[width=4.5cm,clip]{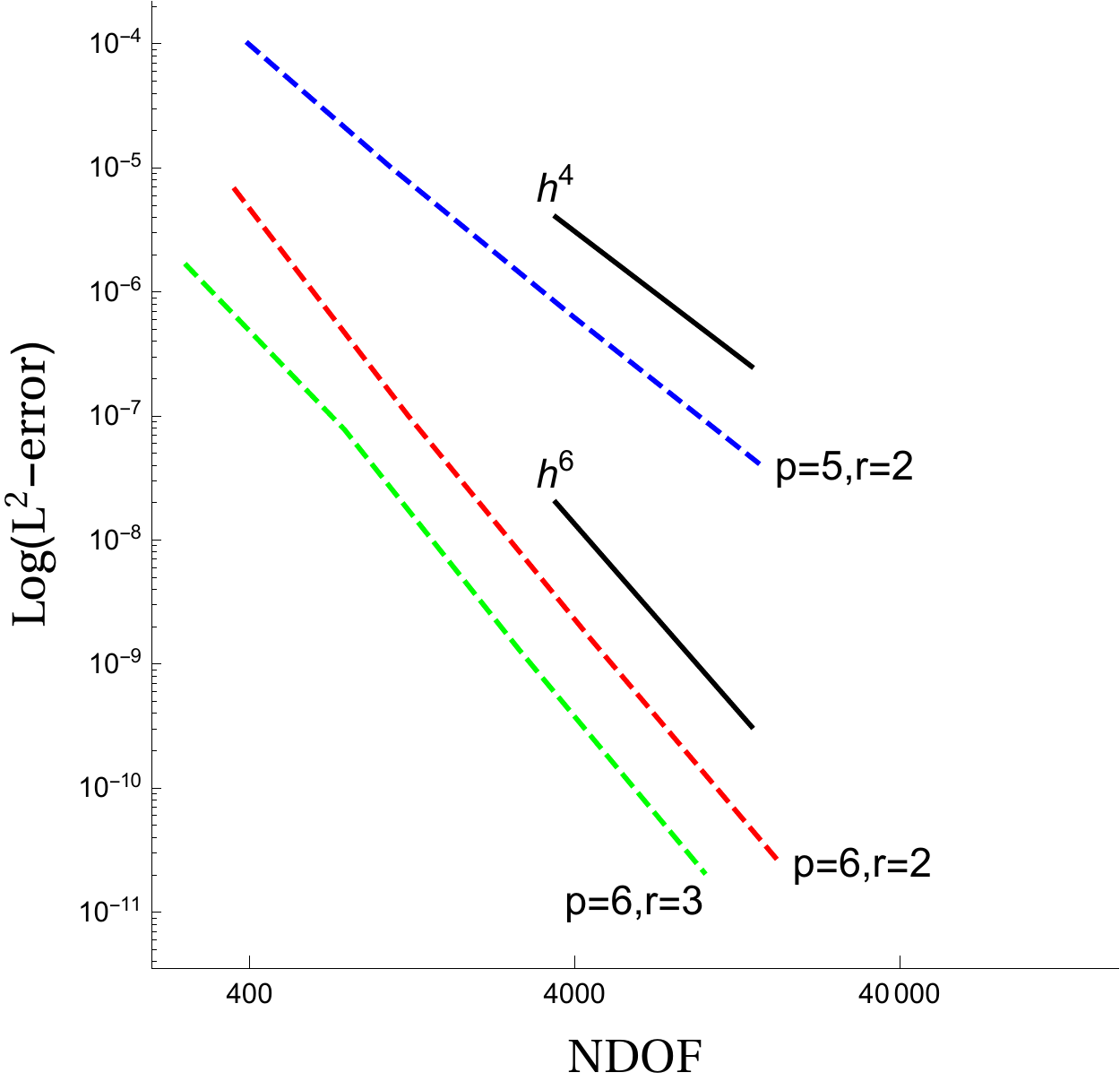} &
\includegraphics[width=4.5cm,clip]{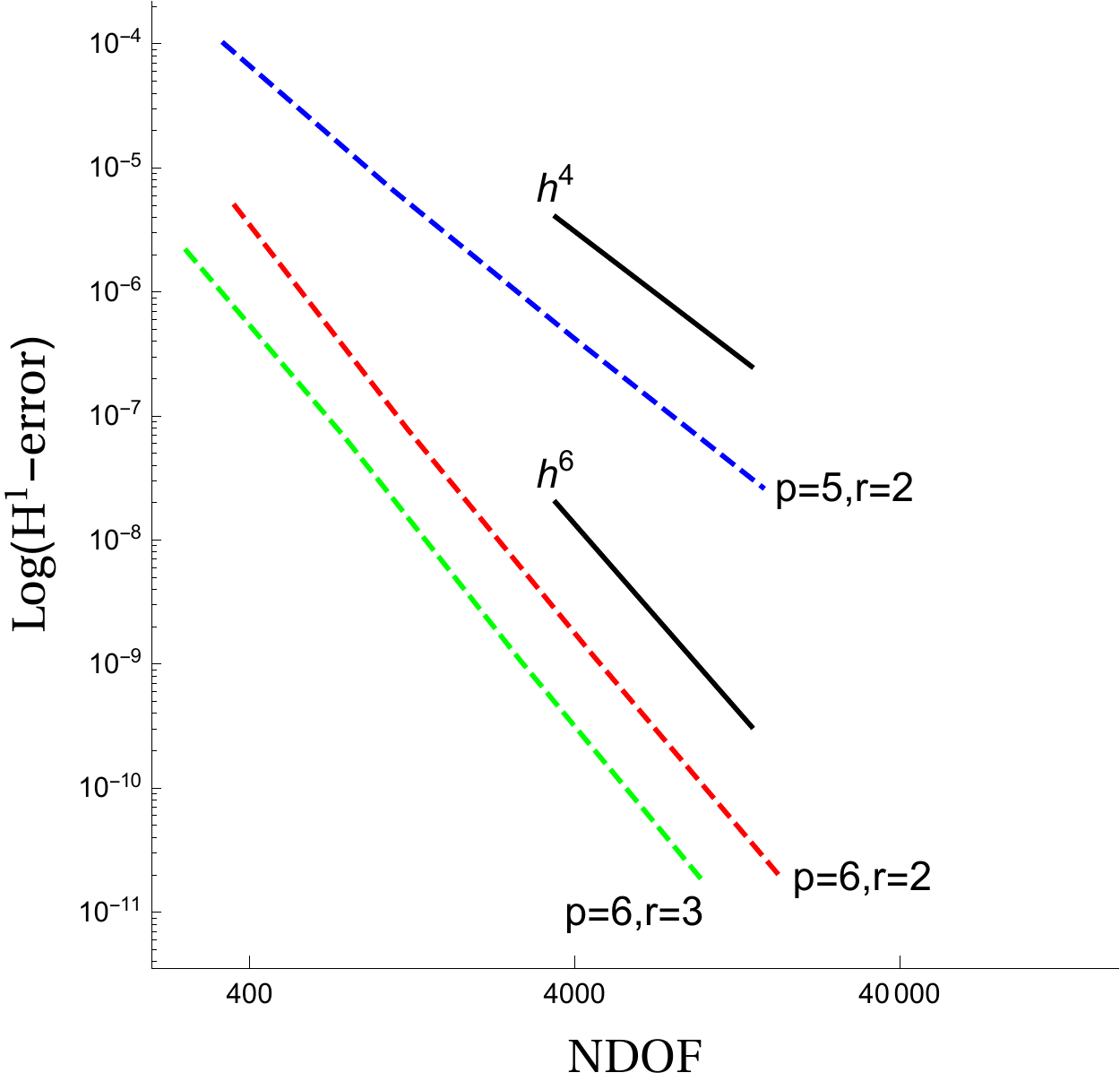} &
\includegraphics[width=4.5cm,clip]{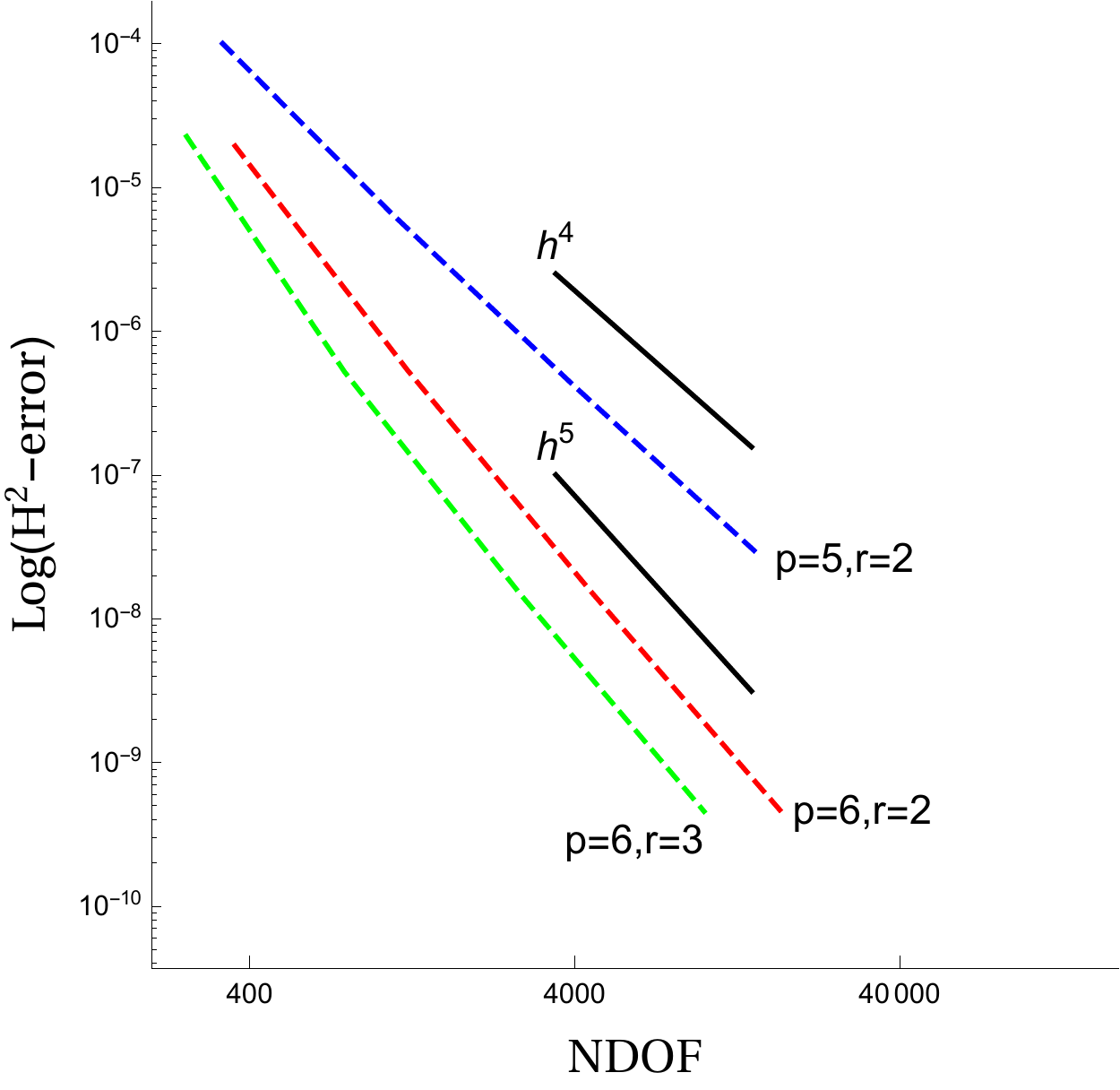}\\
\multicolumn{3}{c}{Relative $L^2$, $H^1$ and $H^2$ errors for collocation at Greville points} \\
\includegraphics[width=4.5cm,clip]{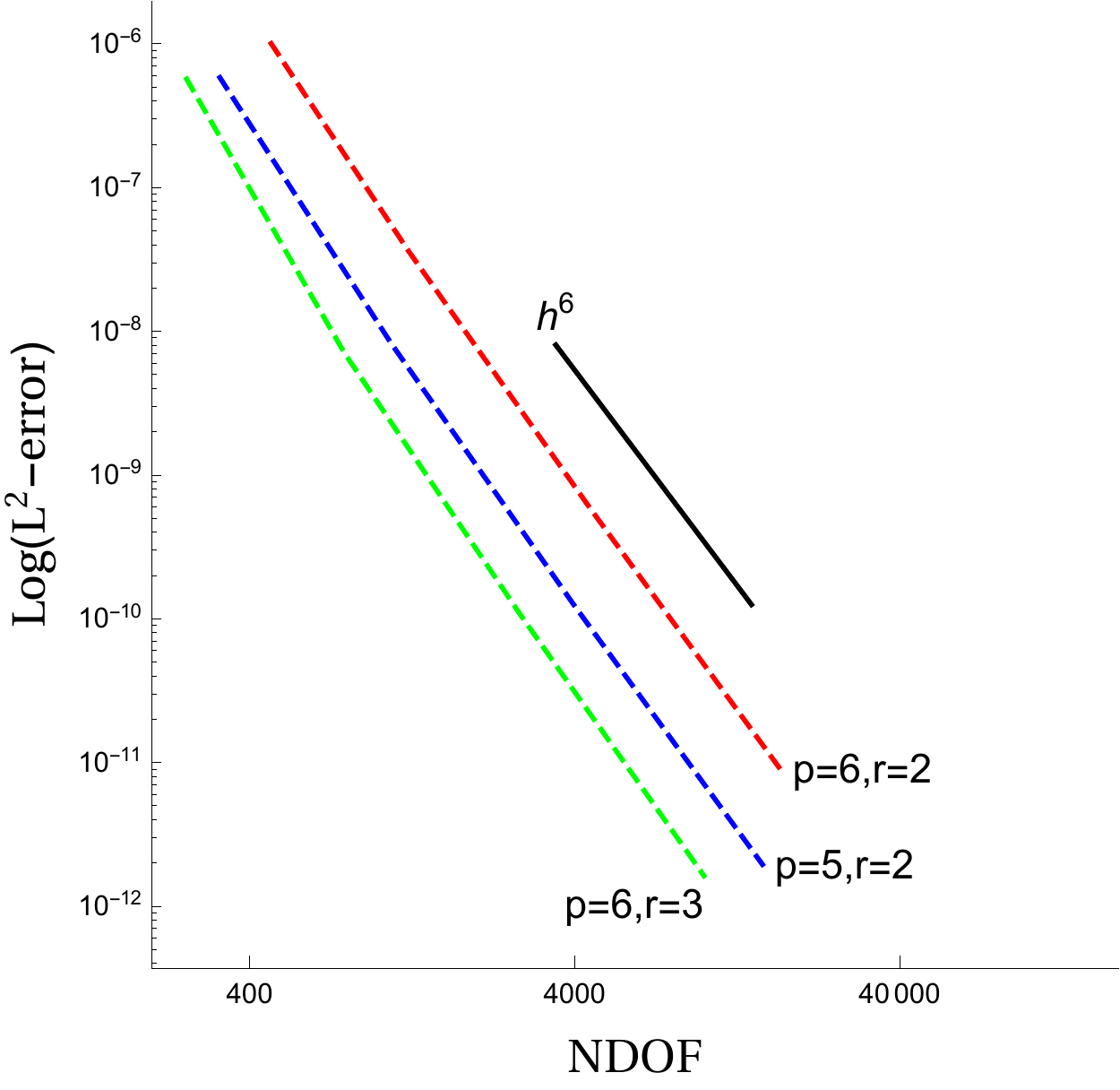} &
\includegraphics[width=4.5cm,clip]{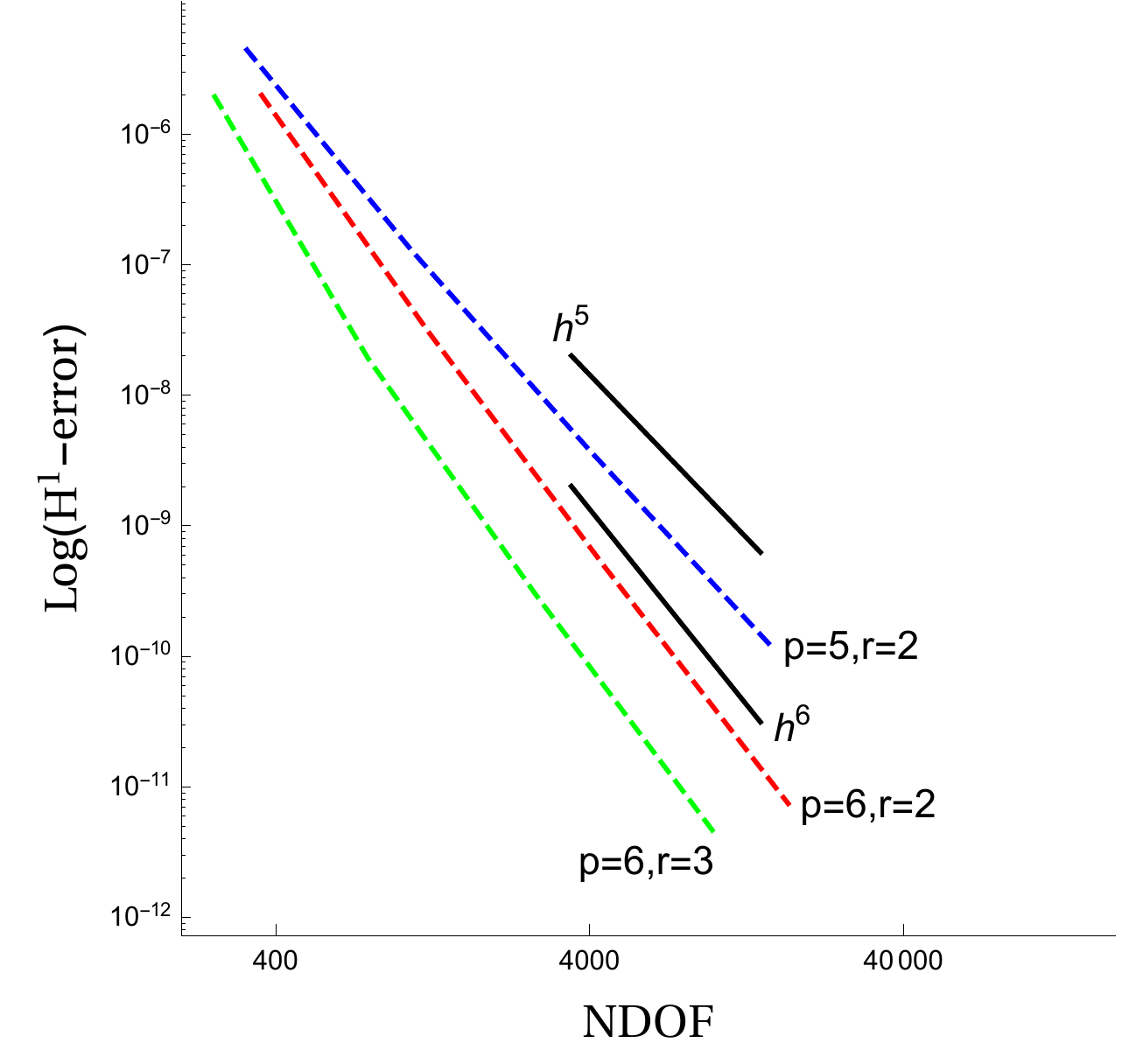} &
\includegraphics[width=4.5cm,clip]{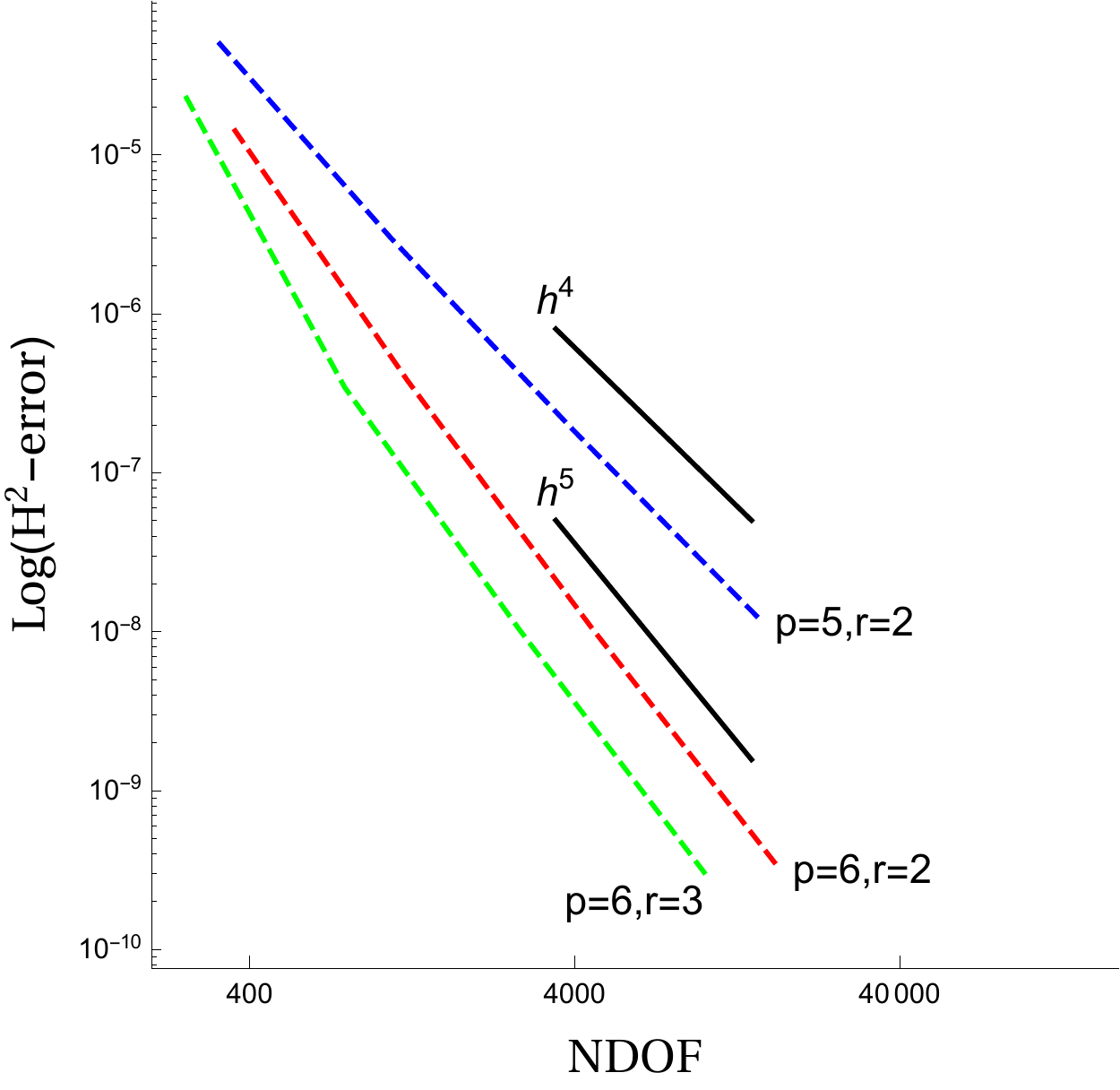}\\
\multicolumn{3}{c}{Relative $L^2$, $H^1$ and $H^2$ errors for collocation at all superconvergent points} \\
\includegraphics[width=4.5cm,clip]{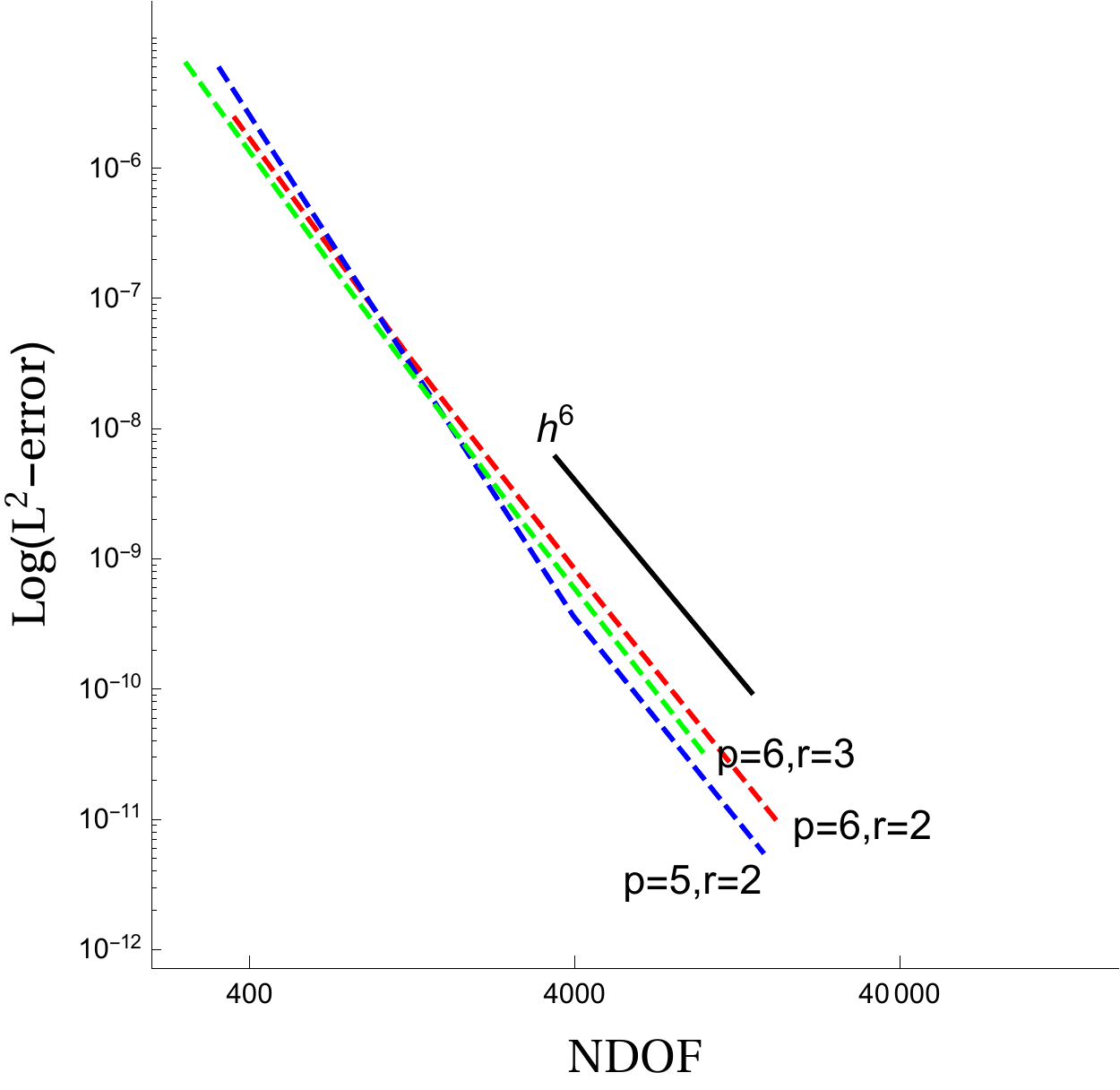} &
\includegraphics[width=4.5cm,clip]{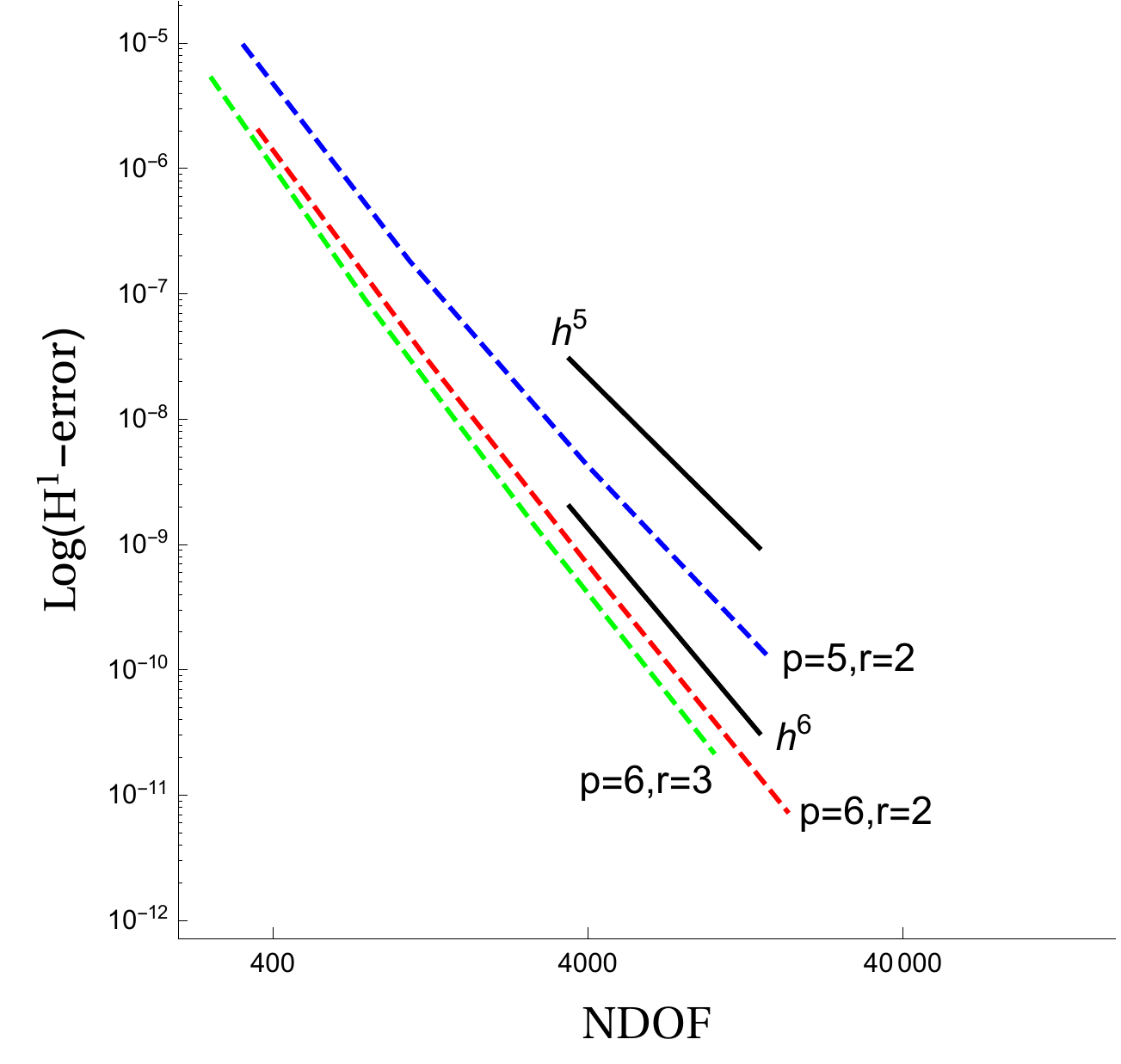} &
\includegraphics[width=4.5cm,clip]{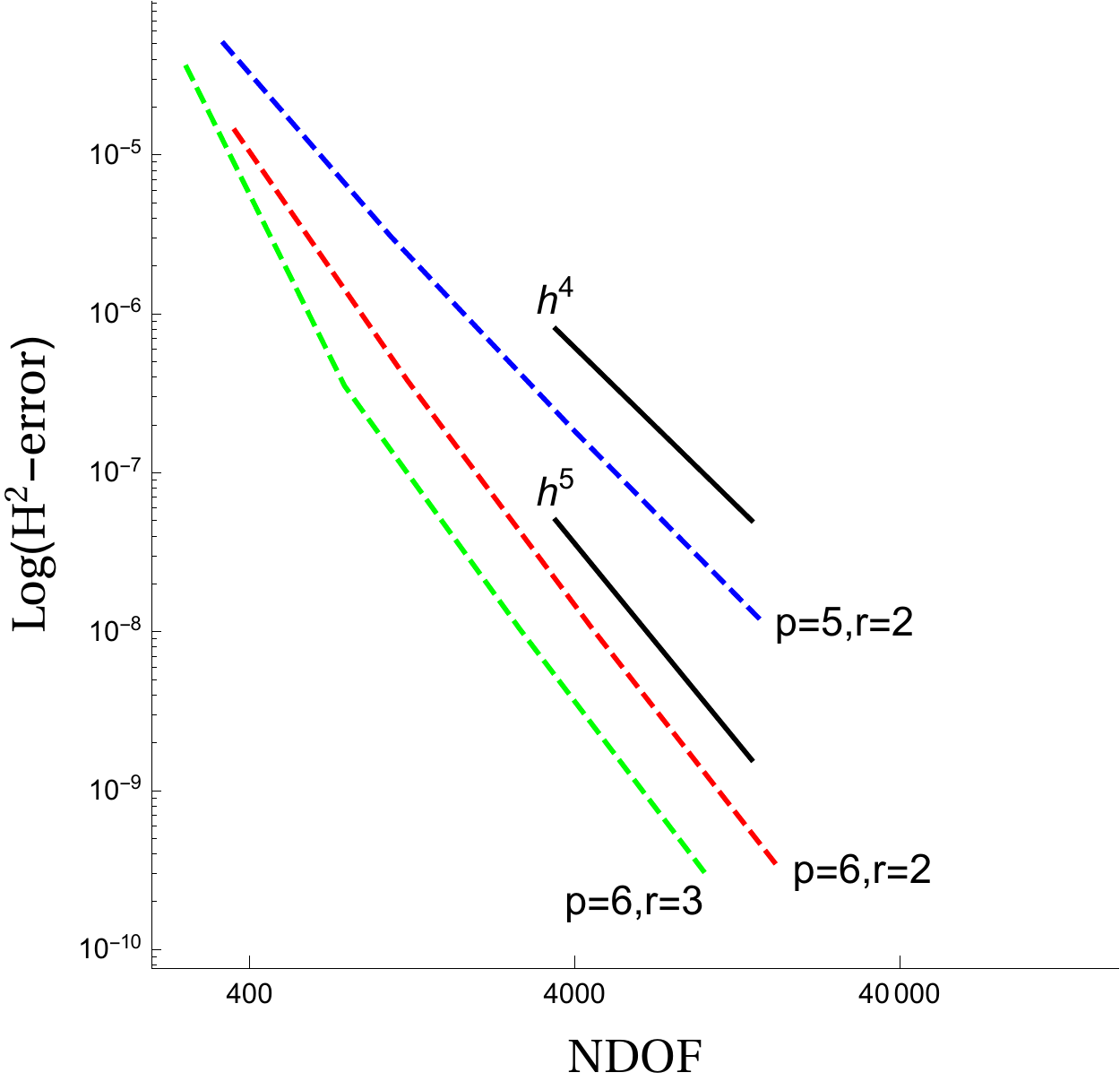}\\
\multicolumn{3}{c}{Relative $L^2$, $H^1$ and $H^2$ errors for collocation at clustered superconvergent points}
\end{tabular}
\caption{Isogeometric collocation on the bilinearly parameterized unit-square using different sets of collocation points, cf. Example.~\ref{ex:onepatch}.}
\label{fig:onepatch_domain}
\end{figure}
\end{ex}
In the case of multi-patch domains, see Section~\ref{section_Numerical_examples}, the convergence rates for the sets of superconvergent points will be the 
same as in the one-patch case, except for one particular case, namely for $(p,r)=(5,2)$ with respect to the $L^2$ norm, where the rate will be 
just of order~$\mathcal{O}(h^p)$ instead of $\mathcal{O}(h^{p+1})$. In the numerical examples in Section~\ref{section_Numerical_examples}, for the sake 
of brevity, we will restrict ourselves just to the case of clustered superconvergent points, since experiments using all superconvergent points have led to 
the same convergence rates, and the number of clustered superconvergent points is dramatically smaller than the number of all superconvergent points for the 
cases~$(p,r)=(5,2)$ and $(p,r)=(6,3)$. However, in the case of multi-patch domains, also for the clustered superconvergent points the number of points is larger 
than the number of degrees of freedom. More precisely, the number of clustered superconvergent points is equal to the number of Greville abscissae, 
and therefore again, the quotient of the number of collocation points and of the number of degrees of freedom converges to one when the number~$k$ of different 
inner knots grows. For solving the resulting overdetermined linear system, we use as for the Greville abscissae the least-squares method.

\section{Numerical examples}  \label{section_Numerical_examples}

We will perform isogeometric collocation on several multi-patch domains by using the Greville points and the clustered superconvergent points as collocation points 
(cf.~Section~\ref{subsub:Greville} and Section~\ref{subsub:superconvergent}), and will study the convergence under mesh refinement. Since the considered sets 
of collocation points will lead to an overdetermined system of linear equations, we will use the least-squares approach to solve the system, see 
Section~\ref{sec:SelectionCollocPoints}. Below, we will denote the $C^2$-smooth space~$\mathcal{W}$ for a specific mesh size $h$ by $\mathcal{W}_h$. 

\begin{ex} \label{ex:multipatchdomain}
We consider the four bilinearly parameterized multi-patch domains~(a)-(d) shown in Fig.~\ref{fig:bilinear_domains}~(first row). For all domains we generate a sequence 
of $C^2$-smooth spaces $\mathcal{W}_h$ with mesh sizes $h=\frac{1}{5},\frac{1}{10},\frac{1}{20},\frac{1}{40}$ for $p=5$ and $r=2$ and with mesh sizes 
$h=\frac{1}{4},\frac{1}{8},\frac{1}{16},\frac{1}{32}$ for $p=6$ and $r=2,3$. The generated $C^2$-smooth spaces are then used to perform isogeometric collocation on 
the different multi-patch domains for a right side function $f$ obtained from the exact solutions~$u$ visualized in Fig.~\ref{fig:bilinear_domains}~(a)-(d), which are 
given for the single domains~(a)-(d) by
\[
 u_a(x_1,x_2)= -4 \cos \left(\frac{2x_1}{3}\right) \sin \left(\frac{2x_2}{3}\right),
\]
\[
 u_b(x_1,x_2)= -4 \cos \left(\frac{x_1}{2}\right) \sin \left(\frac{x_2-2}{2}\right),
\]
\[
 u_c(x_1,x_2)= -4 \cos \left(\frac{x_1+3}{2}\right) \sin \left(\frac{x_2-1}{2}\right),
\]
and
\[
u_d(x_1,x_2)= -4 \cos \left(\frac{x_1}{3}\right) \sin \left(\frac{x_2}{3}\right),
\]
respectively. The resulting relative $L^2$, $H^1$ and $H^2$ errors by using Greville points and clustered superconvergent points are shown in Fig.~\ref{fig:results_Greville} and 
\ref{fig:results_Cauchy}, respectively. In case of Greville points the estimated convergence rates are the same as for the one-patch instance in Example~\ref{ex:onepatch}, 
that is for $p=5$ the same rates of order~$\mathcal{O}(h^{4})$ in the $L^2$, $H^1$ and $H^2$ norm, and for $p=6$ the rates of orders $\mathcal{O}(h^{6})$, 
$\mathcal{O}(h^{6})$ and $\mathcal{O}(h^{5})$ in the $L^2$, $H^1$ and $H^2$ norm, respectively. In case of clustered superconvergent points, we observe for both spline 
degrees $p=5,6$ convergence rates of orders~$\mathcal{O}(h^{p})$, $\mathcal{O}(h^{p})$ and $\mathcal{O}(h^{p-1})$ in the $L^2$, $H^1$ and $H^2$ norm, respectively, which 
means for $p=5$ a reduction of the order by~$1$ with respect to the $L^2$ norm compared to the one-patch instance in Example~\ref{ex:onepatch}. 
Note that using all superconvergent points would lead to the same reduced rate for $p=5$ as employing clustered superconvergent points, and therefore will not be shown for the 
sake of brevity.

\begin{figure}[htp]
\centering\footnotesize
\begin{tabular}{cccc}
(a) & (b) & (c) & (d) \\
\includegraphics[width=3.8cm,clip]{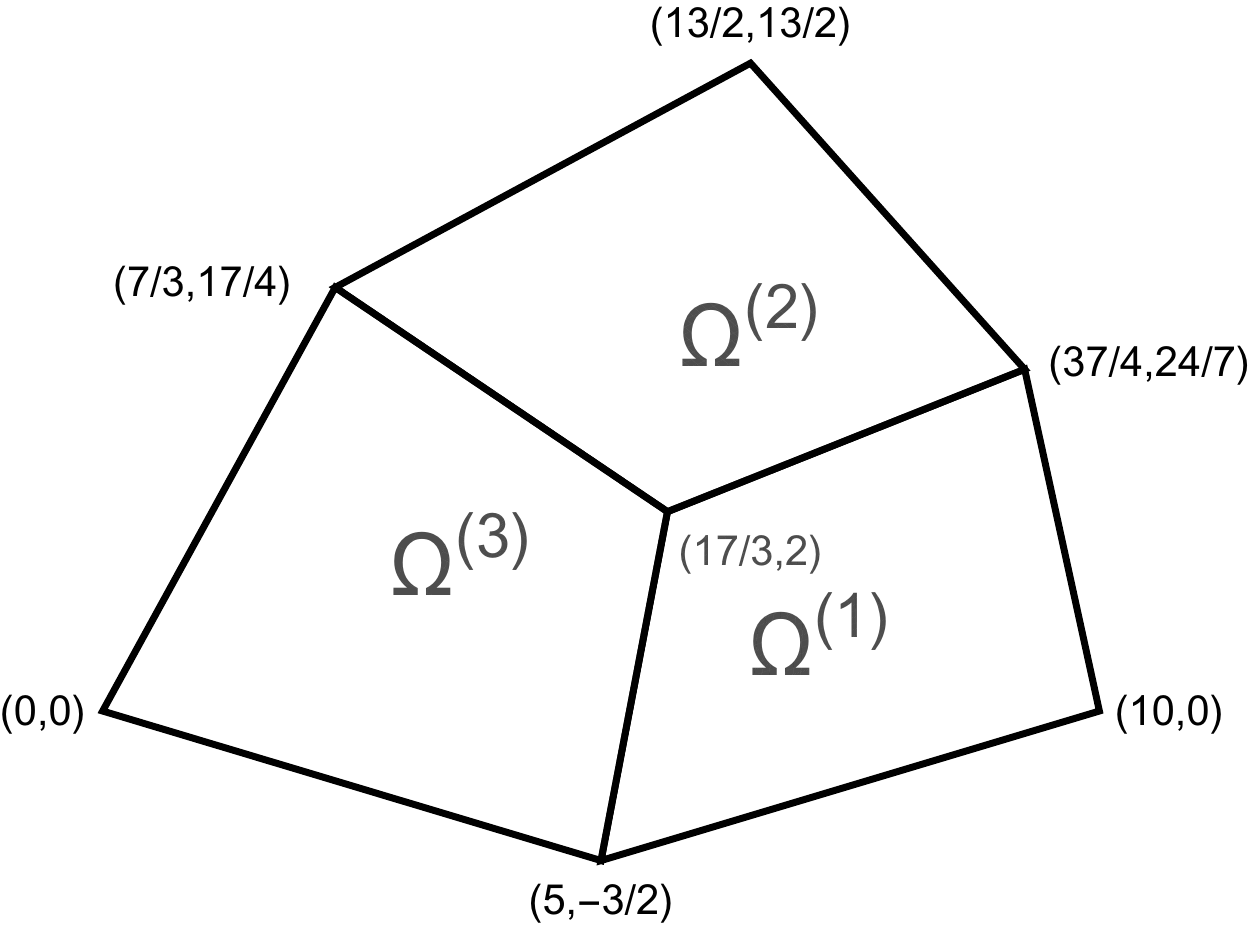} &
\includegraphics[width=3.7cm,clip]{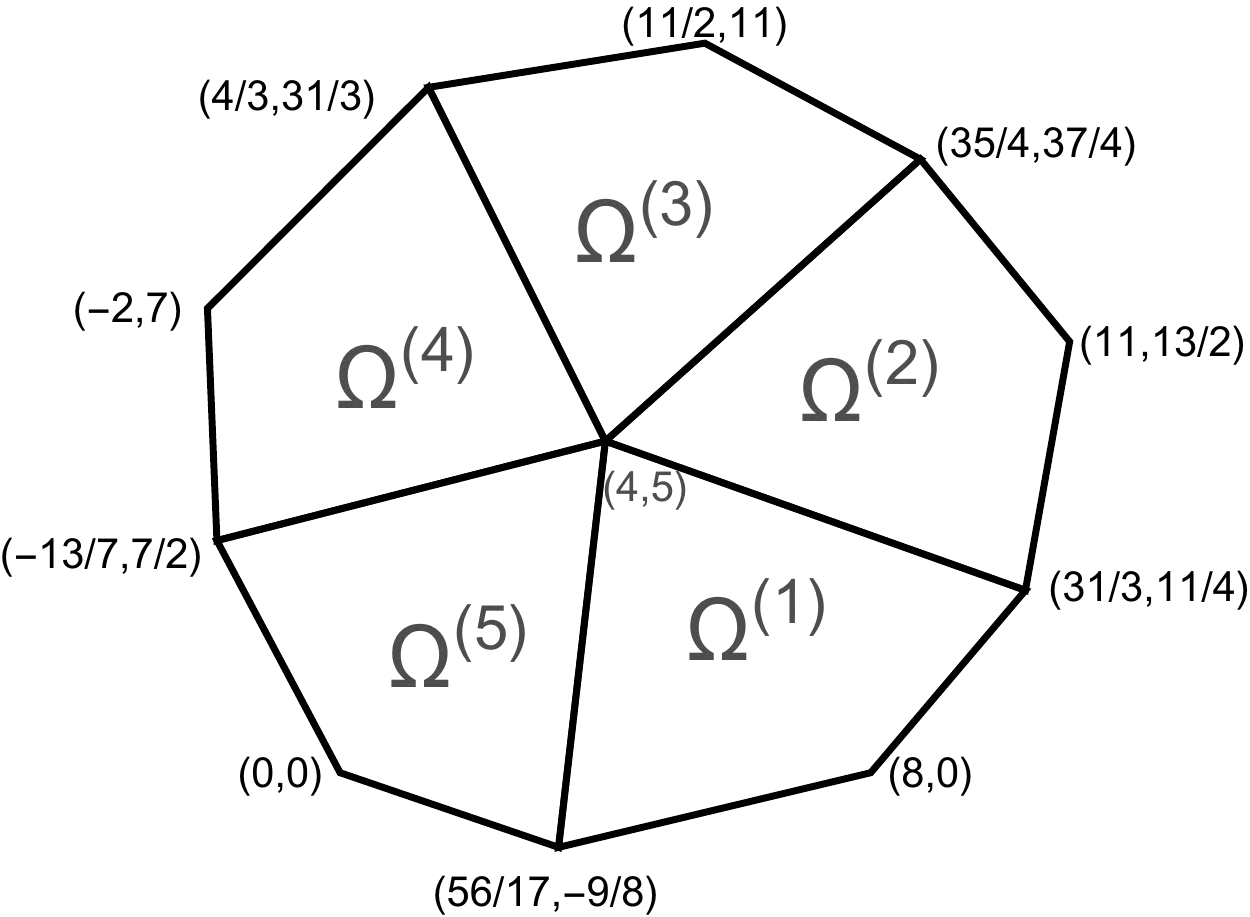} &
\includegraphics[width=3.4cm,clip]{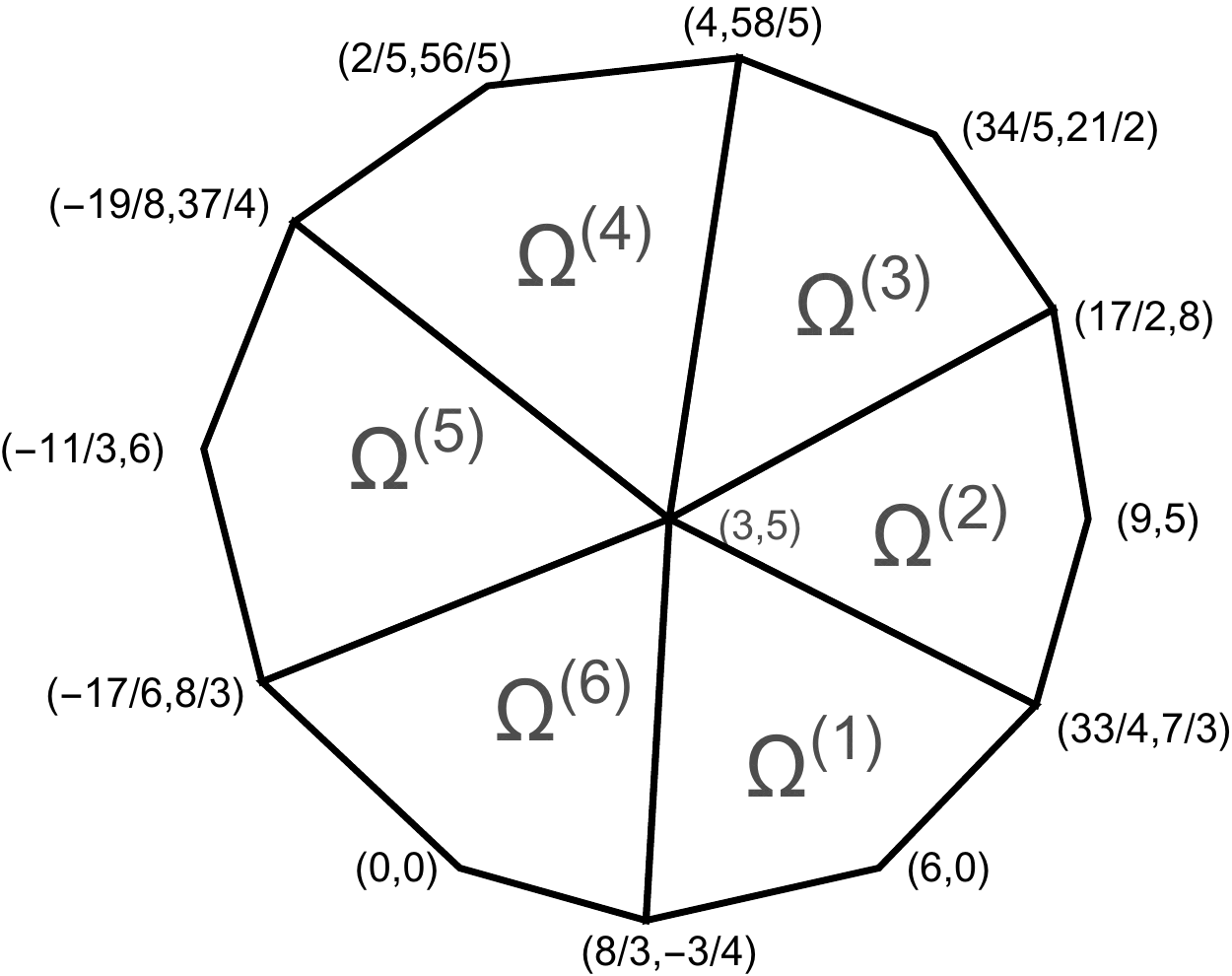} &
\includegraphics[width=3.2cm,clip]{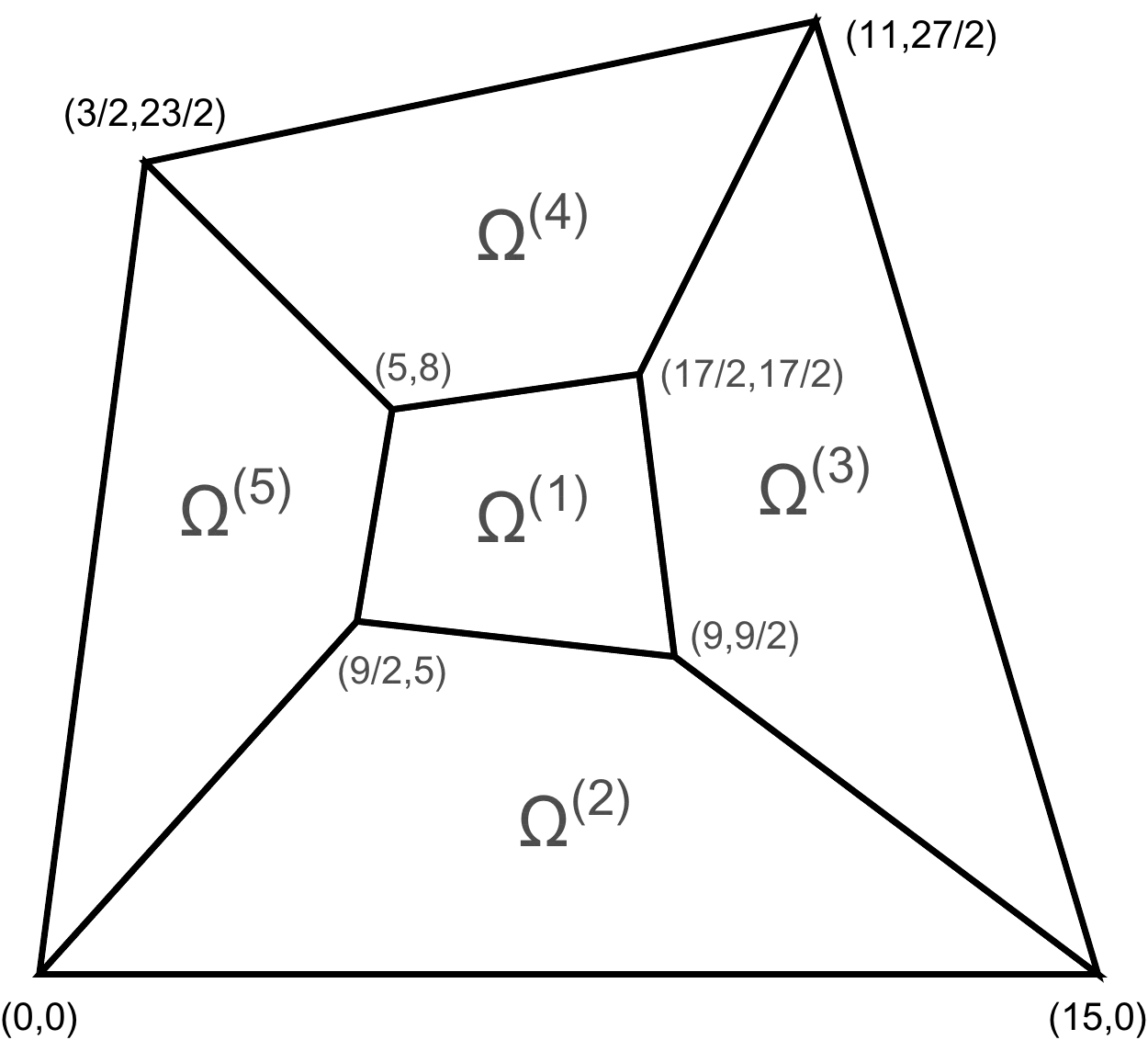} \\
\multicolumn{4}{c}{Computational domains} \\
\includegraphics[width=3.8cm,clip]{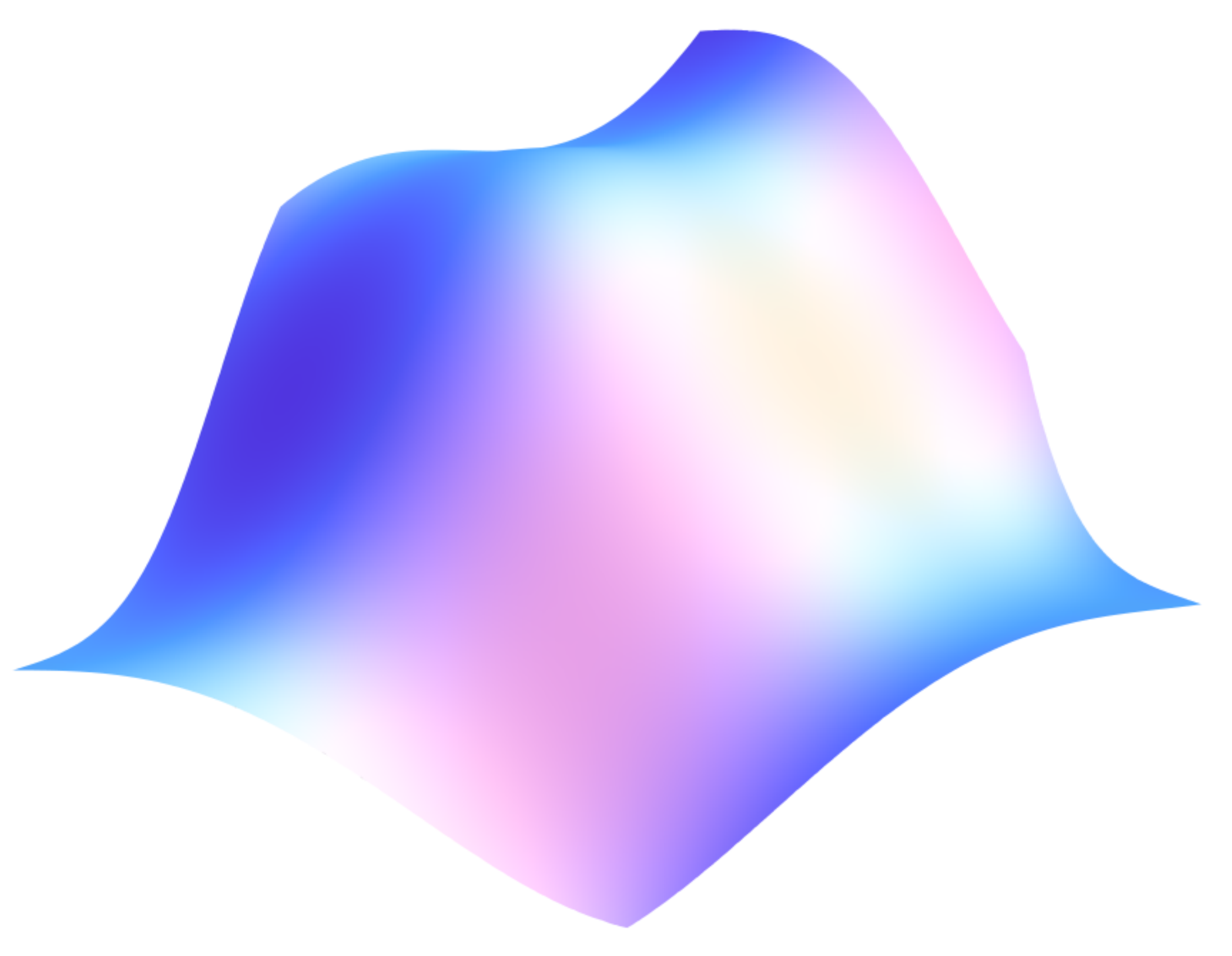} &
\includegraphics[width=3.7cm,clip]{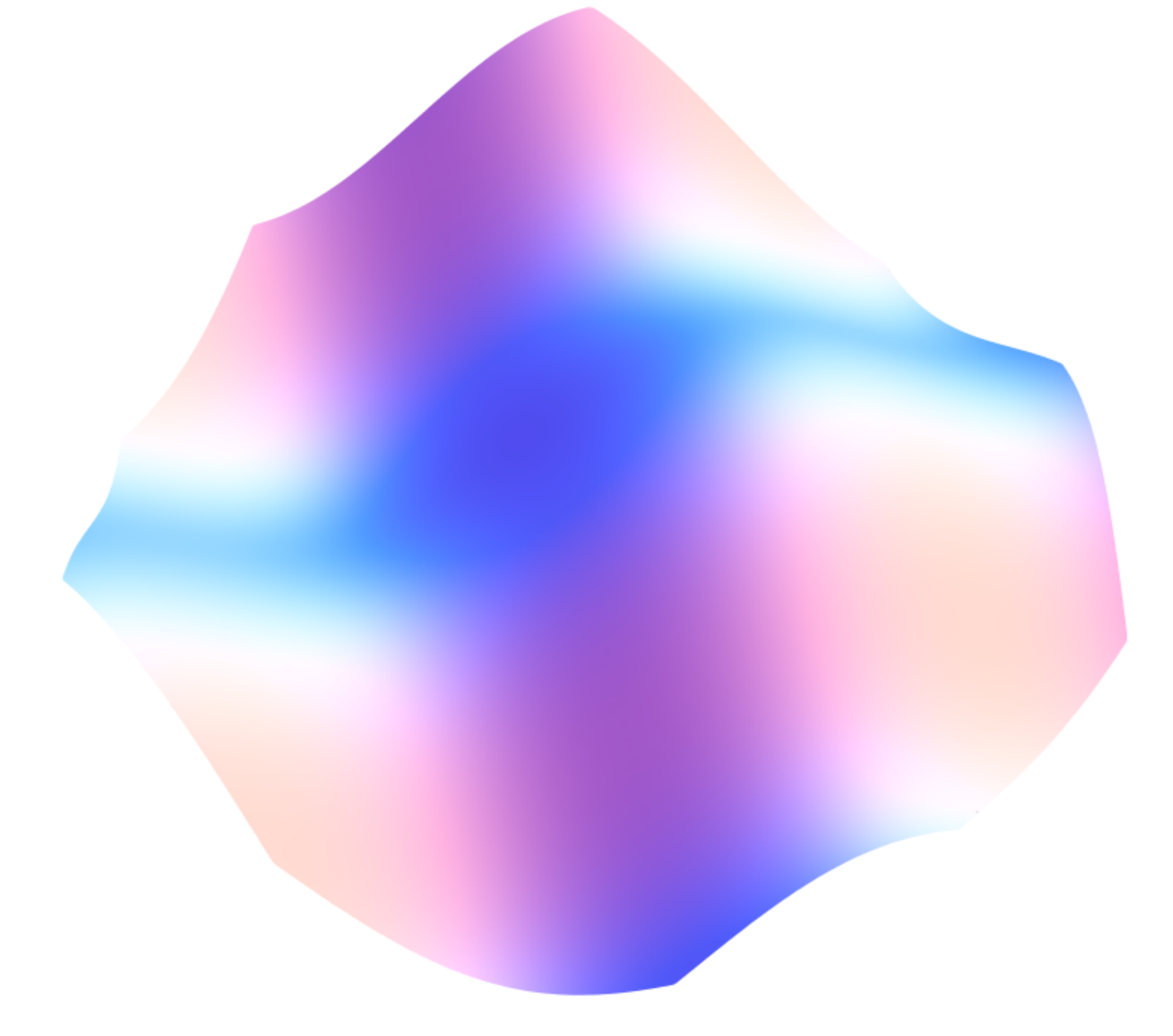} &
\includegraphics[width=3.2cm,clip]{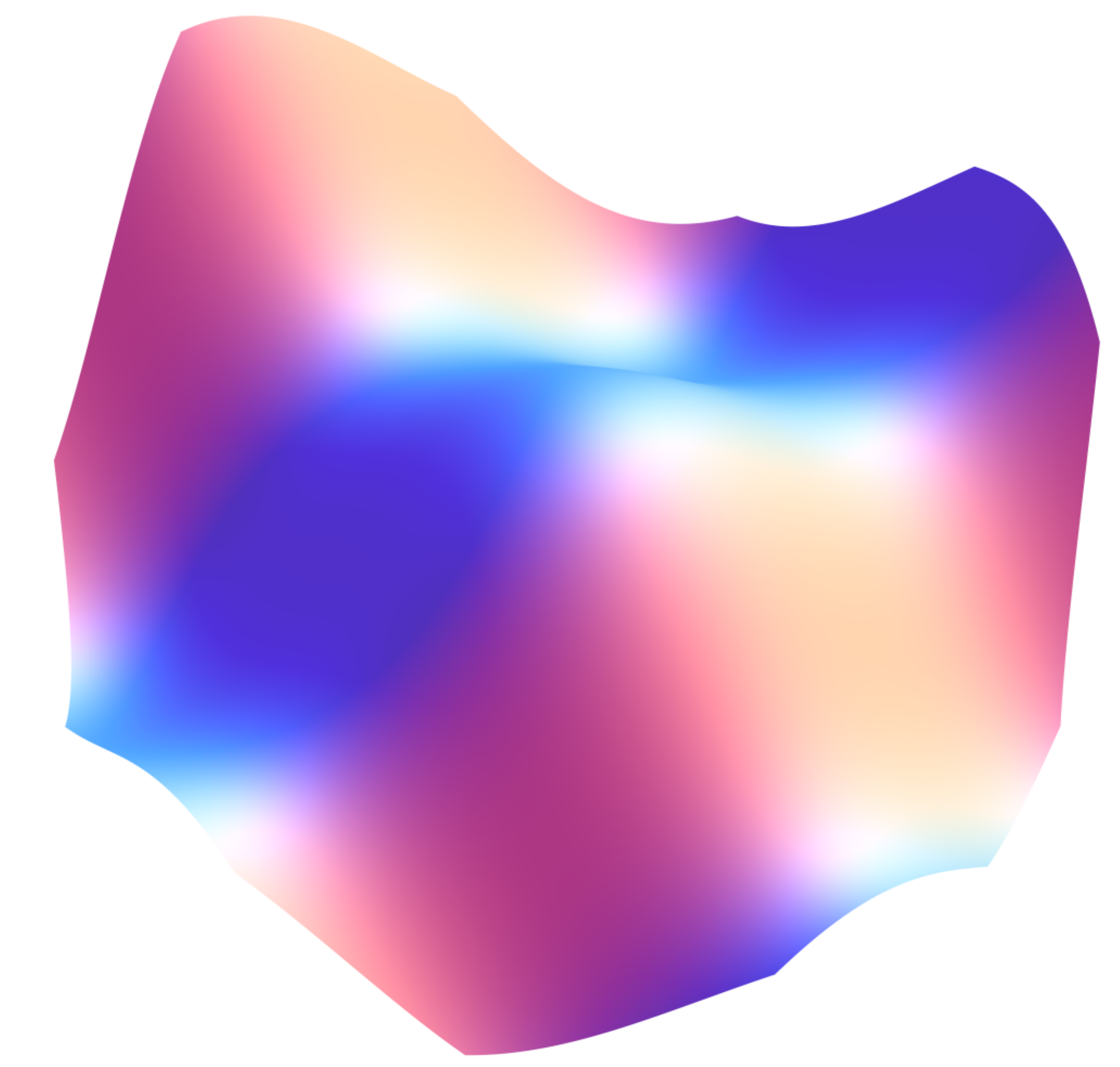} &
\includegraphics[width=3.6cm,clip]{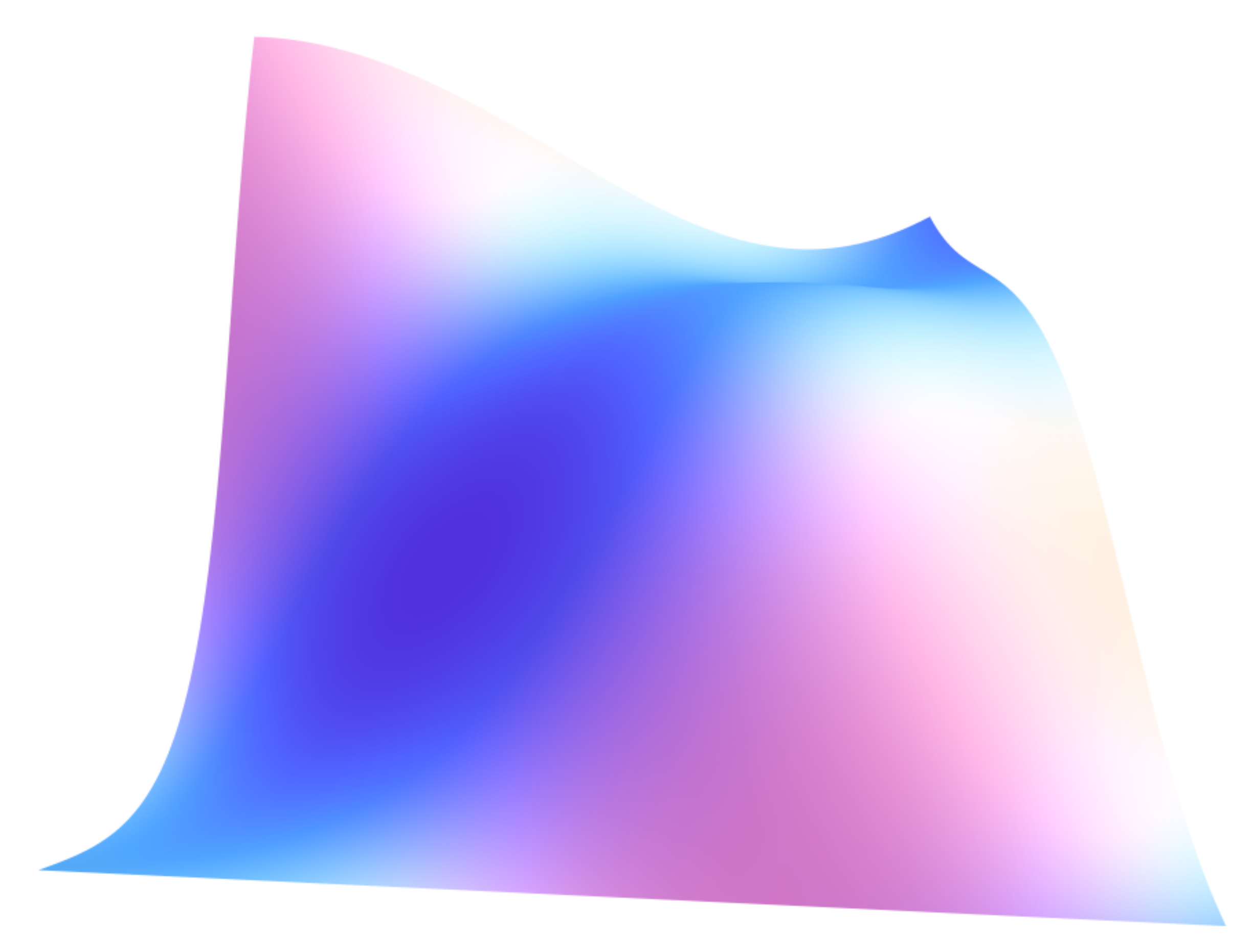} \\
\multicolumn{4}{c}{Exact solutions}
\end{tabular}
\caption{Computational domains and exact solutions for performing isogeometric collocation in Example~\ref{ex:multipatchdomain} (see also~Fig.~\ref{fig:results_Greville} 
and \ref{fig:results_Cauchy}).}
\label{fig:bilinear_domains}
\end{figure}

\begin{figure}[htp]
\centering\footnotesize
\begin{tabular}{ccc}
Relative $L^2$ error & Relative $H^1$ error & Relative $H^2$ error \\
\includegraphics[width=4.5cm,clip]{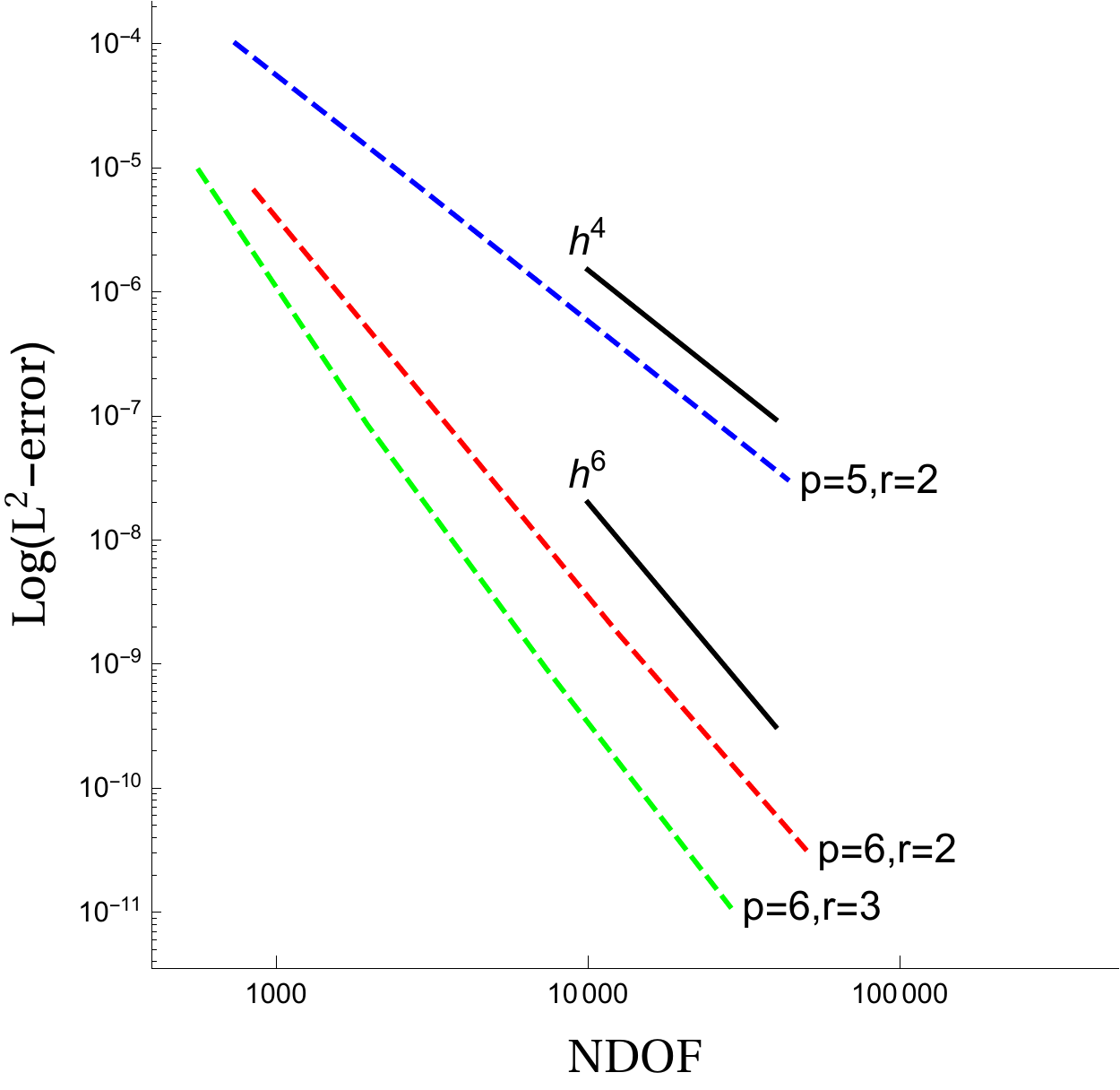} &
\includegraphics[width=4.5cm,clip]{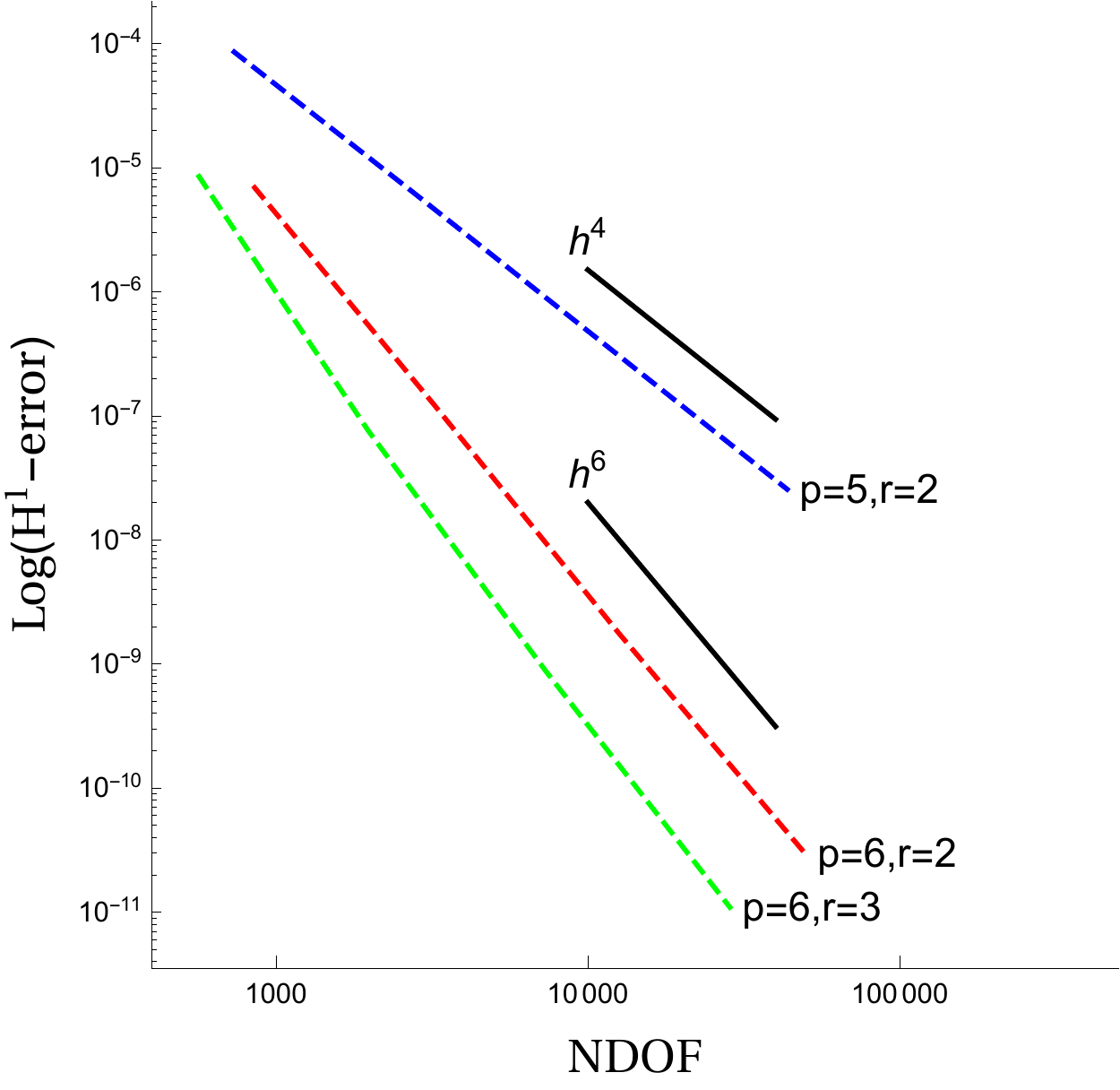} &
\includegraphics[width=4.5cm,clip]{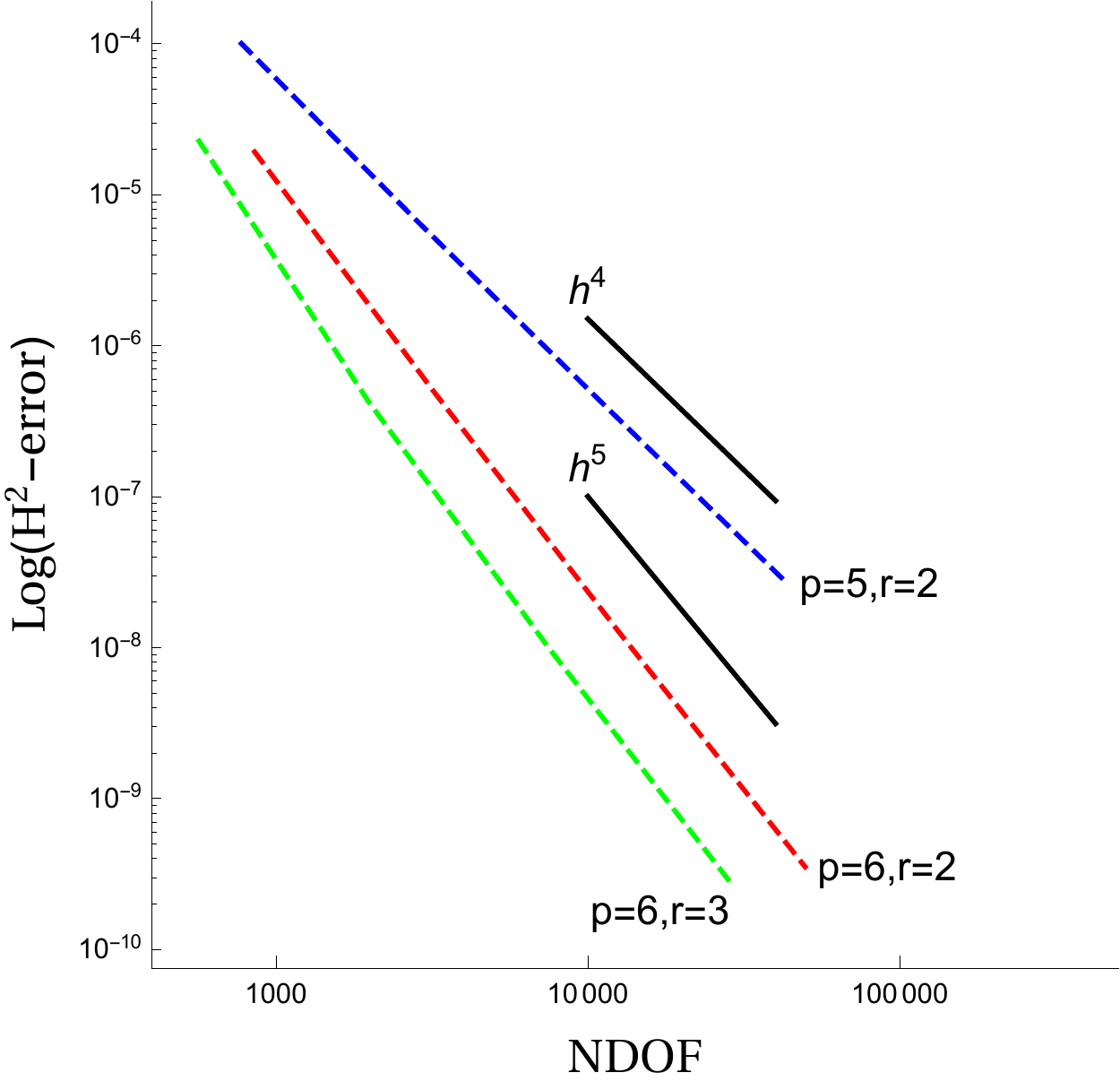} \\
\multicolumn{3}{c}{Three-patch domain~(a)} \\
\includegraphics[width=4.5cm,clip]{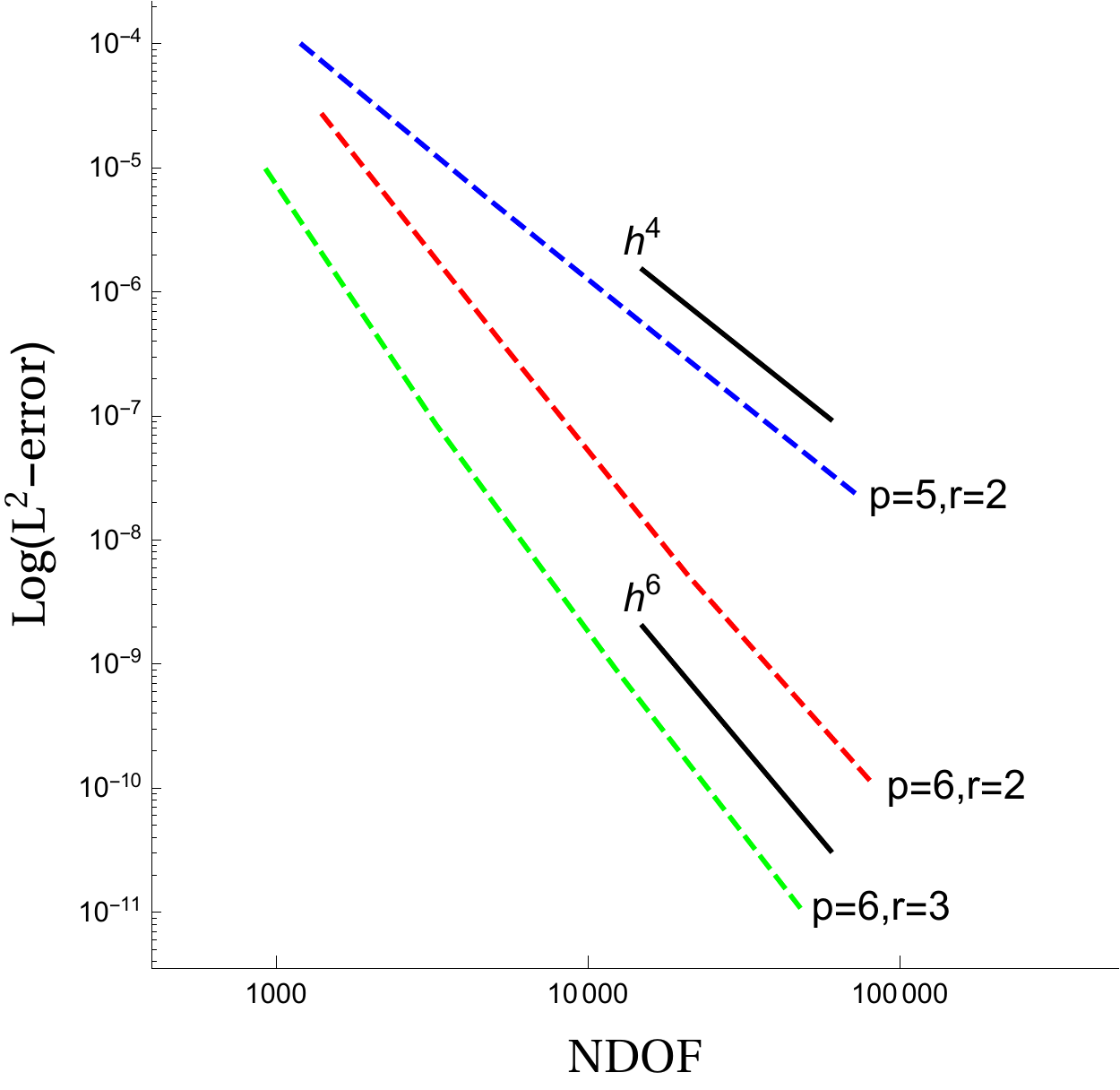} &
\includegraphics[width=4.5cm,clip]{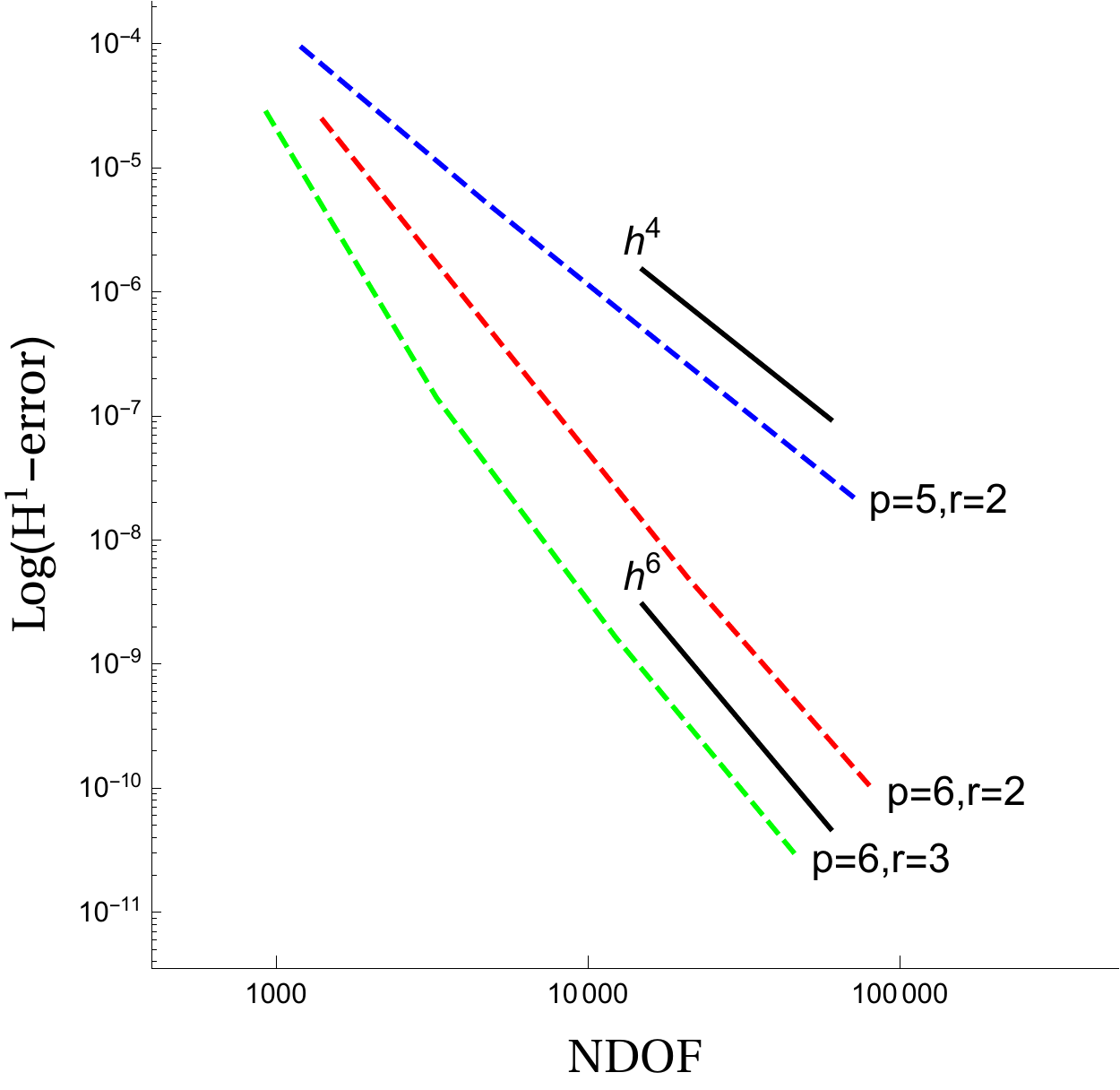} &
\includegraphics[width=4.5cm,clip]{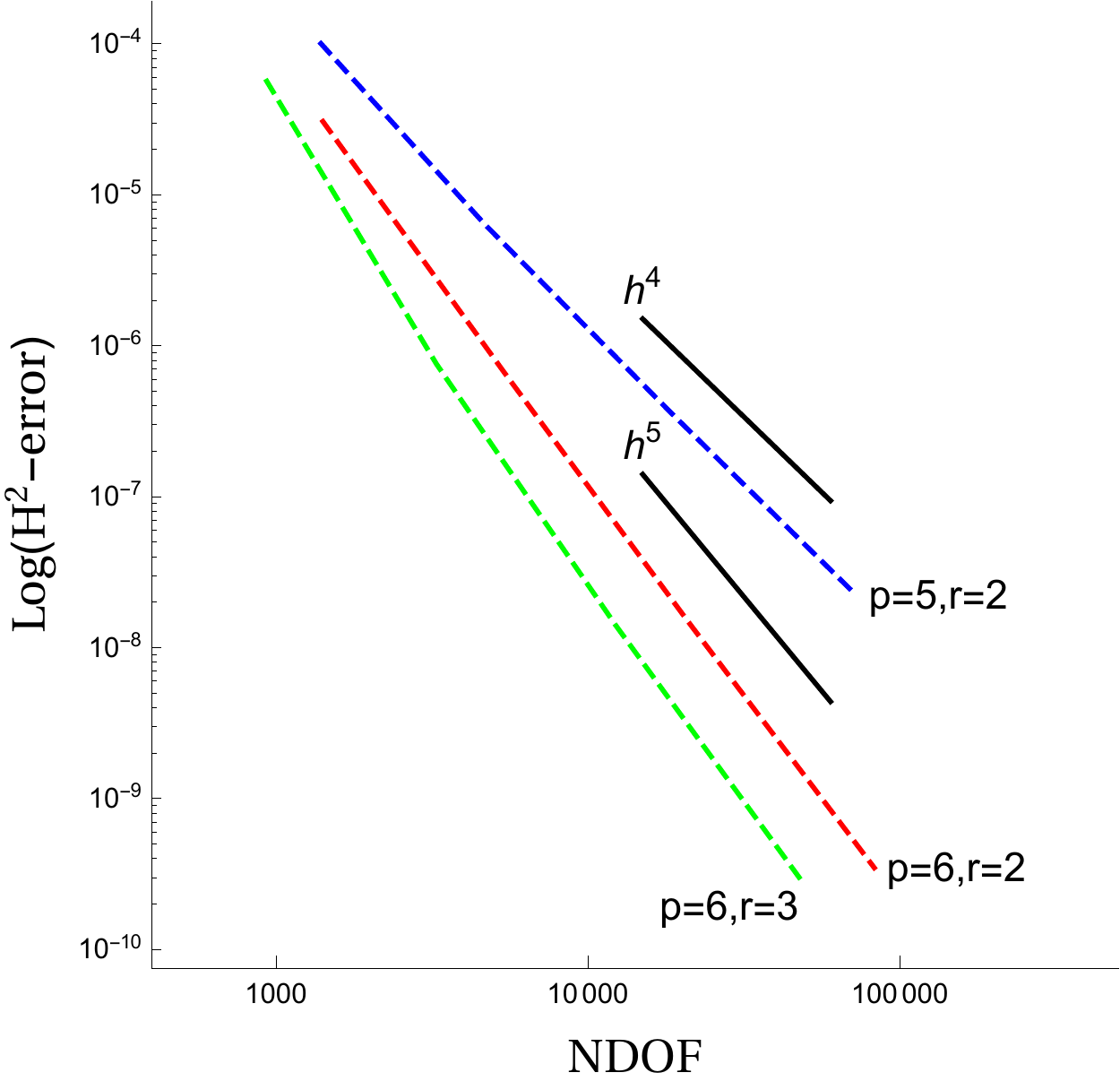} \\
\multicolumn{3}{c}{Five-patch domain~(b)} \\
\includegraphics[width=4.5cm,clip]{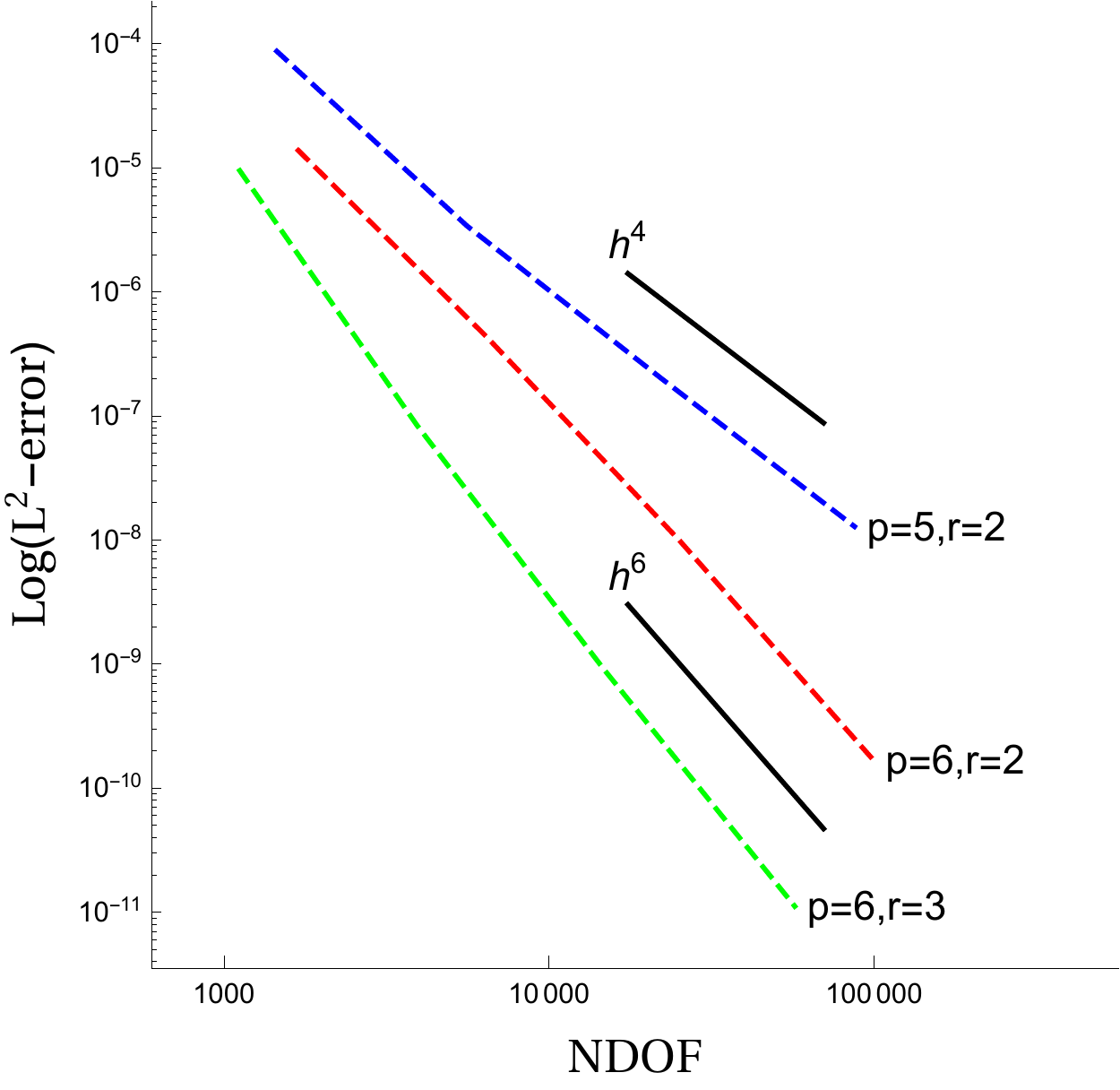} &
\includegraphics[width=4.5cm,clip]{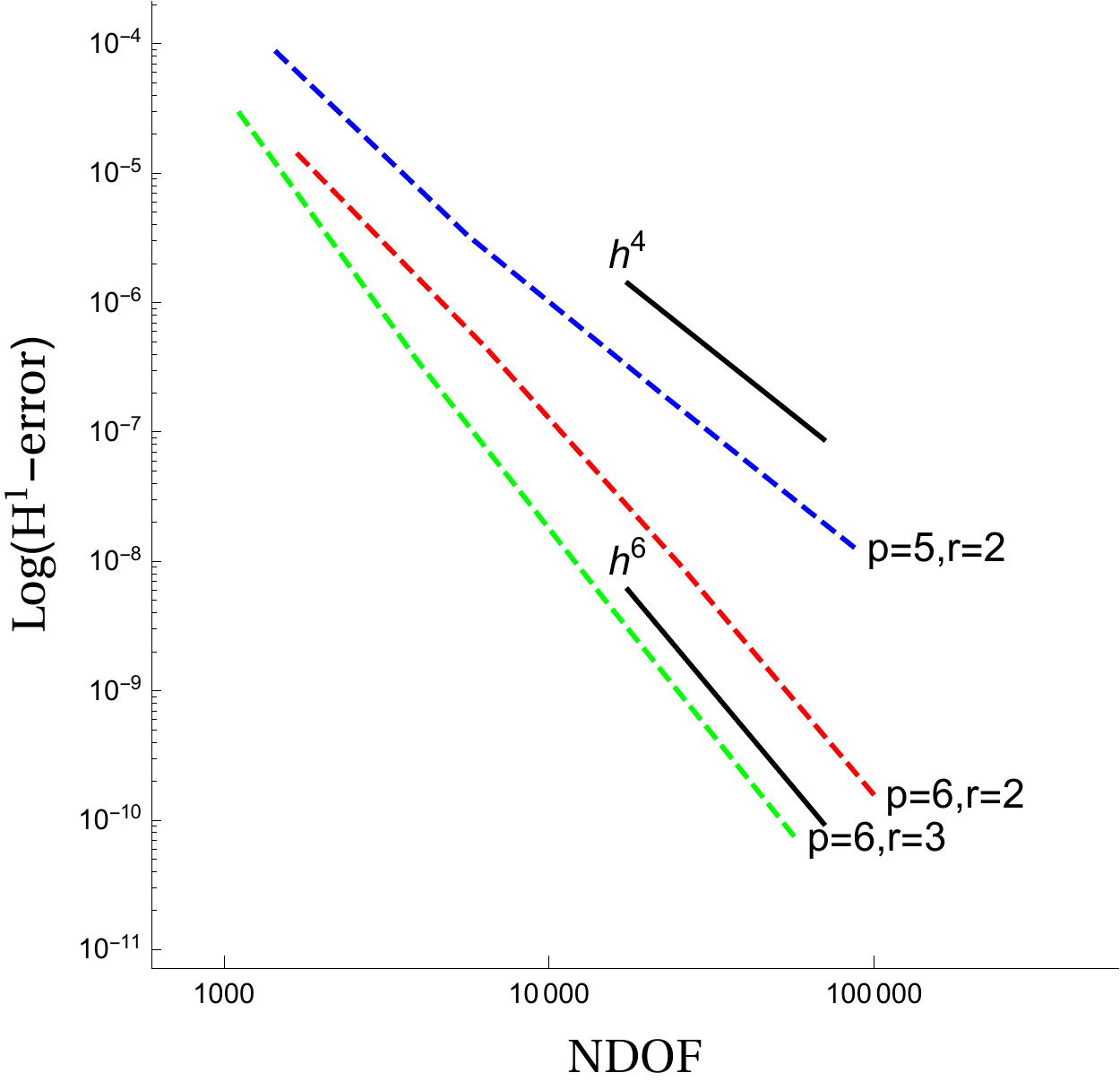} &
\includegraphics[width=4.5cm,clip]{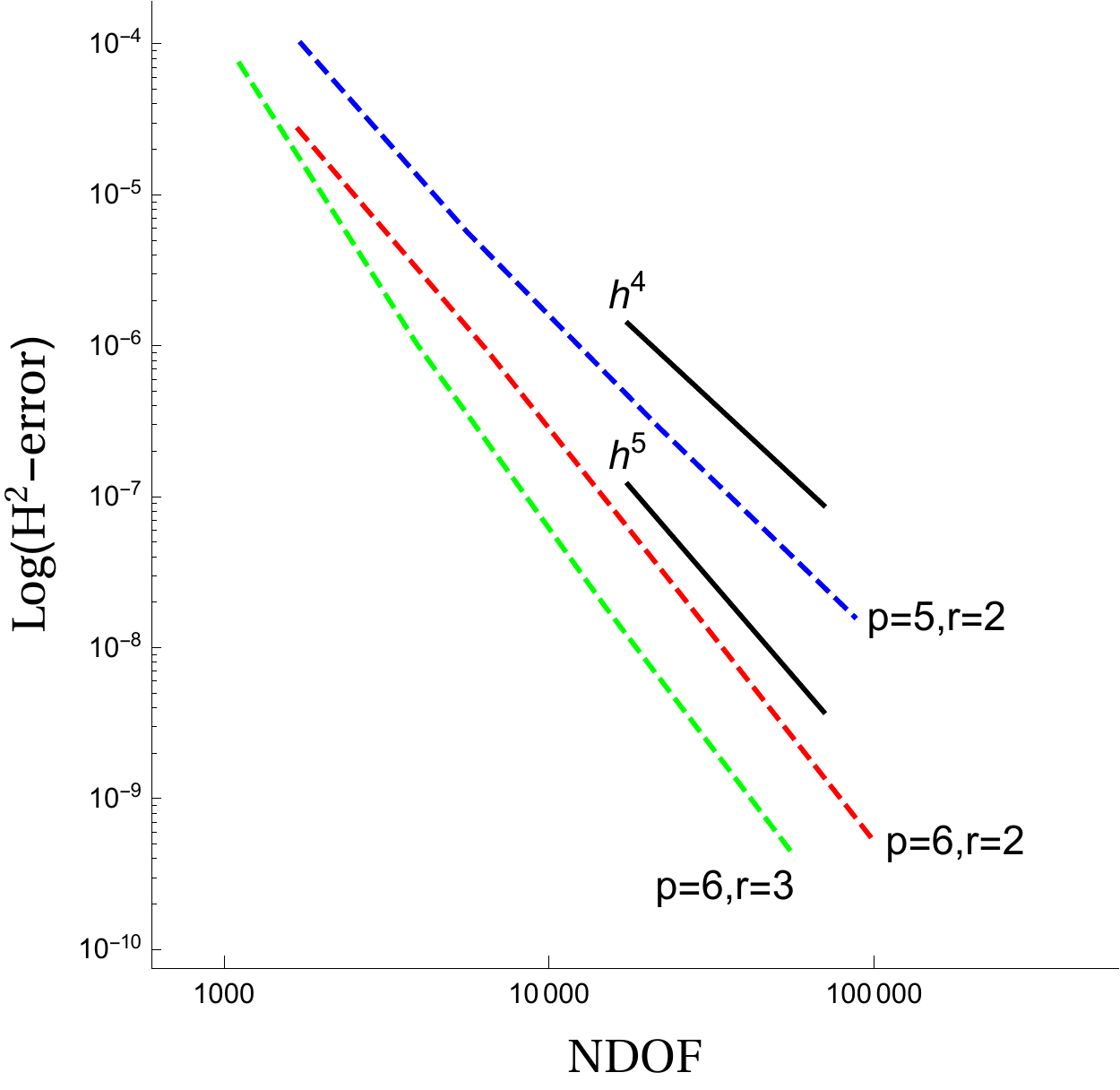} \\
\multicolumn{3}{c}{Six-patch domain~(c)} \\
\includegraphics[width=4.5cm,clip]{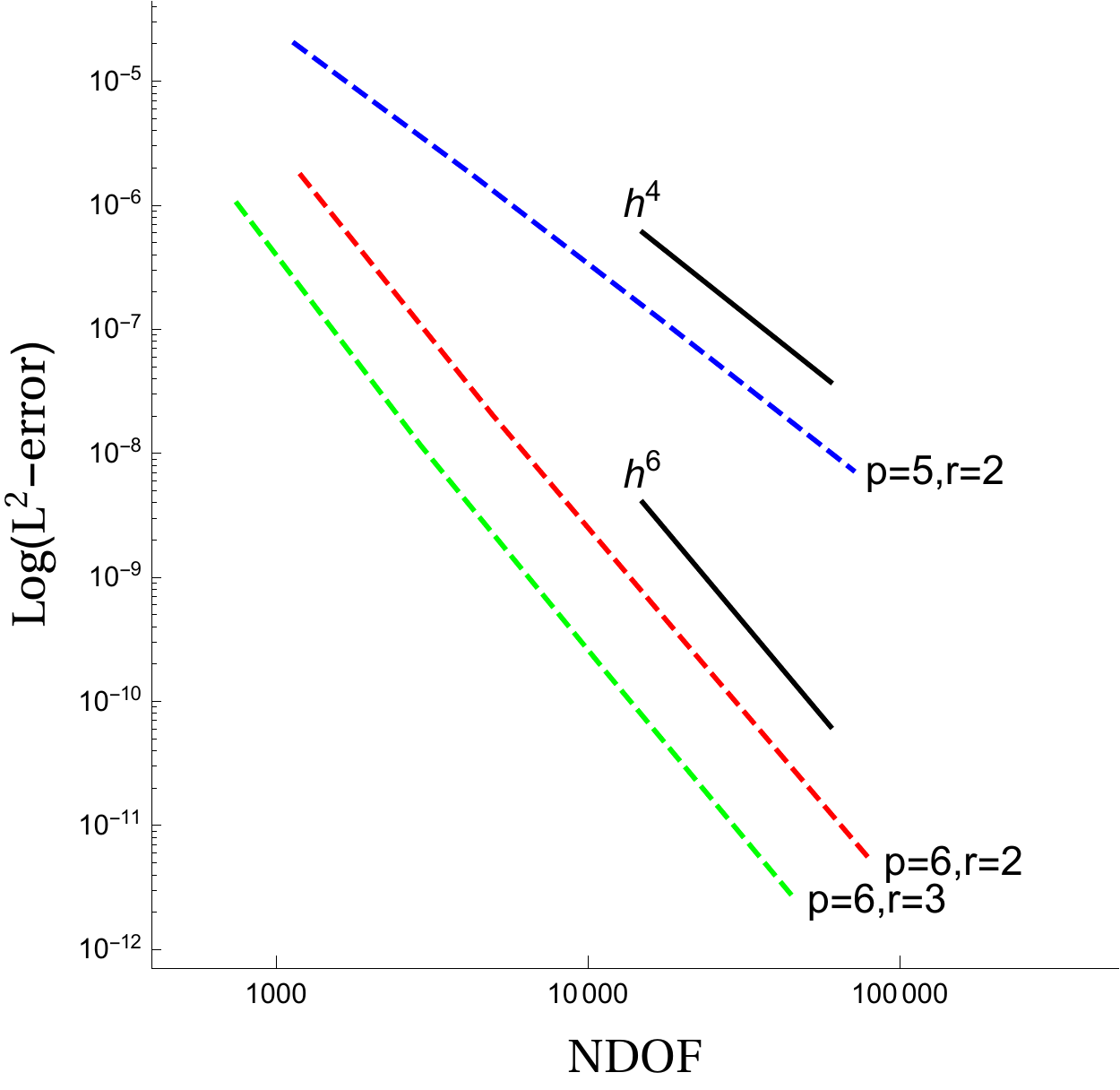} &
\includegraphics[width=4.5cm,clip]{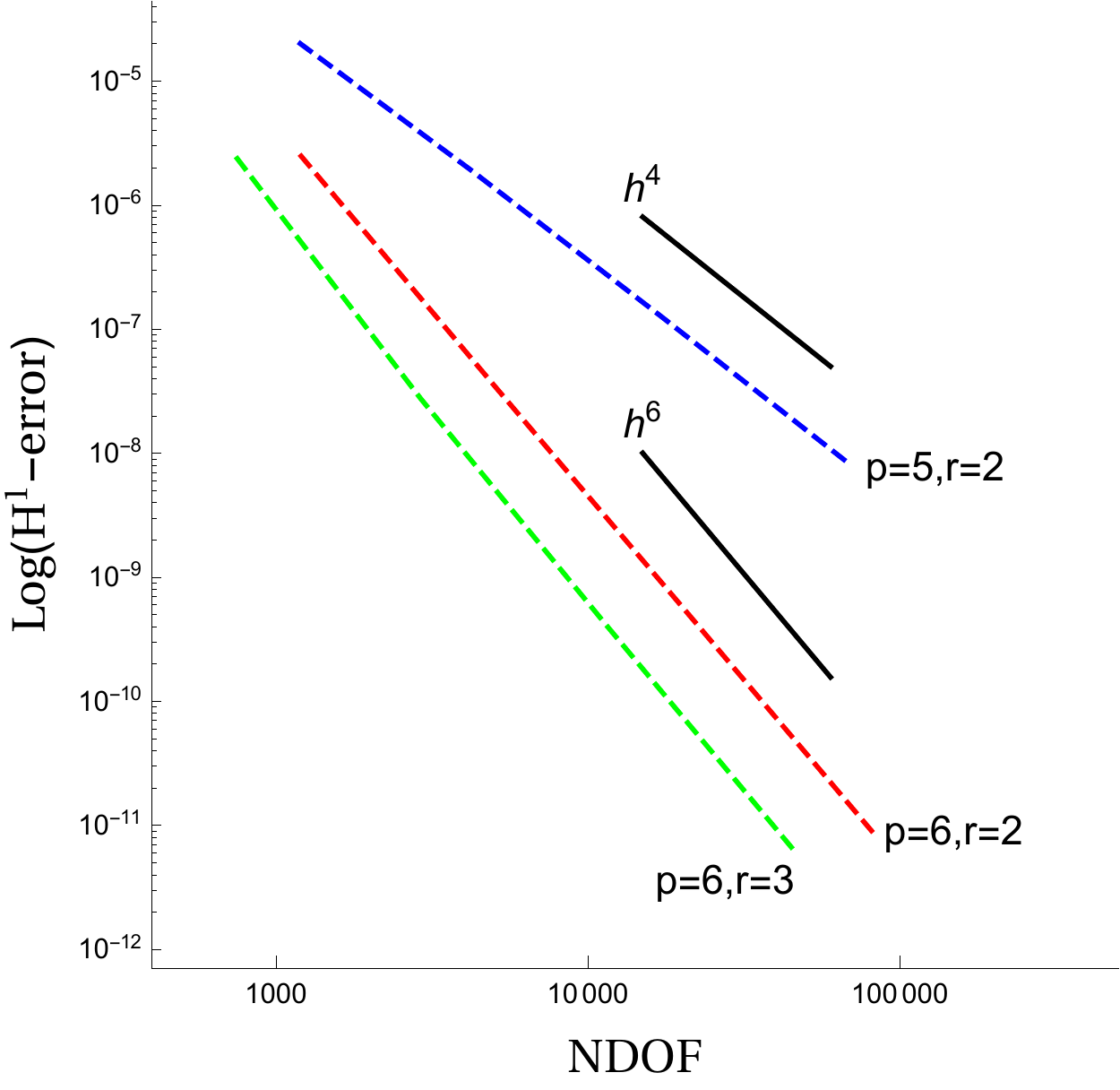} &
\includegraphics[width=4.5cm,clip]{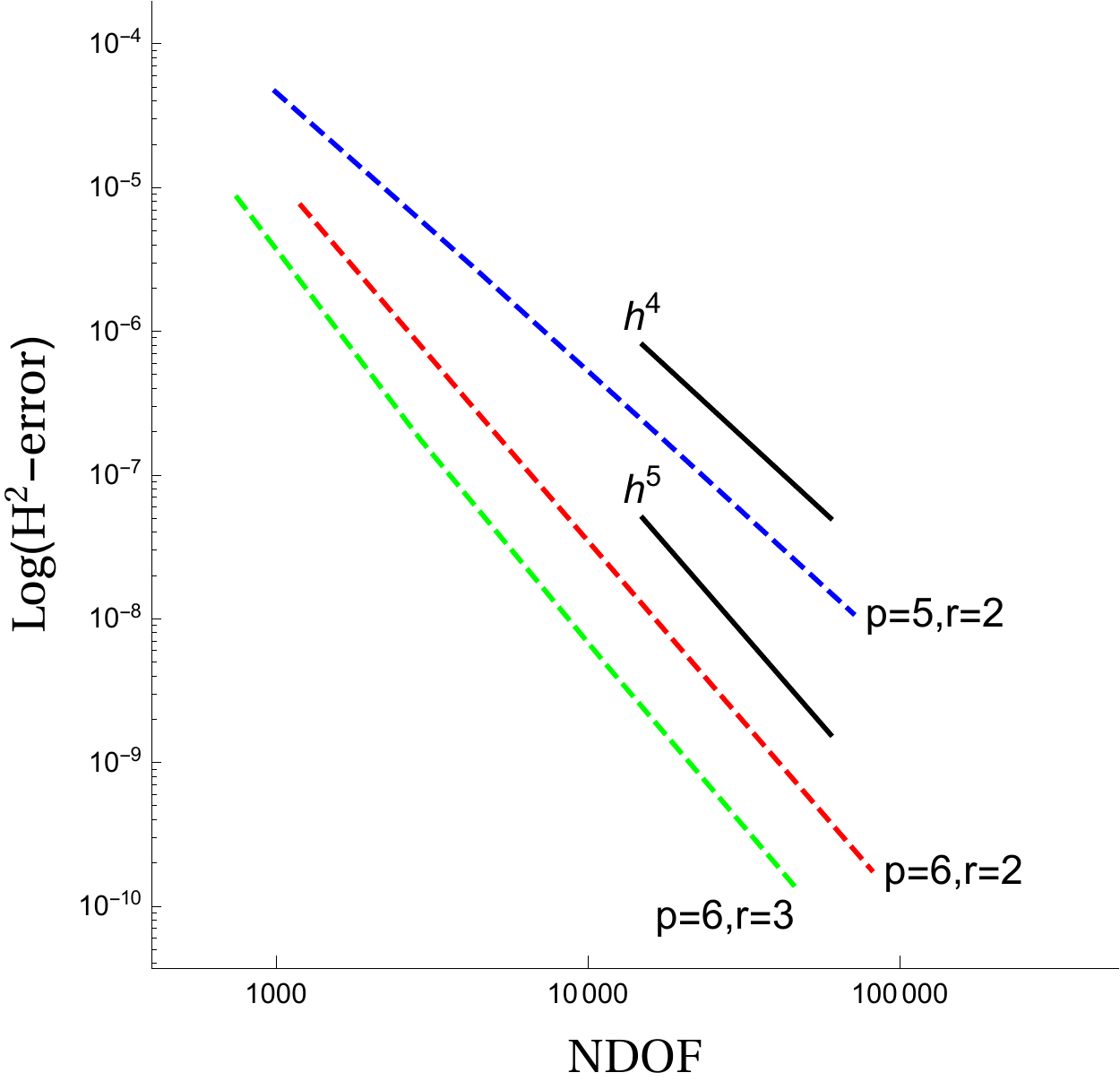} \\
\multicolumn{3}{c}{Five-patch domain~(d)} \\
\end{tabular}
\caption{Error plots w.r.t. the number of degrees of freedom (NDOF) of performing isogeometric collocation on the different multi-patch domains given in 
Fig.~\ref{fig:bilinear_domains}~(first row) for right side functions obtained by the exact solutions shown in Fig~\ref{fig:bilinear_domains}~(second row) using the 
Greville points as collocation points.}
\label{fig:results_Greville}
\end{figure}

\begin{figure}[htp]
\centering\footnotesize
\begin{tabular}{ccc}
Relative $L^2$ error & Relative $H^1$ error & Relative $H^2$ error \\
\includegraphics[width=4.5cm,clip]{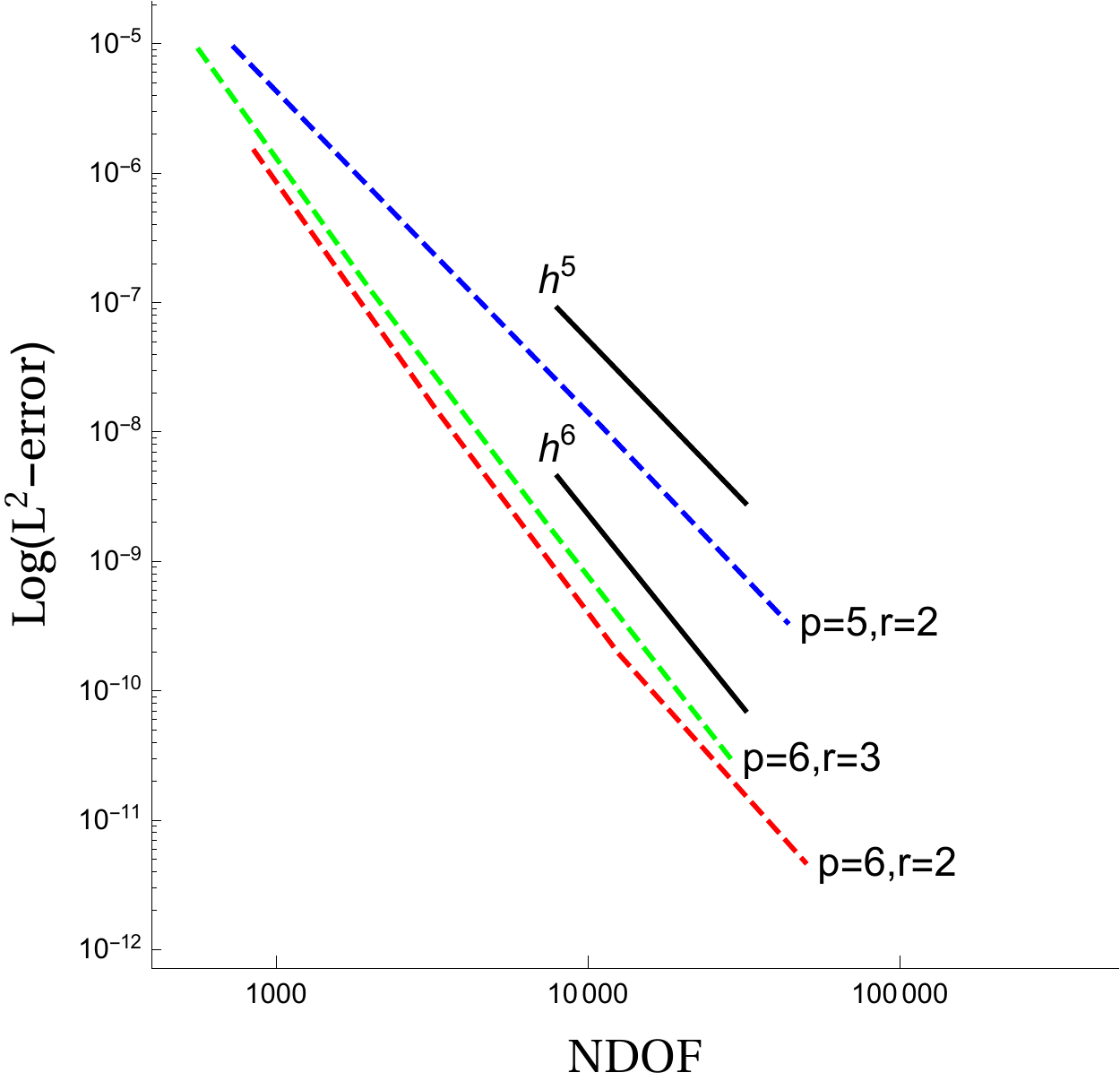} &
\includegraphics[width=4.5cm,clip]{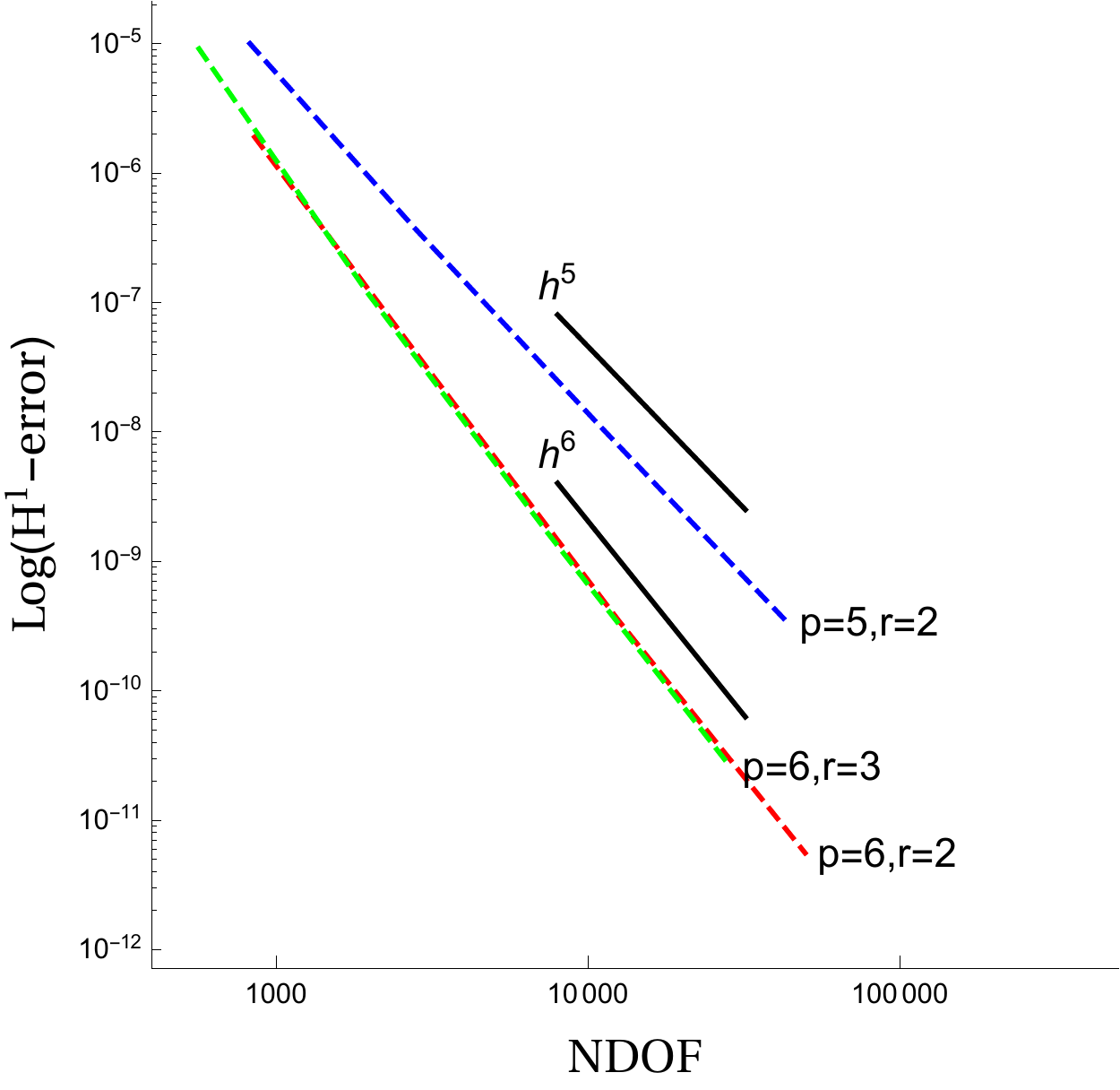} &
\includegraphics[width=4.5cm,clip]{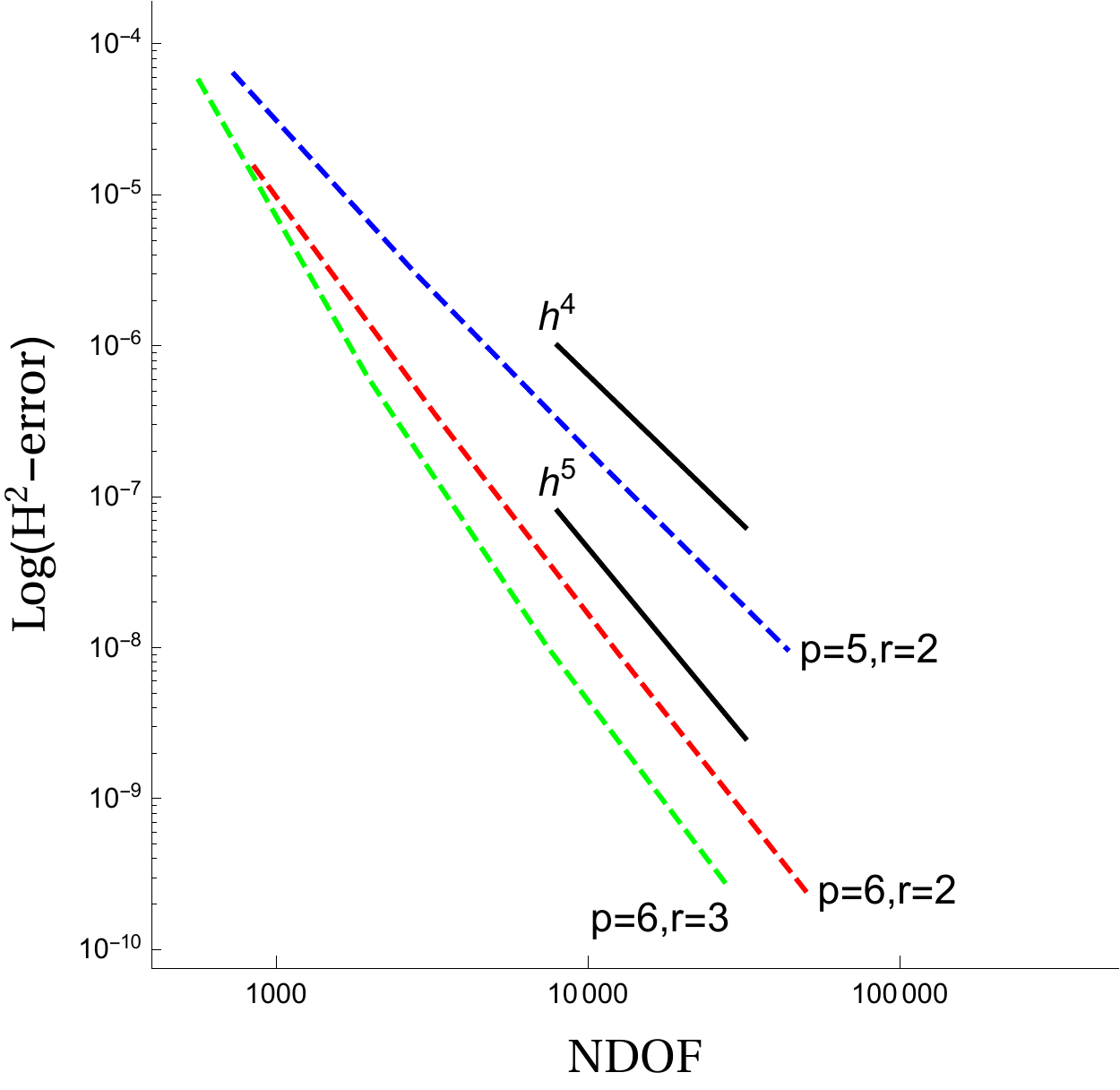} \\
\multicolumn{3}{c}{Three-patch domain~(a)} \\
\includegraphics[width=4.5cm,clip]{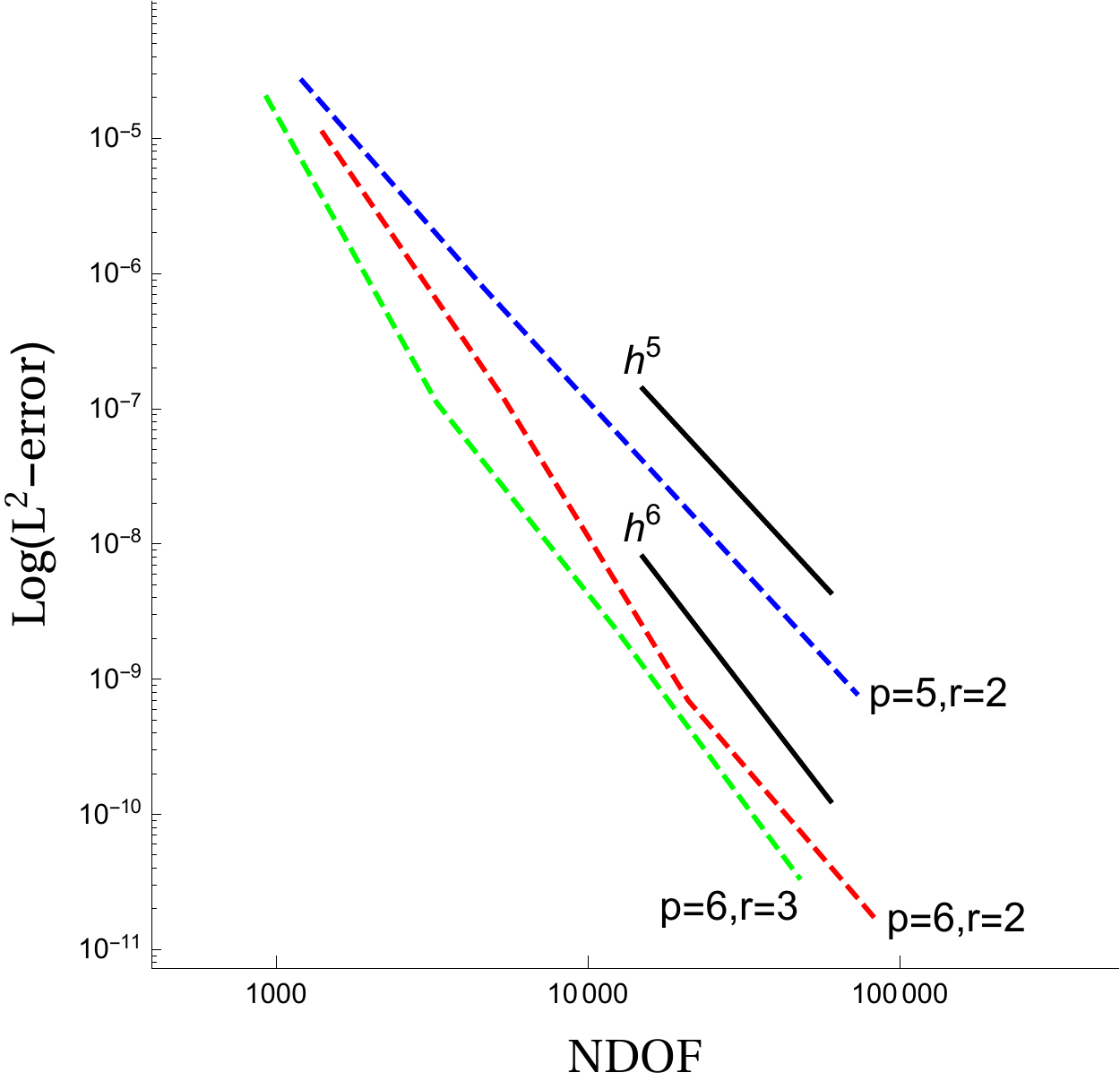} &
\includegraphics[width=4.5cm,clip]{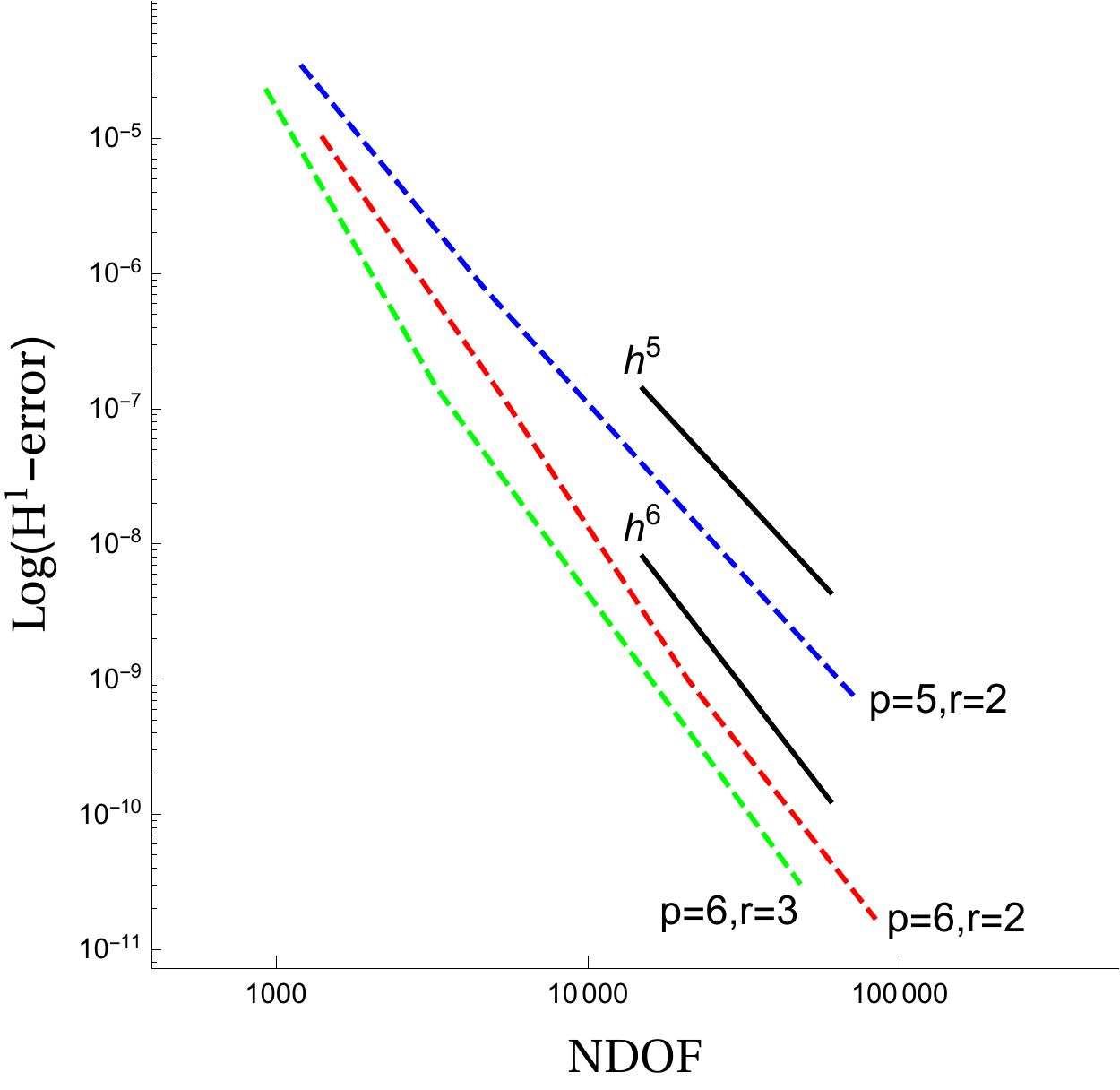} &
\includegraphics[width=4.5cm,clip]{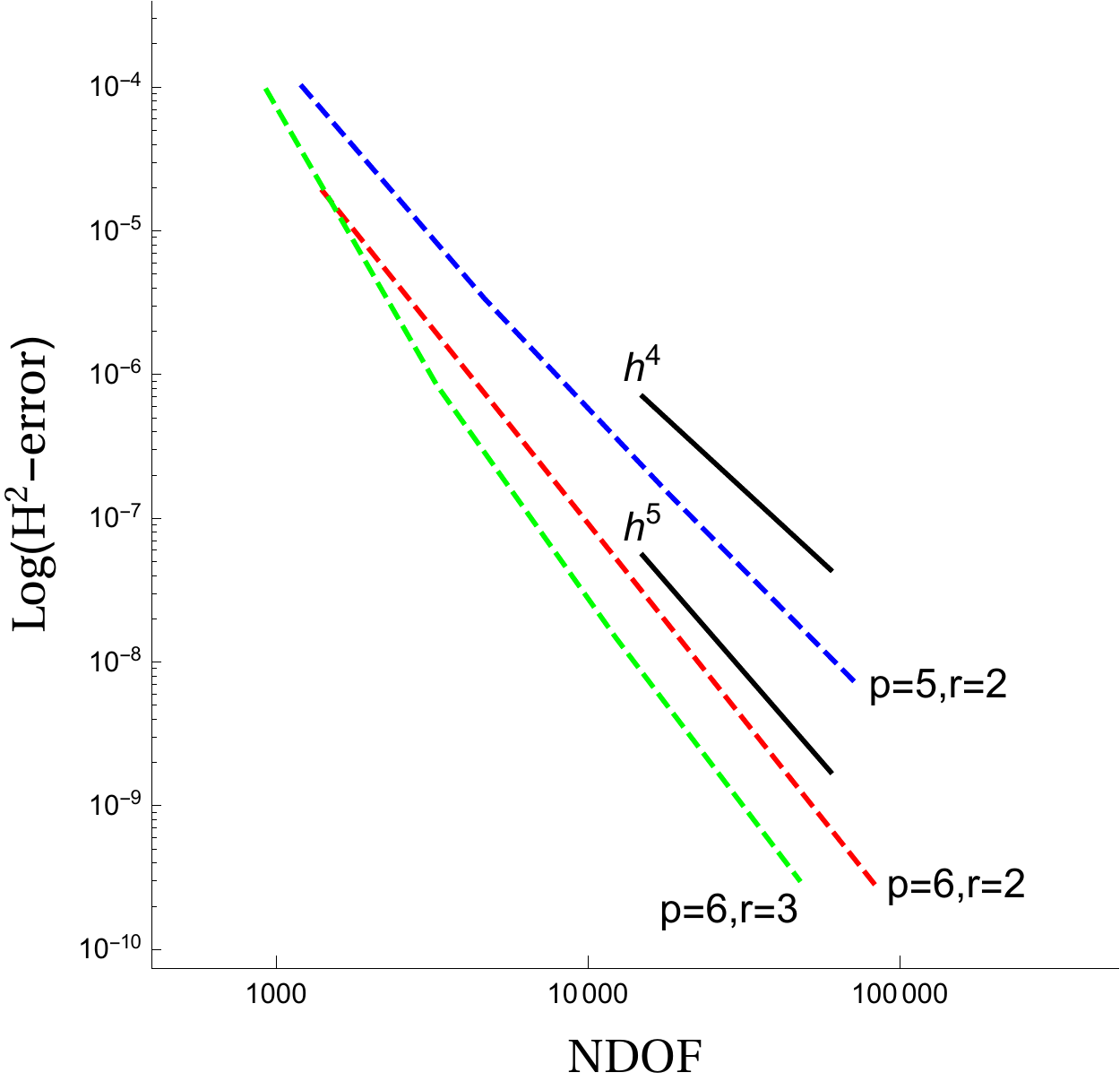} \\
\multicolumn{3}{c}{Five-patch domain~(b)} \\
\includegraphics[width=4.5cm,clip]{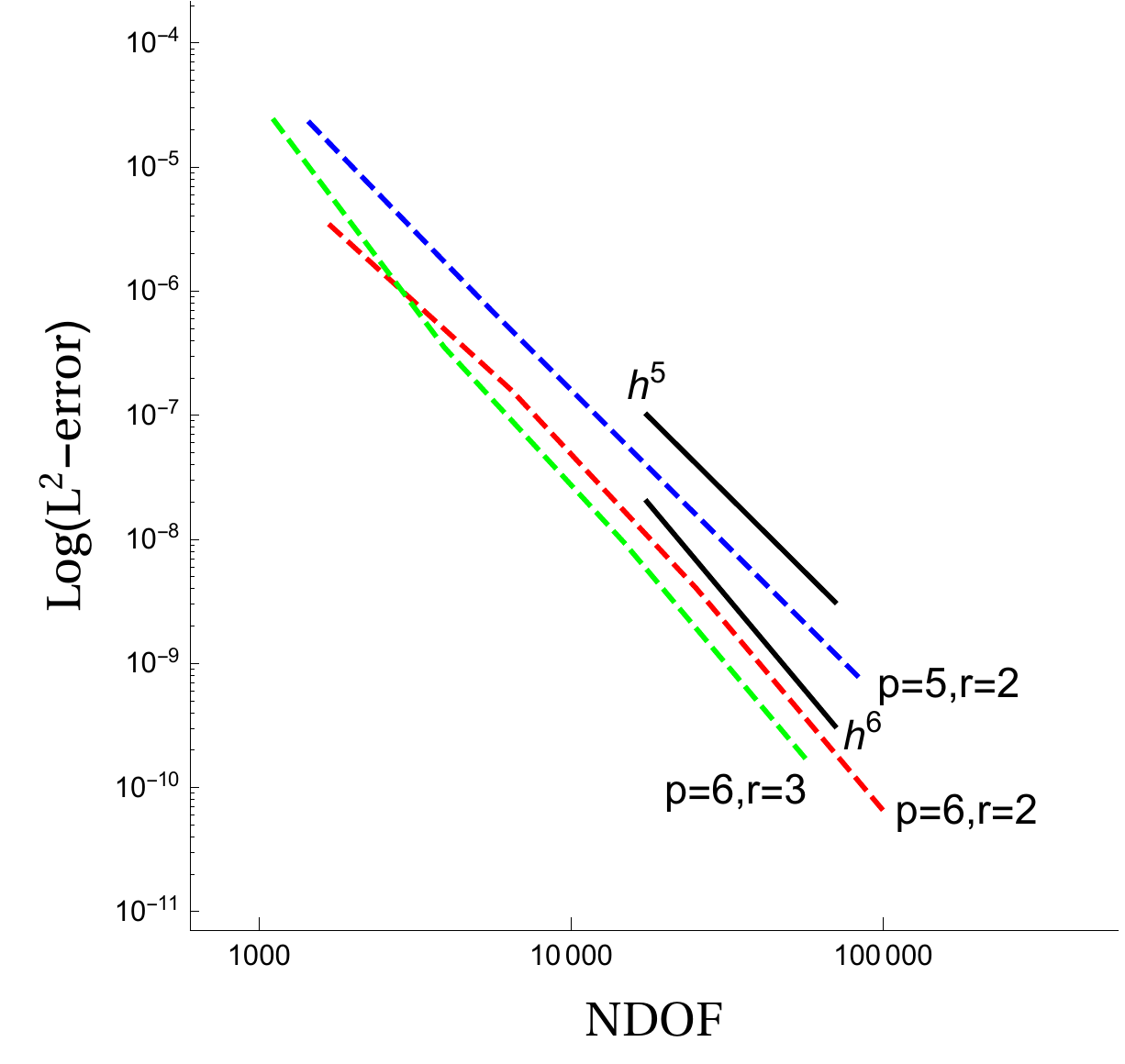} &
\includegraphics[width=4.5cm,clip]{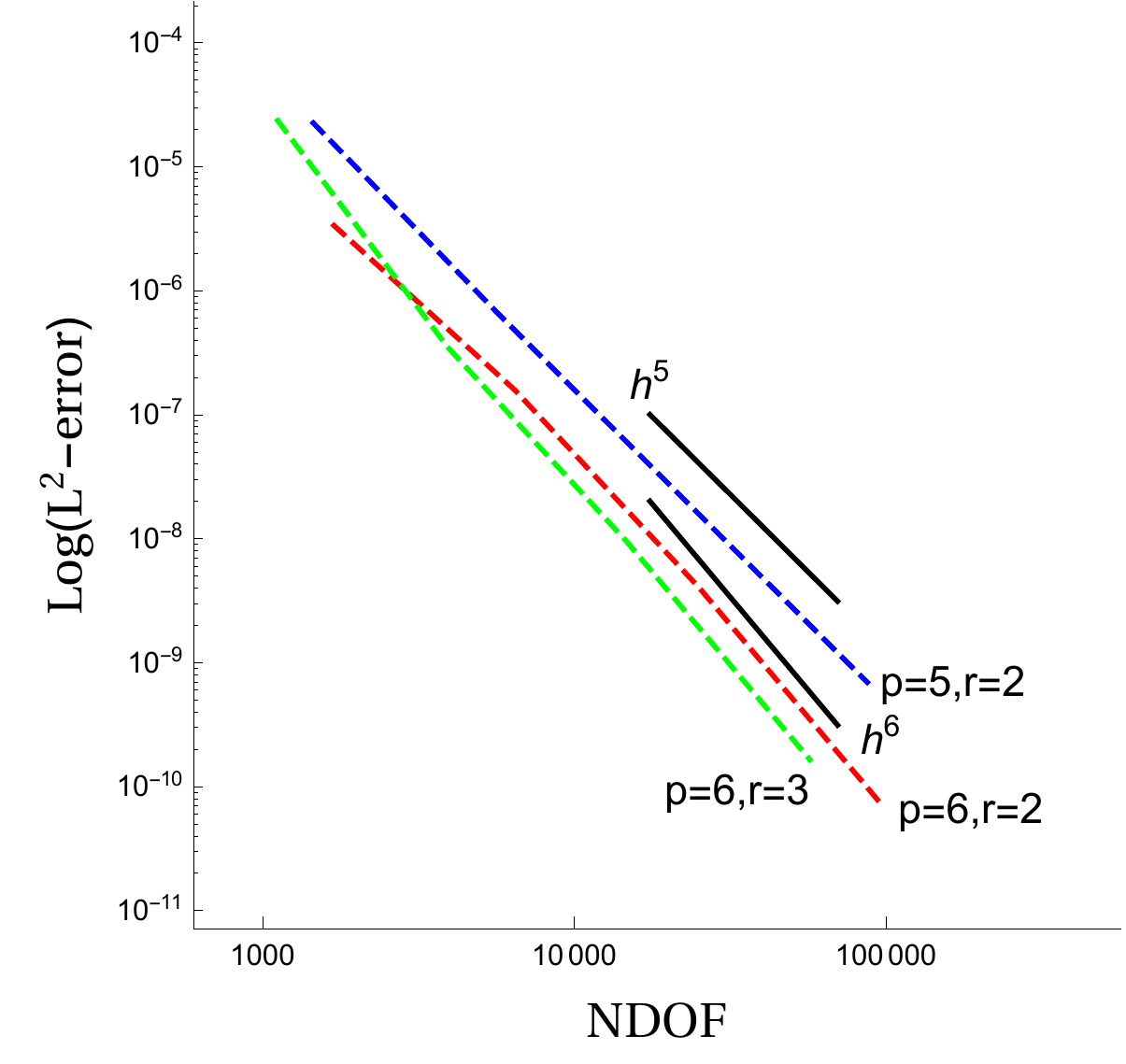} &
\includegraphics[width=4.5cm,clip]{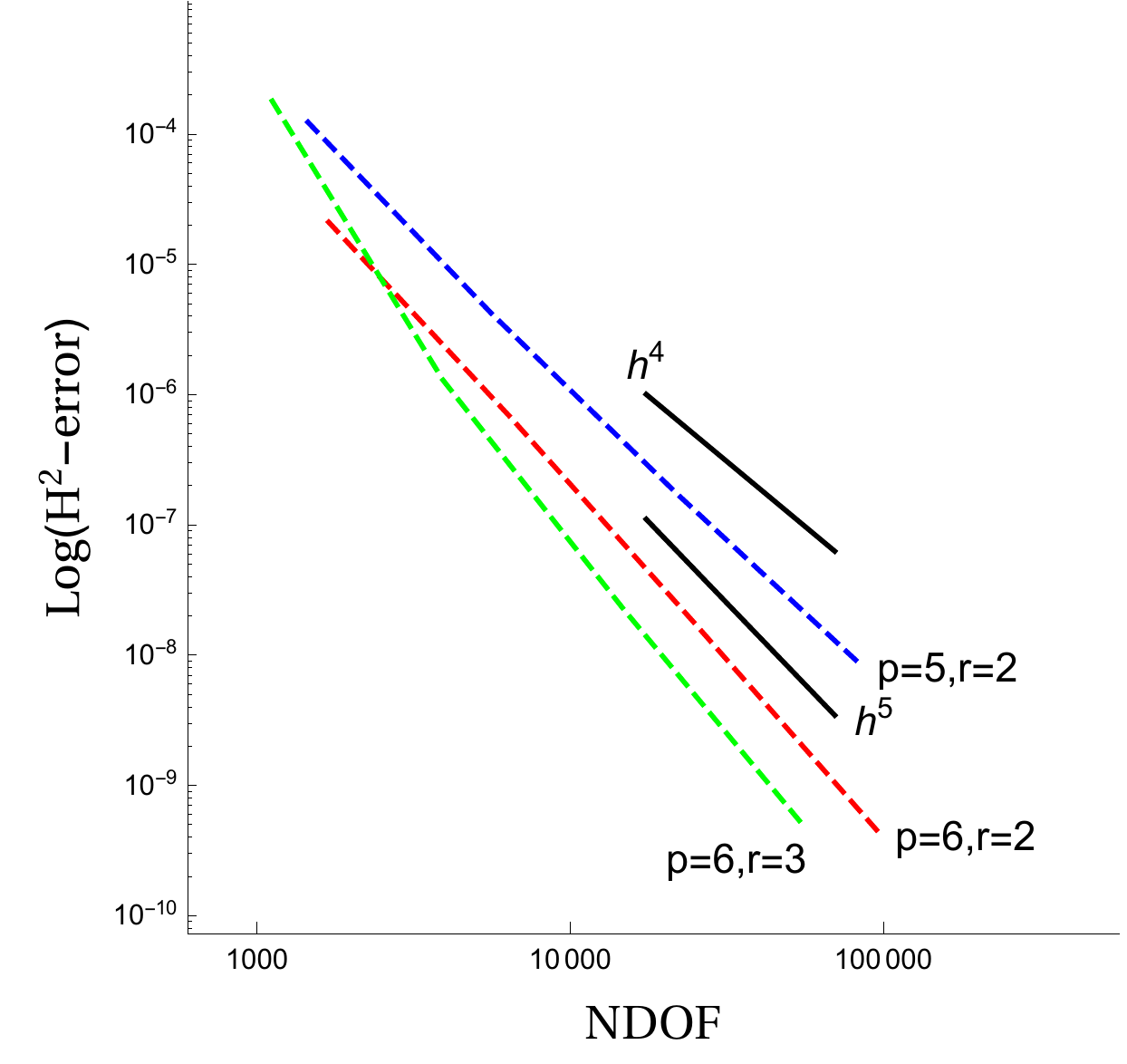} \\
\multicolumn{3}{c}{Six-patch domain~(c)} \\
\includegraphics[width=4.5cm,clip]{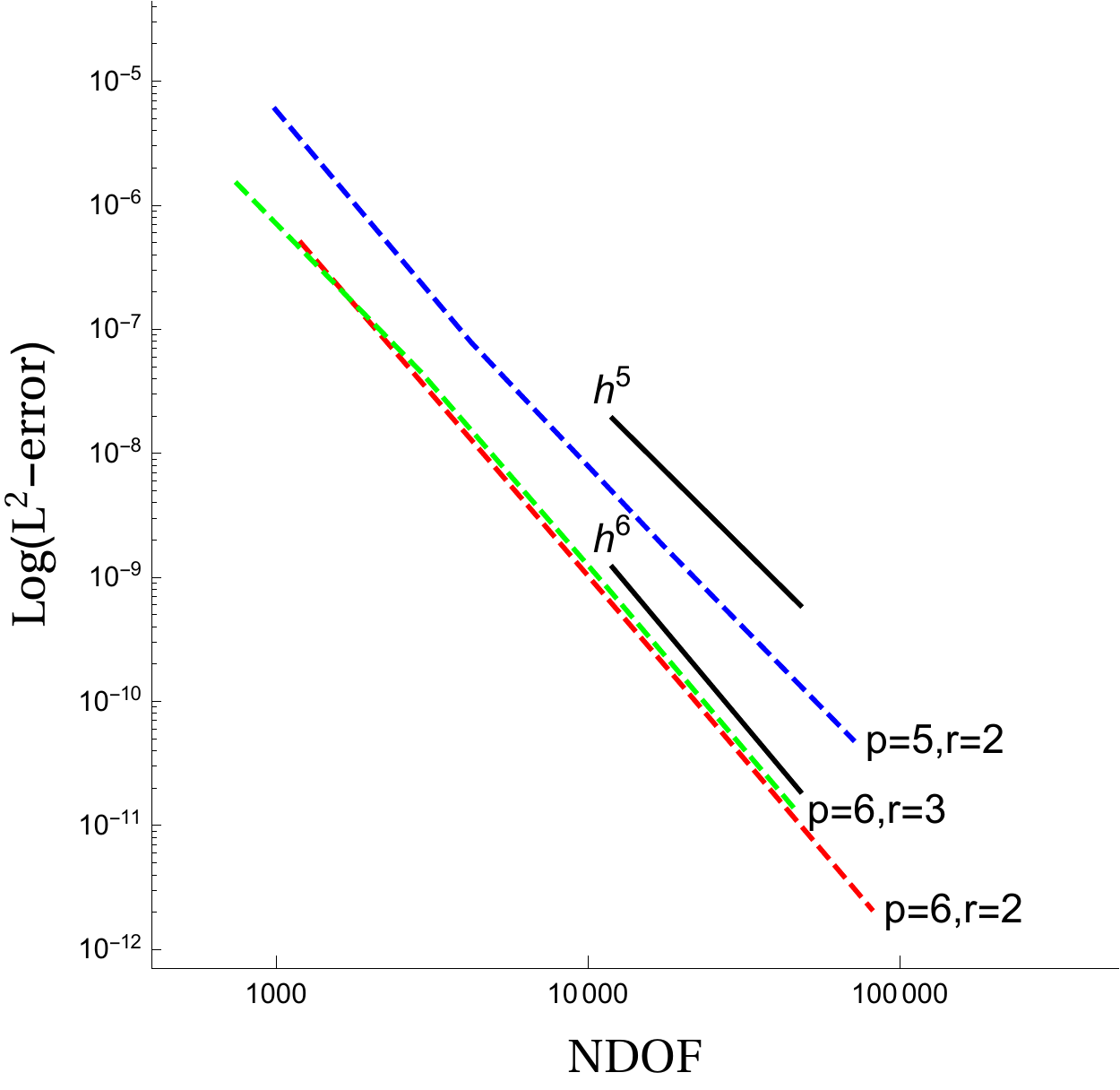} &
\includegraphics[width=4.5cm,clip]{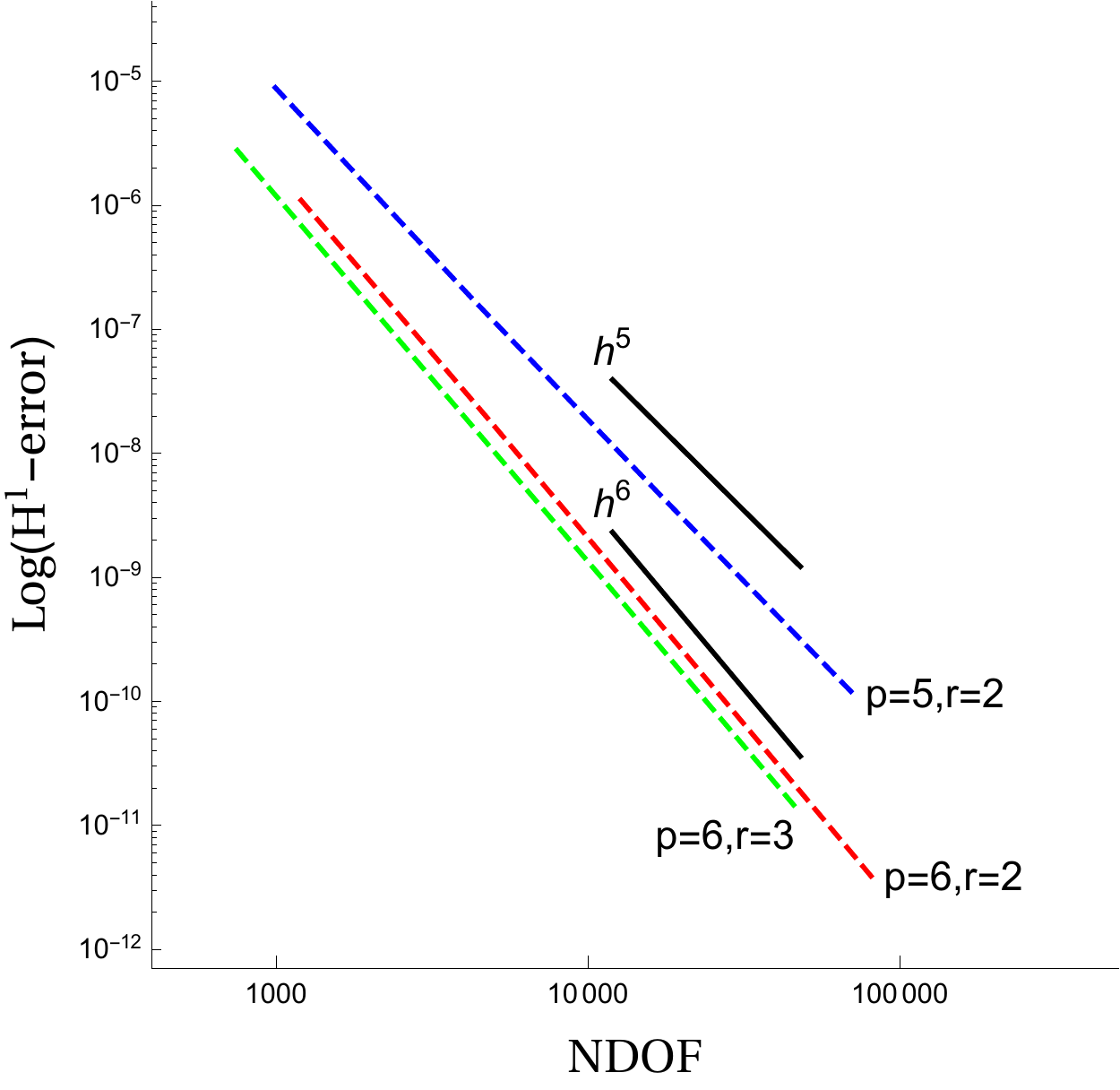} &
\includegraphics[width=4.5cm,clip]{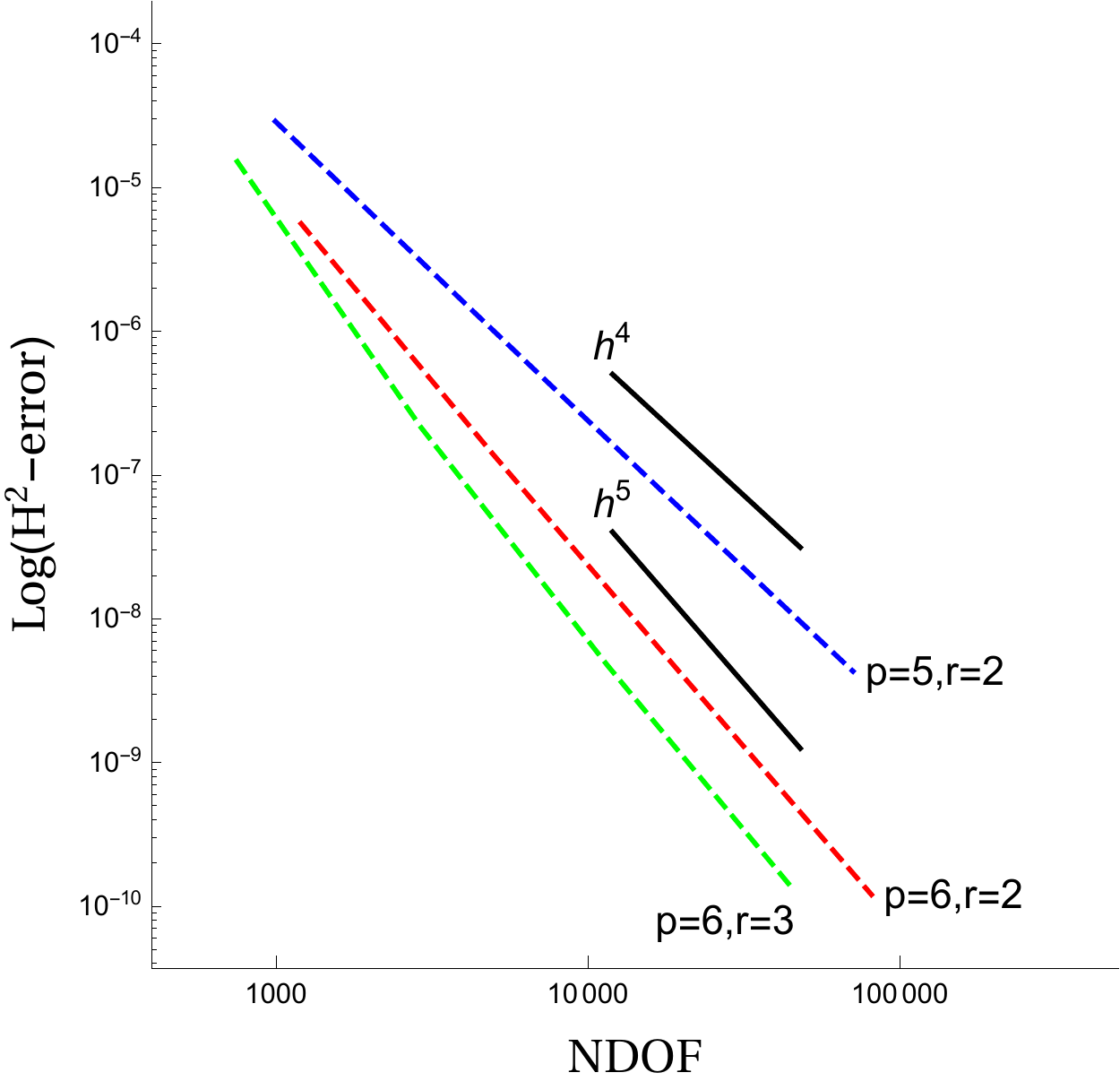} \\
\multicolumn{3}{c}{Five-patch domain~(d)} \\
\end{tabular}
\caption{Error plots w.r.t. the number of degrees of freedom (NDOF) of performing isogeometric collocation on the different multi-patch domains given in 
Fig.~\ref{fig:bilinear_domains}~(first row) for right side functions obtained by the exact solutions shown in Fig~\ref{fig:bilinear_domains}~(second row) using the 
clustered superconvergent points as collocation points.}
\label{fig:results_Cauchy}
\end{figure}
\end{ex}

\begin{ex} \label{ex:bilinearlike}
We consider the bicubic three-patch spline geometry taken from~\cite[Example~4]{KaVi19a}, 
which is also visualized in Fig.~\ref{fig:bilinearlike_domain}~(first row). Since the geometry is parameterized in such a way that it is bilinear close to the 
patch interfaces (cf. \ref{app:three_patch} for the spline control points of the geometry), our approach can be also applied to this more general multi-patch domain. 
We perform isogeometric collocation on this domain for the right side function~$f$ obtained from the exact solution
\[
 u(x_1,x_2)=-4 \cos \left(\frac{2x_1}{3}\right) \sin \left(\frac{2x_2}{3}\right),
\]
see Fig.~\ref{fig:bilinear_domains}~(first row), by employing $C^2$-smooth spaces~$\mathcal{W}_h$ with mesh sizes $h=\frac{1}{4},\frac{1}{8},\frac{1}{16},\frac{1}{32}$ for $p=5,6$ and 
$r=2$. Again, we compare the resulting relative $L^2$, $H^1$ and $H^2$ errors by using Greville points, see Fig.~\ref{fig:bilinearlike_domain}~(second row), and 
clustered superconvergent points, see Fig.~\ref{fig:bilinearlike_domain}~(third row), for collocation, and observe the same rates of convergence as in 
Example~\ref{ex:multipatchdomain}.  
\begin{figure}
\centering\footnotesize
\begin{tabular}{cc}
\includegraphics[width=6.0cm,clip]{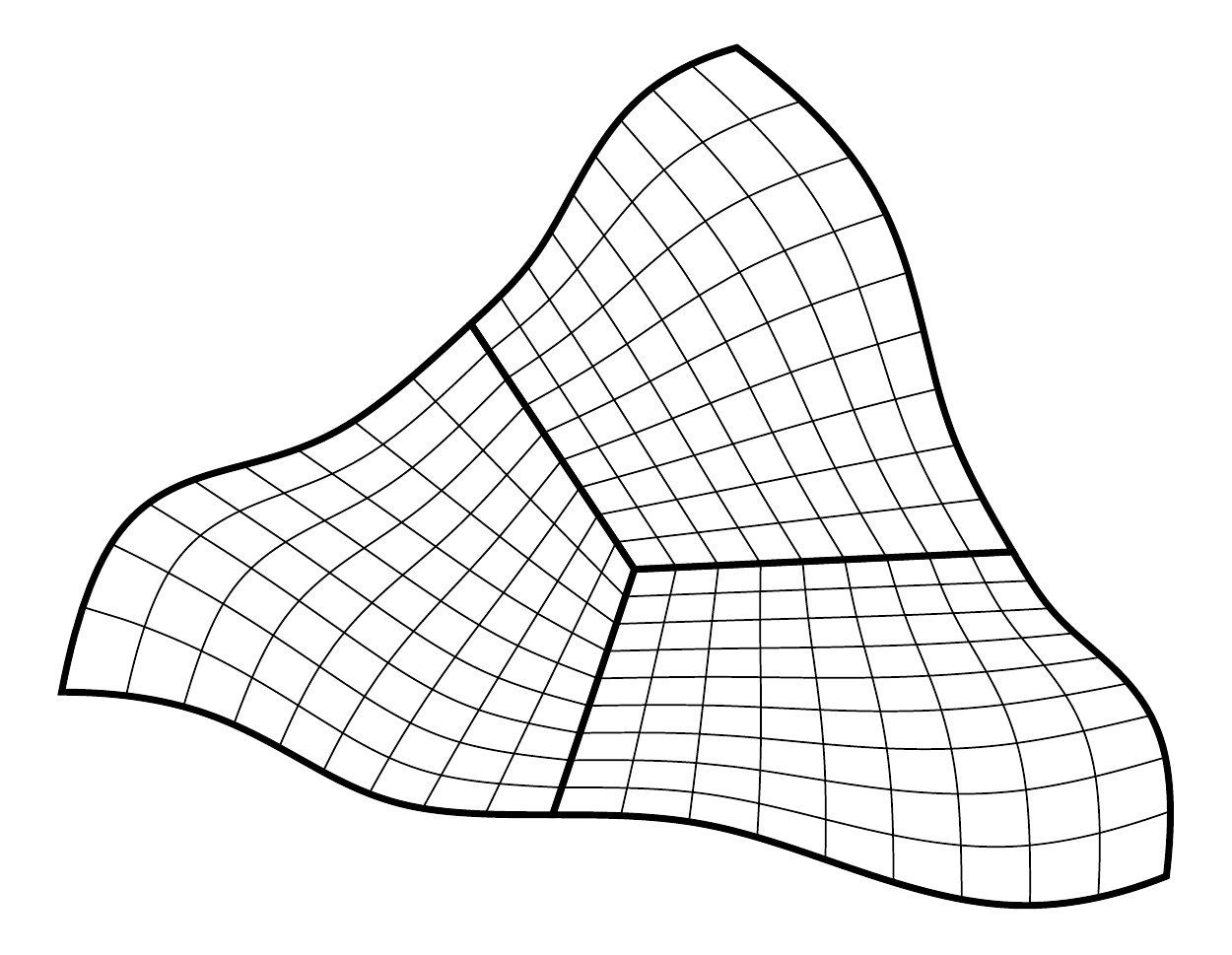} &
\includegraphics[width=5.0cm,clip]{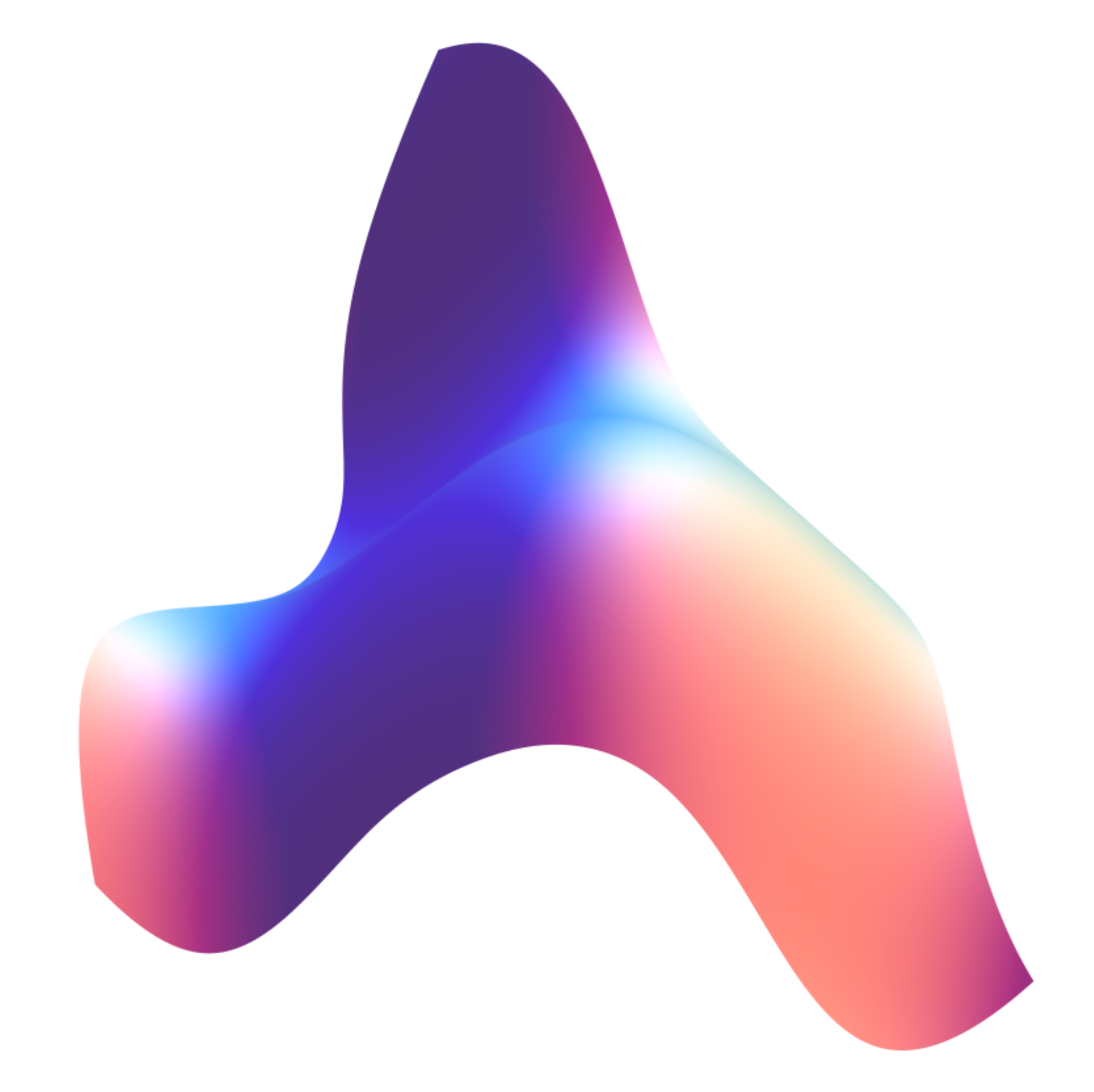} \\
Computational domain & Exact solution 
\end{tabular}
\begin{tabular}{ccc}
\includegraphics[width=4.5cm,clip]{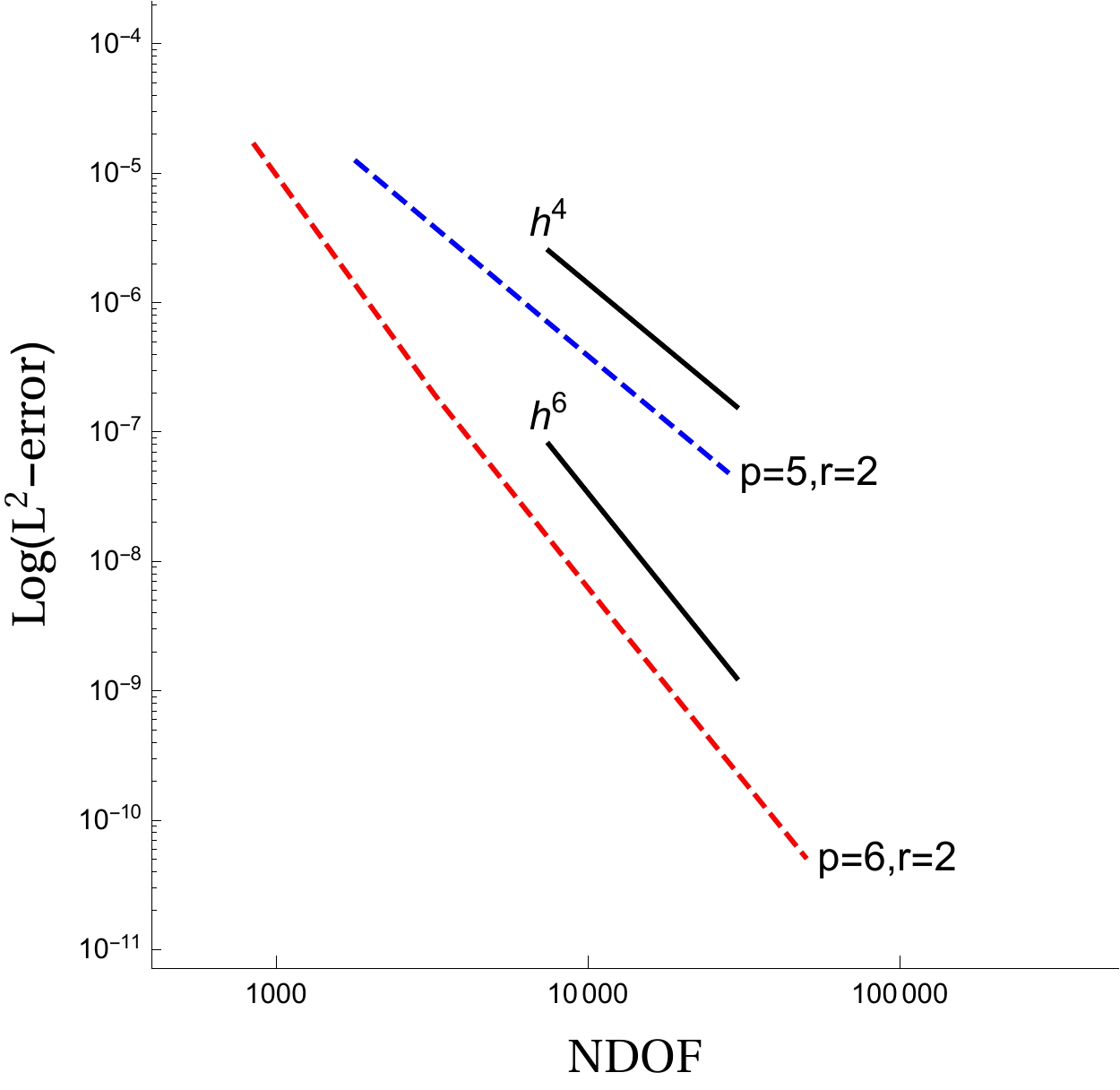} &
\includegraphics[width=4.5cm,clip]{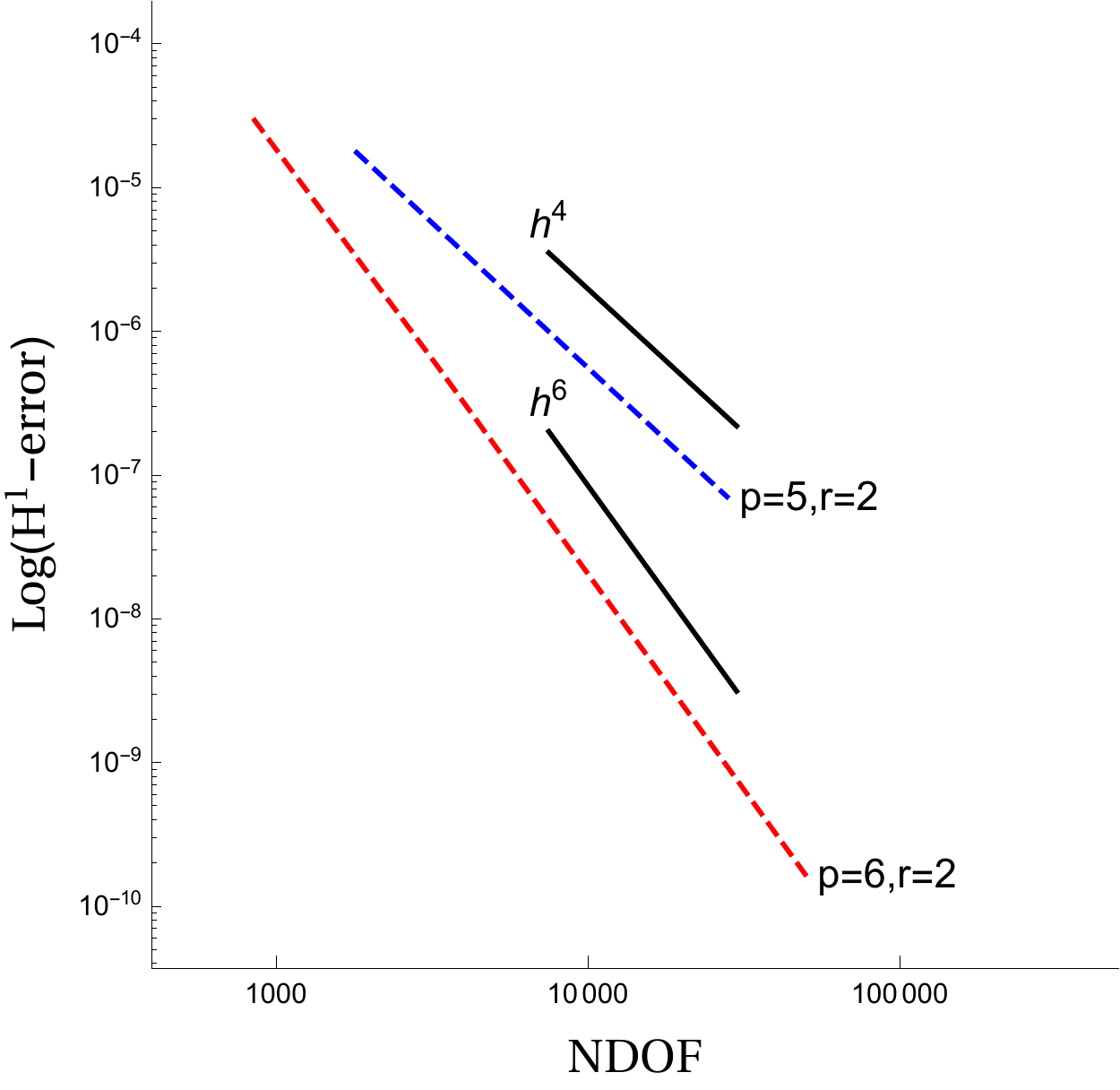} &
\includegraphics[width=4.5cm,clip]{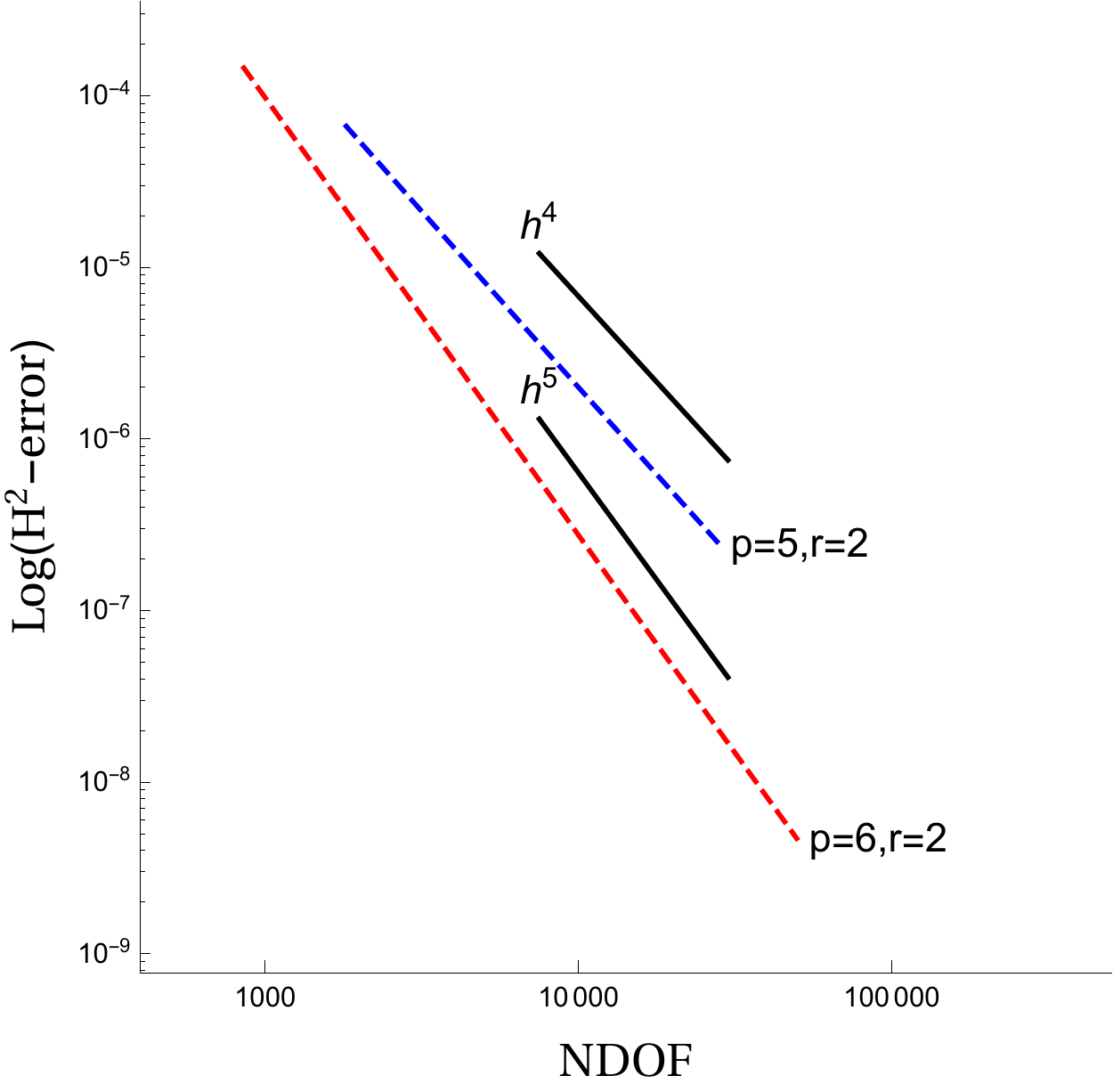}\\
\multicolumn{3}{c}{Relative $L^2$, $H^1$ and $H^2$ errors for collocation at Greville points} \\
\includegraphics[width=4.5cm,clip]{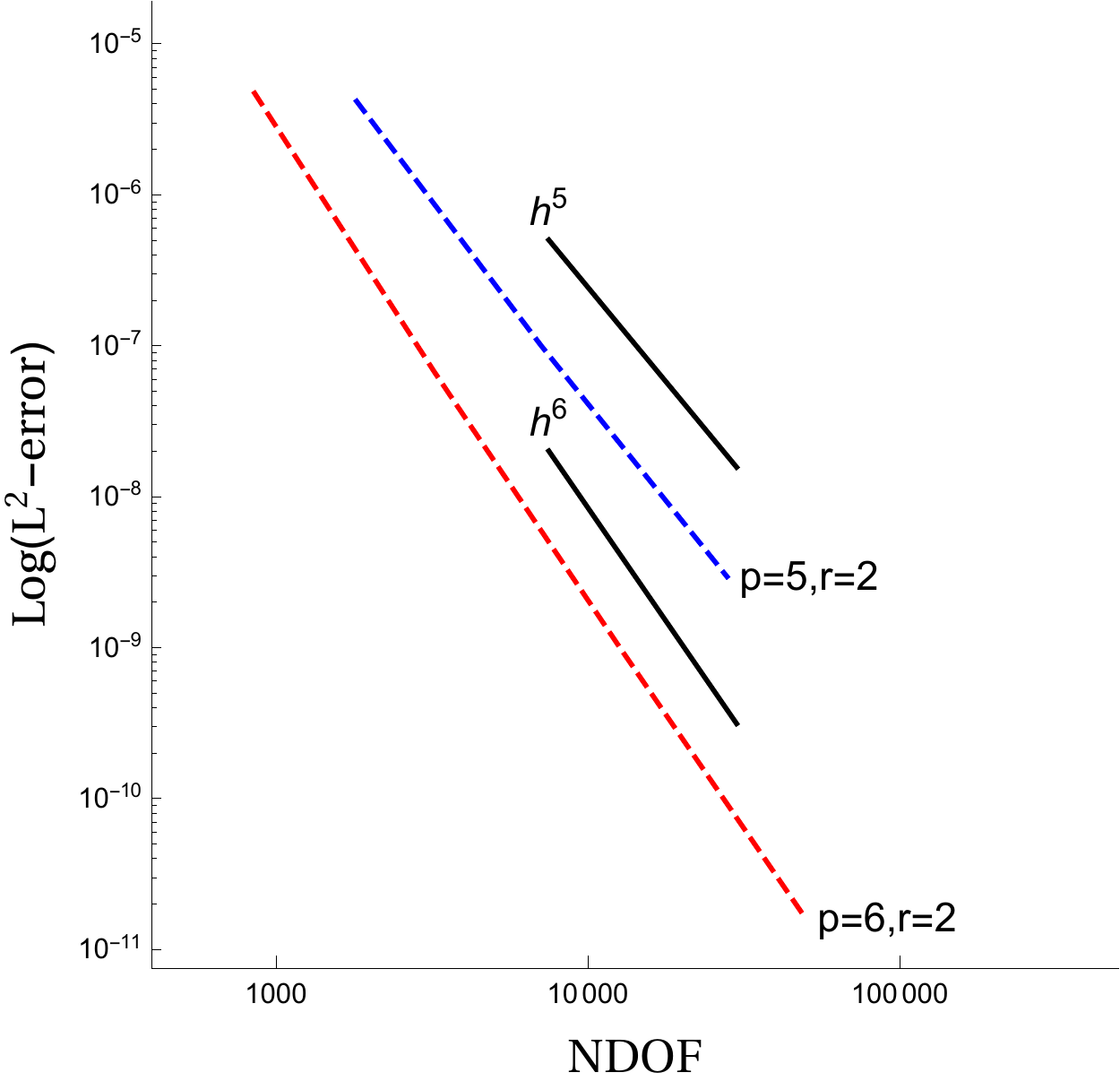} &
\includegraphics[width=4.5cm,clip]{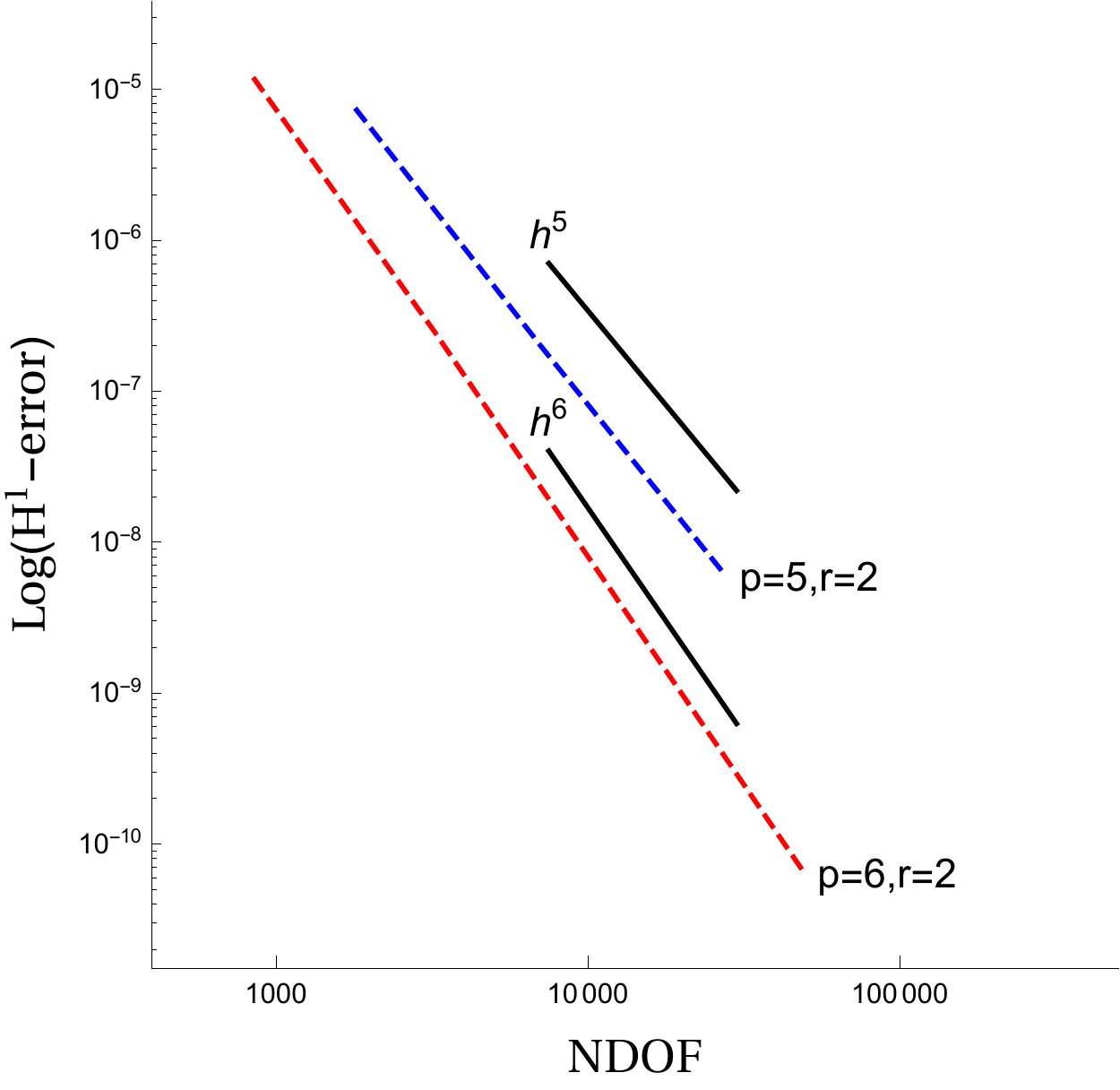} &
\includegraphics[width=4.5cm,clip]{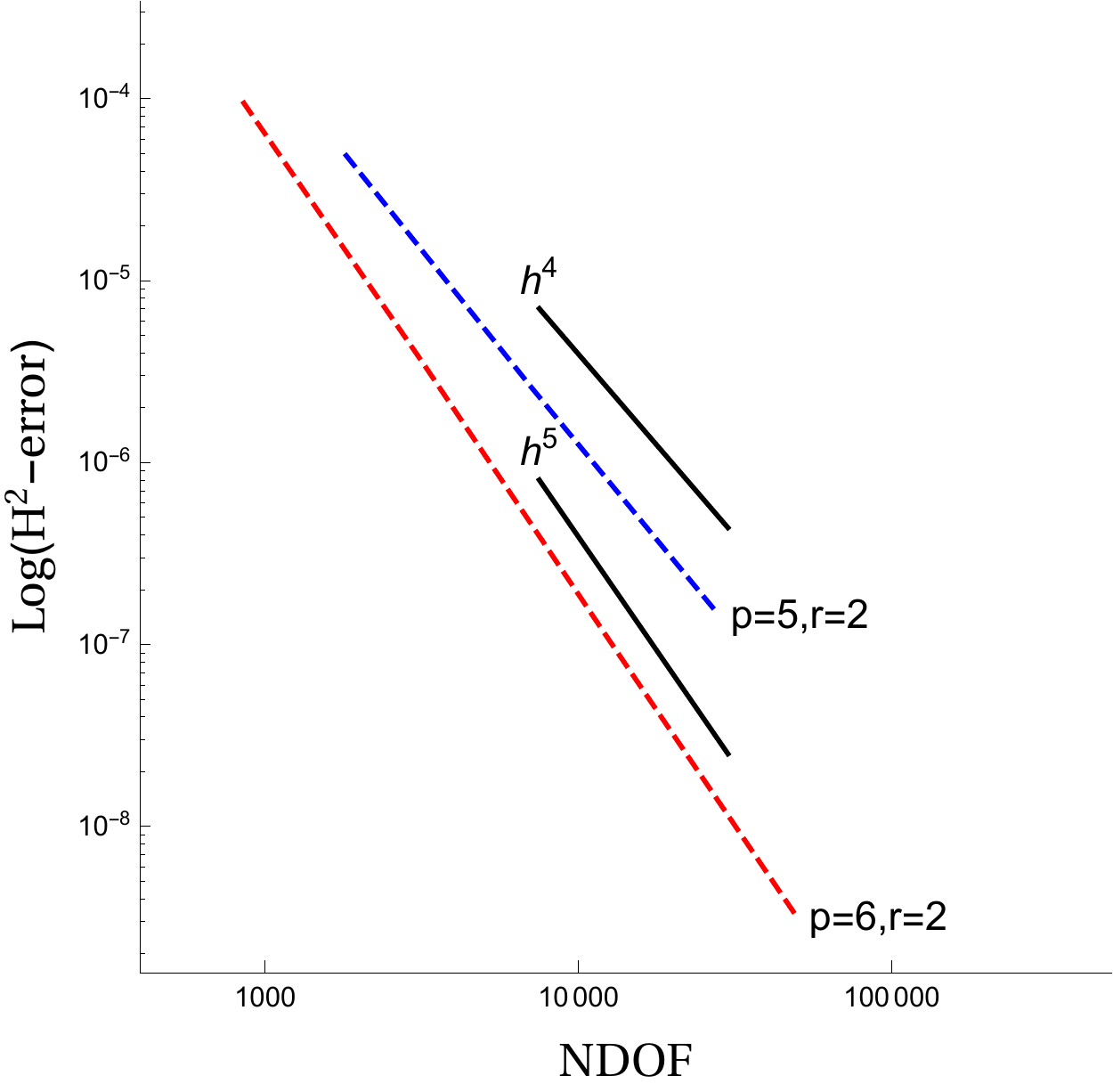}\\
\multicolumn{3}{c}{Relative $L^2$, $H^1$ and $H^2$ errors for collocation at clustered superconvergent points}
\end{tabular}
\caption{Isogeometric collocation on a bicubic three-patch spline geometry using different sets of collocation points, cf. Example.~\ref{ex:bilinearlike}.}
\label{fig:bilinearlike_domain}
\end{figure}
\end{ex}

\section{Conclusion}

We have presented a method for computing a globally $C^2$-smooth approximation of the solution of
the Poisson's equation 
over planar bilinearly parameterized multi-patch domains. Our technique is based on the 
concept of isogeometric collocation and 
on the use of a globally $C^2$-smooth isogeometric spline space as discretization space. The constructed $C^2$-smooth space is a 
particular subspace of the space of all globally $C^2$-smooth isogeometric spline functions on the considered multi-patch domain, and can be generated as the direct 
sum of simpler spaces corresponding to the individual patches, edges and vertices of the multi-patch domain. Moreover, the $C^2$-smooth space possesses a dimension 
independent of the initial geometry, and the construction of a basis with locally supported functions is simple and 
works uniformly for all possible multi-patch configurations.

Two different approaches for the choice of the collocation points have been proposed, where both strategies can be seen as the extension of successfully tested 
techniques in the one-patch case. On the one hand, we employ the tensor-product Greville points, and on the other hand, we generalize the concept of
superconvergent points (cf.~\cite{SuperConvergent2015,GomezLorenzisVariationalCollocation,MonSanTam2017}) to the case of planar multi-patch domains. While 
in the one-patch case splines of maximal smoothness, i.e. $r=p-1$, are normally used, the $C^2$-smooth multi-patch spline spaces require underlying 
spline spaces~$\mathcal{S}_{h}^{p,r}$ of regularity $2 \leq r \leq p-3$. For the sake of simplicity, we have restricted ourselves to the cases~$(p,r) \in 
\{(5,2),(6,2) ,(6,3)\}$. 

For both choices of the collocation points, the numerical results have shown the same convergence behavior with respect to $L^2$, $H^1$ and $H^2$ norm as in the 
one-patch case except in one particular case for the superconvergent points. Namely, in case of $(p,r) =(5,2)$, the convergence rates with respect to 
$L^2$ norm seems to be always reduced by one order, that is $\mathcal{O}(h^5)$ instead of the possible order of $\mathcal{O}(h^{6})$ in the one-patch case. However, 
already this slightly reduced convergence behavior in the $L^2$ norm is better as the one for the case of Greville points, where the obtained convergence rates 
are just of order~$\mathcal{O}(h^{4})$ with respect to the $L^2$ norm.

So far, our isogeometric collocation method is directly applicable only to elliptic PDEs of second order with Dirichlet boundary conditions. To impose also Neumann boundary 
conditions, further investigations amongst others in connection with the solving of the overdetermined linear system with respect to the fulfillment of the boundary 
conditions are needed, cf. \cite[Remark 2.3.2]{SuperConvergent2015}. Starting point of these studies could be e.g. the publications 
\cite{DeEvHuRe15,JiaAnitescuZhangRabczuk2019,SchBoSt15}, where different strategies are described to deal also with Neumann boundary conditions.

The paper leaves several further open issues, which are worth to study. Since the number of collocation points is slightly larger than the dimension of the $C^2$-smooth 
discretization space, the finding of a set of collocation points with the same cardinality as the dimension of the space is of interest to avoid the necessity of 
the least-squares method for solving the resulting linear system. Clearly, collocation points which lead to optimal convergence rates are of interest, too. 
However, this is also still an open problem for even spline degree in the one-patch case. A further possible topic could be the extension of our 
approach to more general multi-patch domains such as bilinear-like parameterizations, cf.~\cite{KaVi17c}, or to multi-patch structured shells and 
multi-patch volumes.

\paragraph*{\bf Acknowledgment}

The authors wish to thank the anonymous reviewers for their comments that helped to improve the paper.
V.~Vitrih was partially supported by the Slovenian Research Agency (research program P1-0404 and research projects J1-9186, J1-1715).
This support is gratefully acknowledged.

\appendix

\section{Spline control points of the geometry from Example~\ref{ex:bilinearlike}} \label{app:three_patch}

The spline geometry from Fig.~\ref{fig:bilinearlike_domain}~(first row), which was taken from \cite[Example~4]{KaVi19a} and which is also used in Example~\ref{ex:bilinearlike} 
in this paper, consists of three geometry mappings~$\ab{F}^{(i)} \in  \mathcal{S}_{h}^{\ab{p},\ab{r}}([0,1]^2) \times \mathcal{S}^{\ab{p},\ab{r}}_{h}([0,1]^2) $, 
$i=1,2,3$, with $\ab{p}=(3,3)$, $r=(2,2)$ and $h=1/4$. Therefore, each geometry mapping~$\ab{F}^{(i)}$, $i=1,2,3$, possesses a spline parameterization of the form
\[
 \ab{F}^{(i)} (\xi_1,\xi_2) = \sum_{j_1=0}^{6} \sum_{j_2=0}^{6} \ab{c}^{(i)}_{j_1,j_2} N^{\ab{p},\ab{r}}_{j_1,j_2}(\xi_1,\xi_2),
\]
where the single control points~$\ab{c}^{(i)}_{j_1,j_2}$ are given in Table~\ref{tab:three-patch}.
\begin{table}[ht]
\footnotesize
\centering
\setlength{\tabcolsep}{0.1em}
\begin{tabular}{|ccccccc|}
\hline
\multicolumn{7}{|c|}{$\f{c}_{j_1,j_2}^{(1)}$} \\
\hline
& & & & & & \\[-0.35cm]
$(\frac{17}{3}, 2)$ & $(\frac{853}{144}, \frac{169}{84})$ & $(\frac{103}{16}, \frac{57}{28})$ & $(\frac{173}{24}, \frac{29}{14})$ & 
$(\frac{383}{48}, \frac{59}{28})$ & $(\frac{1223}{144}, \frac{179}{84})$ & $(\frac{35}{4}, \frac{15}{7})$\\[0.1cm]
$(\frac{101}{18}, \frac{11}{6})$ & $(\frac{10163}{1728}, \frac{1859}{1008})$ & $(\frac{411}{64}, \frac{209}{112})$ & $(\frac{2083}{288}, \frac{319}{168})$ & 
$(\frac{4633}{576}, \frac{649}{336})$ & $(\frac{14833}{1728}, \frac{1969}{1008})$ & $(\frac{425}{48}, \frac{55}{28})$ \\[0.1cm]
$(\frac{11}{2}, \frac{3}{2})$ & $(\frac{371}{64}, \frac{169}{112})$ & $(\frac{409}{64}, \frac{171}{112})$ & $(\frac{233}{32}, \frac{87}{56})$ & 
$(\frac{523}{64}, \frac{177}{112})$ & $(\frac{561}{64}, \frac{179}{112})$ & $(\frac{145}{16}, \frac{45}{28})$ \\[0.1cm]
$(\frac{16}{3}, 1)$ & $(\frac{1633}{288}, \frac{169}{168})$ & $(\frac{203}{32}, \frac{57}{56})$ & $(\frac{37}{5}, \frac{99}{100})$ & $(\frac{849}{100}, \frac{103}{100})$ & 
$(\frac{923}{100}, \frac{57}{50})$ & $(\frac{48}{5}, \frac{121}{100})$ \\[0.1cm]
$(\frac{31}{6}, \frac{1}{2})$ & $(\frac{3193}{576}, \frac{169}{336})$ & $(\frac{403}{64}, \frac{57}{112})$ & $(\frac{749}{100}, \frac{8}{25})$ & 
$(\frac{873}{100}, \frac{29}{100})$ & $(\frac{961}{100}, \frac{27}{50})$ & $(\frac{251}{25}, \frac{13}{20})$ \\[0.1cm]
$(\frac{91}{18}, \frac{1}{6})$ & $(\frac{9433}{1728}, \frac{169}{1008})$ & $(\frac{401}{64}, \frac{19}{112})$ & $(\frac{15}{2}, -\frac{21}{100})$ & 
$(\frac{437}{50}, -\frac{12}{25})$ & $(\frac{961}{100}, -\frac{6}{25})$ & $(\frac{201}{20}, -\frac{1}{14})$ \\[0.1cm]
$(5, 0)$ & $(\frac{65}{12}, 0)$ & $(\frac{25}{4}, 0)$ & $(\frac{15}{2}, -\frac{12}{25})$ & $(\frac{35}{4}, -\frac{21}{25})$ & $(\frac{115}{12}, -\frac{2}{3})$ & 
$(10, -\frac{1}{2})$ \\[0.1cm]
\hline
\hline 
\multicolumn{7}{|c|}{$\f{c}_{j_1,j_2}^{(2)}$} \\
\hline
& & & & & & \\[-0.35cm]
$(\frac{17}{3}, 2)$ & $(\frac{50}{9}, \frac{13}{6})$ & $(\frac{16}{3}, \frac{5}{2})$ & $(5, 3)$ & $(\frac{14}{3}, \frac{7}{2})$ & $(\frac{40}{9}, \frac{23}{6})$ & 
$(\frac{13}{3}, 4)$ \\[0.1cm]
$(\frac{853}{144}, \frac{169}{84})$ & $(\frac{10033}{1728}, \frac{2209}{1008})$ & $(\frac{3209}{576}, \frac{857}{336})$ & $(\frac{167}{32}, \frac{173}{56})$ & 
$(\frac{2803}{576}, \frac{1219}{336})$ & $(\frac{8003}{1728}, \frac{4019}{1008})$ & $(\frac{325}{72}, \frac{25}{6})$ \\[0.1cm]
$(\frac{103}{16}, \frac{57}{28})$ & $(\frac{1211}{192}, \frac{251}{112})$ & $(\frac{387}{64}, \frac{297}{112})$ & $(\frac{181}{32}, \frac{183}{56})$ & 
$(\frac{337}{64}, \frac{435}{112})$ & $(\frac{961}{192}, \frac{481}{112})$ & $(\frac{39}{8}, \frac{9}{2})$ \\[0.1cm]
$(\frac{173}{24}, \frac{29}{14})$ & $(\frac{2033}{288}, \frac{389}{168})$ & $(\frac{649}{96}, \frac{157}{56})$ & $(\frac{127}{20}, \frac{359}{100})$ & 
$(\frac{587}{100}, \frac{441}{100})$ & $(\frac{109}{20}, \frac{99}{20})$ & $(\frac{523}{100}, \frac{523}{100})$ \\[0.1cm]
$(\frac{383}{48}, \frac{59}{28})$ & $(\frac{4499}{576}, \frac{803}{336})$ & $(\frac{1435}{192}, \frac{331}{112})$ & $(\frac{178}{25}, \frac{391}{100})$ & 
$(\frac{661}{100}, \frac{49}{10})$ & $(\frac{597}{100}, \frac{279}{50})$ & $(\frac{142}{25}, \frac{591}{100})$ \\[0.1cm]
$(\frac{1223}{144}, \frac{179}{84})$ & $(\frac{14363}{1728}, \frac{2459}{1008})$ & $(\frac{4579}{576}, \frac{1027}{336})$ & $(\frac{769}{100}, \frac{409}{100})$ & 
$(\frac{73}{10}, \frac{257}{50})$ & $(\frac{33}{5}, \frac{583}{100})$ & $(\frac{31}{5}, \frac{37}{6})$ \\[0.1cm]
$(\frac{35}{4}, \frac{15}{7})$ & $(\frac{137}{16}, \frac{69}{28})$ & $(\frac{131}{16}, \frac{87}{28})$ & $(\frac{799}{100}, \frac{209}{50})$ 
& $(\frac{763}{100}, \frac{263}{50})$ & $(\frac{111}{16}, \frac{83}{14})$ & $(\frac{13}{2}, \frac{25}{4})$ \\[0.1cm]
\hline
\multicolumn{7}{|c|}{$\f{c}_{j_1,j_2}^{(3)}$} \\
\hline
& & & & & & \\[-0.35cm]
$(\frac{17}{3}, 2)$ & $(\frac{101}{18}, \frac{11}{6})$ & $(\frac{11}{2}, \frac{3}{2})$ & $(\frac{16}{3}, 1)$ & $(\frac{31}{6}, \frac{1}{2})$ & 
$(\frac{91}{18}, \frac{1}{6})$ & $(5, 0)$ \\[0.1cm]
$(\frac{50}{9}, \frac{13}{6})$ & $(\frac{2365}{432}, \frac{143}{72})$ & $(\frac{85}{16}, \frac{13}{8})$ & $(\frac{365}{72}, \frac{13}{12})$ & 
$(\frac{695}{144}, \frac{13}{24})$ & $(\frac{2015}{432}, \frac{13}{72})$ & $(\frac{55}{12}, 0)$ \\[0.1cm]
$(\frac{16}{3}, \frac{5}{2})$ & $(\frac{749}{144}, \frac{55}{24})$ & $(\frac{79}{16}, \frac{15}{8})$ & $(\frac{109}{24}, \frac{5}{4})$ & 
$(\frac{199}{48}, \frac{5}{8})$ & $(\frac{559}{144}, \frac{5}{24})$ & $(\frac{15}{4}, 0)$ \\[0.1cm]
$(5, 3)$ & $(\frac{115}{24}, \frac{11}{4})$ & $(\frac{35}{8}, \frac{9}{4})$ & $(\frac{19}{5}, \frac{43}{25})$ & $(\frac{83}{25}, \frac{117}{100})$ & 
$(\frac{153}{50}, \frac{7}{10})$ & $(\frac{59}{20}, \frac{9}{20})$ \\[0.1cm]
$(\frac{14}{3}, \frac{7}{2})$ & $(\frac{631}{144}, \frac{77}{24})$ & $(\frac{61}{16}, \frac{21}{8})$ & $(\frac{61}{20}, \frac{227}{100})$ & 
$(\frac{123}{50}, \frac{181}{100})$ & $(\frac{113}{50}, \frac{6}{5})$ & $(\frac{43}{20}, \frac{9}{10})$ \\[0.1cm]
$(\frac{40}{9}, \frac{23}{6})$ & $(\frac{1775}{432}, \frac{253}{72})$ & $(\frac{55}{16}, \frac{23}{8})$ & $(\frac{251}{100}, \frac{263}{100})$ &
$(\frac{7}{4}, \frac{113}{50})$ & $(\frac{151}{100}, \frac{149}{100})$ & $(\frac{17}{12}, 1)$ \\[0.1cm]
$(\frac{13}{3}, 4)$ & $(\frac{143}{36}, \frac{11}{3})$ & $(\frac{13}{4}, 3)$ & $(\frac{56}{25}, \frac{141}{50})$ & $(\frac{141}{100}, \frac{247}{100})$ & 
$(\frac{10}{9}, \frac{19}{12})$ & $(1, 1)$ \\[0.1cm]
\hline
\end{tabular}
\caption{Control points~$\ab{c}_{j_1,j_2}^{(i)}$, $i=1,2,3$, of the bicubic three-patch spline geometry shown in Fig.~\ref{fig:bilinearlike_domain}~(first row).}
\label{tab:three-patch} 
\end{table}

\end{document}